\documentclass[11pt]{amsart}
\usepackage{fullpage}
\usepackage{amsmath, amssymb, amsthm}
\usepackage{enumitem}
\usepackage[foot]{amsaddr}

\usepackage{tikz}
\usetikzlibrary{calc}
\usepackage{xcolor}
\usepackage[colorlinks=true,linkcolor=red,citecolor=blue,urlcolor=cyan]{hyperref}

\newtheorem{theorem}{Theorem}[section]
\newtheorem{proposition}[theorem]{Proposition}
\newtheorem{lemma}[theorem]{Lemma}

\theoremstyle{remark}
\newtheorem{remark}[theorem]{Remark}

\newcommand{\R}{\mathbb{R}}
\newcommand{\Z}{\mathbb{Z}}
\newcommand{\Hess}{\operatorname{Hess}}
\newcommand{\tr}{\operatorname{tr}}
\renewcommand{\v}{\mathbf{v}}
\newcommand{\x}{\mathbf{x}}
\renewcommand{\P}{\mathcal{P}}
\newcommand{\bzero}{\mathbf{0}}
\newcommand{\bth}{\boldsymbol{\theta}}
\newcommand{\bsigma}{\boldsymbol{\sigma}}
\newcommand{\bdelta}{\boldsymbol{\delta}}
\newcommand{\bepsilon}{\boldsymbol{\varepsilon}}
\newcommand{\btau}{\boldsymbol{\tau}}
\newcommand{\sm}{\mathfrak{sim}}

\newcommand{\calQ}{\mathcal{Q}}
\newcommand{\tP}{\widetilde{\P}}
\newcommand{\tQ}{\widetilde{\calQ}}
\newcommand{\TP}{T\P}
\newcommand{\TQ}{T\calQ}
\newcommand{\calS}{\mathcal{S}}
\newcommand{\cube}{\square_d}
\newcommand{\blam}{\bar\lambda}
\newcommand{\bmu}{\bar\mu}

\newcommand{\bfn}{\mathbf{n}}
\newcommand{\bfm}{\mathbf{m}}
\newcommand{\PiD}{\Pi^{\lambda}}
\newcommand{\SiD}{\Sigma^{\lambda}}
\newcommand{\PiN}{\Pi^{\mu}}
\newcommand{\SiN}{\Sigma^{\mu}}
\newcommand{\ES}{E_S}
\newcommand{\ETU}{E_{TU}}

\usepackage[backend=biber,giveninits=true,sorting=nyt,style=alphabetic,natbib=true,maxcitenames=2,maxbibnames=10,url=false,doi=true,backref=false]{biblatex}
\renewbibmacro{in:}{\ifentrytype{article}{}{\printtext{\bibstring{in}\intitlepunct}}}

\begin{document}
\title{Local optimality of the $d$-cube \\ for Laplacian eigenvalues}

\date{\today}

\author{Braxton Osting}
\email{osting@math.utah.edu}
\address{Department of Mathematics, University of Utah, Salt Lake City, UT}

\author{Robert Viator}
\email{viatorr@denison.edu}
\address{Department of Mathematics, Denison University, Granville, OH}

\thanks{B.O. acknowledges partial support from NSF DMS-2136198 and DMS-2513175.}


\begin{abstract}
For each $d \ge 2$ we prove that the unit $d$-cube $\cube \subset \R^d$ is, on the class $\P$ of convex polytopes combinatorially equivalent to a $d$-cube, (i) a strict local minimizer of the scale-invariant isoperimetric ratio $I(\Omega) = \tfrac{P(\Omega)^d}{|\Omega|^{d-1}}$, with $P$ the perimeter, (ii) a strict local minimizer of the volume-normalized first Dirichlet eigenvalue $\blam(\Omega) = |\Omega|^{2/d}\lambda_1(\Omega)$, and (iii) a strict local maximizer of the volume-normalized first Neumann eigenvalue $\bmu(\Omega) = |\Omega|^{2/d}\mu_1(\Omega)$. The method is a shape-derivative Hessian analysis adapted to the hyperoctahedral symmetry of the cube.  In symmetry-adapted coordinates on the parameter space --- the decomposition into irreducible representations of the hyperoctahedral group --- each Hessian is block-diagonal, with four scalar eigenvalues, which we give in closed form in $\pi$, rational coefficients, and a few universal series-constants. The first Neumann eigenvalue is not differentiable at $\cube$, since its lowest cluster splits under perturbation. The Neumann case therefore works with a smooth auxiliary functional, whose Hessian is computed by second-order perturbation theory of the cluster; the coupling to modes outside the cluster enters at the same order. Finally, we show that extremality of $I$, $\blam$, and $\bmu$ does not hold in high dimensions for a particular multilinear perturbation that leaves $\P$ by breaking convexity. All results are local; we discuss why a global statement, available for quadrilaterals (Endo and Osting, 2026), appears out of reach for $d \ge 3$.
\end{abstract}

\maketitle

%%%%%%%%%%%%%%%%%%%%%%%%%%%%%%%%%%
%%%%%%%%%%%%%%%%%%%%%%%%%%%%%%%%%%
\section{Introduction}\label{s:Intro}

The problem of identifying the domain that minimizes the first Dirichlet--Laplace eigenvalue under a volume constraint traces back to Rayleigh~\cite{Rayleigh1894}, who conjectured that the disk minimizes the principal eigenvalue among planar domains of fixed area. This was settled in dimension $d = 2$ by Faber~\cite{Faber1923} and Krahn~\cite{Krahn1925}, and extended by Krahn~\cite{Krahn1926} to all $d \ge 2$: the ball is the unique minimizer of $\blam(\Omega) := |\Omega|^{2/d}\lambda_1(\Omega)$ among bounded open sets in $\R^d$. 
The analogous problem for the first nonzero Neumann eigenvalue was resolved by Szeg\H{o}~\cite{Szego1954} in the plane and Weinberger~\cite{Weinberger1956} in all $d \ge 2$: the ball is the unique maximizer of $\bmu(\Omega) := |\Omega|^{2/d}\mu_1(\Omega)$. 
Without a convexity constraint, Bucur and Henrot~\cite{BucurHenrot2019} proved that the second nontrivial Neumann eigenvalue $\mu_2$ is maximized by a union of two equal balls, and Bogosel, Henrot, and Michetti~\cite{BogoselHenrotMichetti2024} analyzed the convexity-constrained problem. In the planar convexity-constrained setting, optimal shapes are frequently polygonal (Lamboley--Novruzi~\cite{LamboleyNovruzi2010}), which motivates carrying out shape optimization directly within classes of polytopes.
For broad surveys of spectral shape optimization, see~\cite{Henrot2006, HenrotSurvey2017, HenrotPierre2018}; for a general introduction to spectral geometry, see~\cite{LevitinMangoubiPolterovich2023}.

Both Faber--Krahn and Szeg\H{o}--Weinberger rely on Schwarz symmetrization and the rotational symmetry of the ball. These tools are unavailable when the admissible class is restricted to a fixed combinatorial type that cannot be symmetrized to a ball. P\'olya conjectured~\cite{PolyaSzego1951} that the regular $n$-gon minimizes $\blam$ among convex $n$-gons. The cases $n = 3, 4$ follow from Steiner symmetrization~\cite{Henrot2006}; for $n \ge 5$ the conjecture remains open in general. Bogosel and Bucur~\cite{BogoselBucur2024} have recently shown that it reduces, for each $n$, to finitely many certified numerical computations, which they carry out for $n = 5, 6, 7, 8$ to establish local minimality; a fully validated, interval-arithmetic certification for $n = 5, 6$ appears in~\cite{BogoselBucur2024b}. 
Two related functionals on convex $n$-gons are optimized by the regular $n$-gon for every $n$: the logarithmic capacity, by Solynin and Zalgaller~\cite{SolyninZalgaller2004}, and the Cheeger constant, by Bucur and Fragal\`a~\cite{BucurFragala2016}.

The natural Neumann analogue --- that the regular $n$-gon maximizes $\bmu$ among convex $n$-gons --- is known for $n = 3$ by the transplantation argument of Laugesen and Siudeja~\cite{LaugesenSiudeja2009}, and for $n = 4$ by~\cite{OstingQuadNeumann2026}. Neither conjecture is granted \emph{a priori}: counter-intuitive symmetry-breaking can occur in related shape problems, as in Graham's largest small hexagon~\cite{Graham1975}, where the regular hexagon is \emph{not} the area maximizer among hexagons of fixed diameter, and, closer to the present setting, in the $p \to 1^+$ limit of the polygonal Szeg\H{o}--Weinberger inequality for the $p$-Laplacian (P\'olya's longest shortest-fence problem), where Bogosel, Bucur, and Fragal\`a~\cite{BogoselBucurFragala2026} show that the optimal quadrilateral is an isosceles trapezium rather than the square.

In dimensions $d \ge 3$ the analogous extremal problems for polytopes are far less explored than the planar ones. Among axis-aligned boxes, where the eigenvalues are explicit, the optimization of low and asymptotic eigenvalues is already non-trivial: see 
Antunes and Freitas~\cite{AntunesFreitas2013} for optimal Dirichlet rectangles, 
van den Berg, Bucur, and Gittins~\cite{vandenBergBucurGittins2016} for Neumann rectangles, and 
van den Berg--Gittins~\cite{vandenBergGittins2017} and Gittins--Larson~\cite{GittinsLarson2017} for cuboids. 
Closest in spirit to the present work is the polytope shape-optimization framework of D\'iaz-Avalos and Laurain~\cite{DiazAvalosLaurain2025}. They parameterize simple polytopes in arbitrary dimension by their supporting hyperplanes, establish stability under small hyperplane perturbations, and compute first- and second-order shape derivatives; their numerics support regular polygons in $2$D and regular polyhedra in $3$D as optima for the polyhedral Saint--Venant functional.  

In this paper we consider the higher-dimensional analogue of these questions for convex polytopes combinatorially equivalent to $d$-cubes. Let $\cube := [-\tfrac12, \tfrac12]^d \subset \R^d$ denote the unit $d$-cube centered at the origin, the natural higher-dimensional counterpart of the regular $n$-gon. We parameterize via the bounding hyperplanes: $\cube$ is the intersection of $2d$ half-spaces, $\cube = \bigcap_{k=1}^{d}\bigl\{\x \in \R^d : -\tfrac{1}{2} \le x_k \le \tfrac{1}{2}\bigr\}$. A small \emph{convex perturbation} is obtained by slightly tilting each hyperplane. The resulting intersection is automatically convex, and for sufficiently small tilts is a polytope with $2d$ flat facets meeting at $2^d$ vertices. We denote this class of \emph{combinatorial $d$-cubes near $\cube$} by $\tP$. The functionals below are invariant under similarity transformations, so we may normalize each body by a similarity; $\P$ denotes the normalized class, made precise in Section~\ref{s:CubeParam}, and every body in $\tP$ near $\cube$ is similar to one in $\P$.

We consider three shape functionals on $\P$:
\begin{equation}\label{e:three-functionals-intro}
I(\Omega) := \frac{P(\Omega)^d}{|\Omega|^{d-1}}, \qquad \blam(\Omega) := |\Omega|^{2/d}\lambda_1(\Omega), \qquad \bmu(\Omega) := |\Omega|^{2/d}\mu_1(\Omega),
\end{equation}
the scale-invariant isoperimetric ratio (with $P(\Omega)$ the surface area), and the volume-normalized fundamental Dirichlet and Neumann eigenvalues.  Our main result is:

\begin{theorem}\label{t:Main}
For all $d \ge 2$, the cube is a critical point of $I$ and $\blam$ on $\P$, by the symmetry of the cube.  Moreover, the $d$-cube is
\begin{enumerate}
\item[(a)] a strict local minimizer of $I$ on $\P$, with $I(\cube) = (2d)^d$;
\item[(b)] a strict local minimizer of $\blam$ on $\P$, with $\blam(\cube) = d\pi^2$;
\item[(c)] a strict local maximizer of $\bmu$ on $\P$, with $\bmu(\cube) = \pi^2$.
\end{enumerate}
\end{theorem}

\noindent The functional $\bmu$ is not differentiable at $\cube$, but every one-sided directional derivative of $\bmu$ at $\cube$ is non-positive (Remark~\ref{r:one-sided}). 
Strictness in Theorem~\ref{t:Main} is understood modulo similarity: every $\Omega \in \P$ near $\cube$ that is not similar to $\cube$ satisfies the strict inequality.

Part (a) of Theorem~\ref{t:Main} is a purely geometric statement, framed by two classical results. For $d = 3$, Fejes T\'oth~\cite{FejesToth1948} proved the stronger, global statement that the cube minimizes the isoperimetric ratio among all convex polyhedra with six faces. In every dimension, Lindel\"of's theorem (see~\cite[\S 18.3]{Gruber2007}) solves the isoperimetric problem when the facet normals are held fixed: the minimizer is circumscribed about a ball, so the cube is optimal among axis-aligned boxes. The class $\P$ allows the normals to tilt, and part (a) provides the local statement on all of $\P$ in every dimension. Extremum problems of this type for polytopes are surveyed in~\cite{Florian1993} and~\cite[\S\S 11.2, 18.3]{Gruber2007}. Part (b) extends to arbitrary dimension the planar Steiner-symmetrization result that the square strictly minimizes $\blam$ among convex quadrilaterals~\cite{PolyaSzego1951}, complementing the polygonal Faber--Krahn analysis of Bogosel--Bucur~\cite{BogoselBucur2024}. Part (c) extends the planar local-maximality result of~\cite{OstingQuadNeumann2026} from $d = 2$ to all $d \ge 2$.

\subsection*{Proof strategy for Theorem~\ref{t:Main}}

All three statements are proved by computing the Hessian of the functional at $\cube$ in symmetry-adapted coordinates, a parameterization adapted to the symmetry group of the $d$-cube, the hyperoctahedral group $B_d := (\Z/2\Z)^d \rtimes S_d$ of order $2^d\cdot d!$. Parts (a) and (b) reduce to shape-derivative computations. Part (c) is handled separately: the lowest Neumann eigenvalue $\mu_1$ is not differentiable at $\cube$ (the cluster splits), so we work with a smooth auxiliary functional and apply degenerate perturbation theory to the eigenvalue cluster. After factoring out the similarity transformations (translations, rotations, and dilations), the tangent space $\TP$ to $\P$ at $\cube$ has dimension
\begin{equation}\label{e:Nd}
N(d) \;:=\; \tfrac{d(3d-1)}{2} \;-\; 1, \qquad d \ge 2.
\end{equation}
Section~\ref{s:CubeParam} constructs $\TP \cong \R^{N(d)}$ inside a linear parameterization $\bdelta \mapsto \Omega_{\bdelta}$ of the larger class $\calQ$ of multilinear perturbations of the cube.

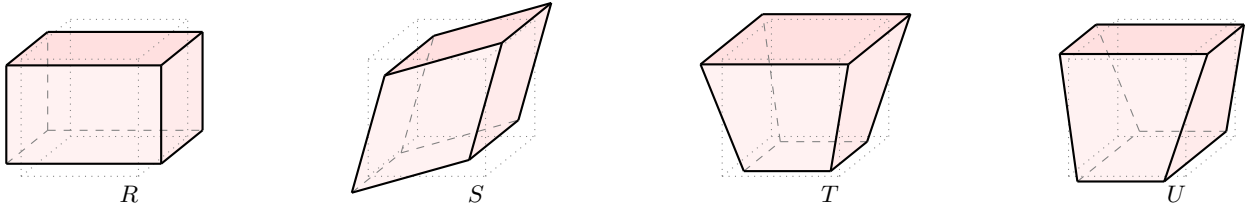
\begin{figure}[t]\centering
\begin{tikzpicture}[scale=1.55, x={(1cm,0cm)},y={(0cm,1cm)},z={(0.42cm,0.34cm)}]
\coordinate (vmmm) at (-0.160,0.080,0.080);
\coordinate (vmmp) at (-0.160,0.920,0.080);
\coordinate (vmpm) at (-0.160,0.080,0.920);
\coordinate (vmpp) at (-0.160,0.920,0.920);
\coordinate (vpmm) at (1.160,0.080,0.080);
\coordinate (vpmp) at (1.160,0.920,0.080);
\coordinate (vppm) at (1.160,0.080,0.920);
\coordinate (vppp) at (1.160,0.920,0.920);
\fill[red!12] (vmmp)--(vpmp)--(vppp)--(vmpp)--cycle;
\fill[red!5] (vmmm)--(vpmm)--(vpmp)--(vmmp)--cycle;
\fill[red!8] (vpmm)--(vppm)--(vppp)--(vpmp)--cycle;
\draw[gray,dotted] (0.000,0.000,0.000) -- (0.000,1.000,0.000);
\draw[gray,dotted] (0.000,0.000,0.000) -- (0.000,0.000,1.000);
\draw[gray,dotted] (0.000,0.000,0.000) -- (1.000,0.000,0.000);
\draw[gray,dotted] (0.000,1.000,0.000) -- (0.000,1.000,1.000);
\draw[gray,dotted] (0.000,1.000,0.000) -- (1.000,1.000,0.000);
\draw[gray,dotted] (0.000,0.000,1.000) -- (0.000,1.000,1.000);
\draw[gray,dotted] (0.000,0.000,1.000) -- (1.000,0.000,1.000);
\draw[gray,dotted] (0.000,1.000,1.000) -- (1.000,1.000,1.000);
\draw[gray,dotted] (1.000,0.000,0.000) -- (1.000,1.000,0.000);
\draw[gray,dotted] (1.000,0.000,0.000) -- (1.000,0.000,1.000);
\draw[gray,dotted] (1.000,1.000,0.000) -- (1.000,1.000,1.000);
\draw[gray,dotted] (1.000,0.000,1.000) -- (1.000,1.000,1.000);
\draw[gray,dashed] (vmmm) -- (vmpm);
\draw[gray,dashed] (vmpm) -- (vmpp);
\draw[gray,dashed] (vmpm) -- (vppm);
\draw[thick] (vmmm) -- (vmmp);
\draw[thick] (vmmm) -- (vpmm);
\draw[thick] (vmmp) -- (vmpp);
\draw[thick] (vmmp) -- (vpmp);
\draw[thick] (vmpp) -- (vppp);
\draw[thick] (vpmm) -- (vpmp);
\draw[thick] (vpmm) -- (vppm);
\draw[thick] (vpmp) -- (vppp);
\draw[thick] (vppm) -- (vppp);
\node[font=\footnotesize] at (0.71,-0.32,0.5) {$R$};
\end{tikzpicture}\hfill
\begin{tikzpicture}[scale=1.55, x={(1cm,0cm)},y={(0cm,1cm)},z={(0.42cm,0.34cm)}]
\coordinate (vmmm) at (-0.140,-0.140,0.000);
\coordinate (vmmp) at (0.140,0.860,0.000);
\coordinate (vmpm) at (-0.140,-0.140,1.000);
\coordinate (vmpp) at (0.140,0.860,1.000);
\coordinate (vpmm) at (0.860,0.140,0.000);
\coordinate (vpmp) at (1.140,1.140,0.000);
\coordinate (vppm) at (0.860,0.140,1.000);
\coordinate (vppp) at (1.140,1.140,1.000);
\fill[red!12] (vmmp)--(vpmp)--(vppp)--(vmpp)--cycle;
\fill[red!5] (vmmm)--(vpmm)--(vpmp)--(vmmp)--cycle;
\fill[red!8] (vpmm)--(vppm)--(vppp)--(vpmp)--cycle;
\draw[gray,dotted] (0.000,0.000,0.000) -- (0.000,1.000,0.000);
\draw[gray,dotted] (0.000,0.000,0.000) -- (0.000,0.000,1.000);
\draw[gray,dotted] (0.000,0.000,0.000) -- (1.000,0.000,0.000);
\draw[gray,dotted] (0.000,1.000,0.000) -- (0.000,1.000,1.000);
\draw[gray,dotted] (0.000,1.000,0.000) -- (1.000,1.000,0.000);
\draw[gray,dotted] (0.000,0.000,1.000) -- (0.000,1.000,1.000);
\draw[gray,dotted] (0.000,0.000,1.000) -- (1.000,0.000,1.000);
\draw[gray,dotted] (0.000,1.000,1.000) -- (1.000,1.000,1.000);
\draw[gray,dotted] (1.000,0.000,0.000) -- (1.000,1.000,0.000);
\draw[gray,dotted] (1.000,0.000,0.000) -- (1.000,0.000,1.000);
\draw[gray,dotted] (1.000,1.000,0.000) -- (1.000,1.000,1.000);
\draw[gray,dotted] (1.000,0.000,1.000) -- (1.000,1.000,1.000);
\draw[gray,dashed] (vmmm) -- (vmpm);
\draw[gray,dashed] (vmpm) -- (vmpp);
\draw[gray,dashed] (vmpm) -- (vppm);
\draw[thick] (vmmm) -- (vmmp);
\draw[thick] (vmmm) -- (vpmm);
\draw[thick] (vmmp) -- (vmpp);
\draw[thick] (vmmp) -- (vpmp);
\draw[thick] (vmpp) -- (vppp);
\draw[thick] (vpmm) -- (vpmp);
\draw[thick] (vpmm) -- (vppm);
\draw[thick] (vpmp) -- (vppp);
\draw[thick] (vppm) -- (vppp);
\node[font=\footnotesize] at (0.71,-0.32,0.5) {$S$};
\end{tikzpicture}\hfill
\begin{tikzpicture}[scale=1.55, x={(1cm,0cm)},y={(0cm,1cm)},z={(0.42cm,0.34cm)}]
\coordinate (vmmm) at (0.130,0.000,0.130);
\coordinate (vmmp) at (-0.130,1.000,-0.130);
\coordinate (vmpm) at (0.130,0.000,0.870);
\coordinate (vmpp) at (-0.130,1.000,1.130);
\coordinate (vpmm) at (0.870,0.000,0.130);
\coordinate (vpmp) at (1.130,1.000,-0.130);
\coordinate (vppm) at (0.870,0.000,0.870);
\coordinate (vppp) at (1.130,1.000,1.130);
\fill[red!12] (vmmp)--(vpmp)--(vppp)--(vmpp)--cycle;
\fill[red!5] (vmmm)--(vpmm)--(vpmp)--(vmmp)--cycle;
\fill[red!8] (vpmm)--(vppm)--(vppp)--(vpmp)--cycle;
\draw[gray,dotted] (0.000,0.000,0.000) -- (0.000,1.000,0.000);
\draw[gray,dotted] (0.000,0.000,0.000) -- (0.000,0.000,1.000);
\draw[gray,dotted] (0.000,0.000,0.000) -- (1.000,0.000,0.000);
\draw[gray,dotted] (0.000,1.000,0.000) -- (0.000,1.000,1.000);
\draw[gray,dotted] (0.000,1.000,0.000) -- (1.000,1.000,0.000);
\draw[gray,dotted] (0.000,0.000,1.000) -- (0.000,1.000,1.000);
\draw[gray,dotted] (0.000,0.000,1.000) -- (1.000,0.000,1.000);
\draw[gray,dotted] (0.000,1.000,1.000) -- (1.000,1.000,1.000);
\draw[gray,dotted] (1.000,0.000,0.000) -- (1.000,1.000,0.000);
\draw[gray,dotted] (1.000,0.000,0.000) -- (1.000,0.000,1.000);
\draw[gray,dotted] (1.000,1.000,0.000) -- (1.000,1.000,1.000);
\draw[gray,dotted] (1.000,0.000,1.000) -- (1.000,1.000,1.000);
\draw[gray,dashed] (vmmm) -- (vmpm);
\draw[gray,dashed] (vmpm) -- (vmpp);
\draw[gray,dashed] (vmpm) -- (vppm);
\draw[thick] (vmmm) -- (vmmp);
\draw[thick] (vmmm) -- (vpmm);
\draw[thick] (vmmp) -- (vmpp);
\draw[thick] (vmmp) -- (vpmp);
\draw[thick] (vmpp) -- (vppp);
\draw[thick] (vpmm) -- (vpmp);
\draw[thick] (vpmm) -- (vppm);
\draw[thick] (vpmp) -- (vppp);
\draw[thick] (vppm) -- (vppp);
\node[font=\footnotesize] at (0.71,-0.32,0.5) {$T$};
\end{tikzpicture}\hfill
\begin{tikzpicture}[scale=1.55, x={(1cm,0cm)},y={(0cm,1cm)},z={(0.42cm,0.34cm)}]
\coordinate (vmmm) at (0.130,0.000,-0.130);
\coordinate (vmmp) at (-0.130,1.000,0.130);
\coordinate (vmpm) at (0.130,0.000,1.130);
\coordinate (vmpp) at (-0.130,1.000,0.870);
\coordinate (vpmm) at (0.870,0.000,-0.130);
\coordinate (vpmp) at (1.130,1.000,0.130);
\coordinate (vppm) at (0.870,0.000,1.130);
\coordinate (vppp) at (1.130,1.000,0.870);
\fill[red!12] (vmmp)--(vpmp)--(vppp)--(vmpp)--cycle;
\fill[red!5] (vmmm)--(vpmm)--(vpmp)--(vmmp)--cycle;
\fill[red!8] (vpmm)--(vppm)--(vppp)--(vpmp)--cycle;
\draw[gray,dotted] (0.000,0.000,0.000) -- (0.000,1.000,0.000);
\draw[gray,dotted] (0.000,0.000,0.000) -- (0.000,0.000,1.000);
\draw[gray,dotted] (0.000,0.000,0.000) -- (1.000,0.000,0.000);
\draw[gray,dotted] (0.000,1.000,0.000) -- (0.000,1.000,1.000);
\draw[gray,dotted] (0.000,1.000,0.000) -- (1.000,1.000,0.000);
\draw[gray,dotted] (0.000,0.000,1.000) -- (0.000,1.000,1.000);
\draw[gray,dotted] (0.000,0.000,1.000) -- (1.000,0.000,1.000);
\draw[gray,dotted] (0.000,1.000,1.000) -- (1.000,1.000,1.000);
\draw[gray,dotted] (1.000,0.000,0.000) -- (1.000,1.000,0.000);
\draw[gray,dotted] (1.000,0.000,0.000) -- (1.000,0.000,1.000);
\draw[gray,dotted] (1.000,1.000,0.000) -- (1.000,1.000,1.000);
\draw[gray,dotted] (1.000,0.000,1.000) -- (1.000,1.000,1.000);
\draw[gray,dashed] (vmmm) -- (vmpm);
\draw[gray,dashed] (vmpm) -- (vmpp);
\draw[gray,dashed] (vmpm) -- (vppm);
\draw[thick] (vmmm) -- (vmmp);
\draw[thick] (vmmm) -- (vpmm);
\draw[thick] (vmmp) -- (vmpp);
\draw[thick] (vmmp) -- (vpmp);
\draw[thick] (vmpp) -- (vppp);
\draw[thick] (vpmm) -- (vpmp);
\draw[thick] (vpmm) -- (vppm);
\draw[thick] (vpmp) -- (vppp);
\draw[thick] (vppm) -- (vppp);
\node[font=\footnotesize] at (0.71,-0.32,0.5) {$U$};
\end{tikzpicture}
\caption{Representative deformations in the four irreducible blocks of Theorem~\ref{t:irrep-decomp} at $d = 3$, drawn as the deformed cubes $(\mathrm{id} + t\bth)(\cube)$ for a representative field $\bth$ in each block. $R$: $\bth_R = (2x_1, -x_2, -x_3)$, a traceless stretch. $S$: the pure shear $(x_3, 0, x_1)$, carrying the cube to a rhombic prism. $T$: $\bth_T = (4x_1 x_3, 4x_2 x_3, 0)$ --- the top facet dilates, the bottom contracts (a frustum). $U$: $\bth_U = (4x_1 x_3, -4x_2 x_3, 0)$ --- the top facet stretches along $x_1$ and contracts along $x_2$, the bottom facet oppositely.}
\label{fig:modes}
\end{figure}

The tangent space $\TP$ decomposes as a $B_d$-representation into four pairwise non-isomorphic irreducible pieces (Theorem~\ref{t:irrep-decomp}), 
% proved in~\cite{CubeRep2026}, 
\begin{equation}\label{e:irrep-decomp}
\TP \;=\; R^{(d)} \,\oplus\, S^{(d)} \,\oplus\, T^{(d)} \,\oplus\, U^{(d)} \;\cong\; \R^{N(d)} , \qquad d \ge 2.
\end{equation}
Geometrically, $R^{(d)}$ consists of traceless axis-aligned stretches, which deform the cube to \emph{rectangular} boxes; $S^{(d)}$ of traceless pairwise \emph{shears} along the diagonal directions $\tfrac{1}{\sqrt2}(\mathbf e_i \pm \mathbf e_j)$, which deform the cube to rhombic parallelepipeds; and $T^{(d)}, U^{(d)}$ are two distinct families of \emph{trapezoidal} deformations. The dimensions are $\dim R^{(d)} = d-1$,  $\dim S^{(d)} = \tbinom{d}{2}$,  $\dim T^{(d)} = d$,  and $\dim U^{(d)} = d(d-2)$, summing to $N(d)$. 
These deformations are illustrated for $d = 3$ in Figure~\ref{fig:modes}. 
We call coordinates on this decomposition \emph{symmetry-adapted}. 
At $d = 2$ the block $U^{(2)}$ is empty, since $d(d-2) = 0$. 

By Schur's lemma, the Hessian of any $B_d$-invariant smooth functional $F$ at $\cube$ is block-diagonal in symmetry-adapted coordinates, acting as a scalar $\nu^F_X$ on each block $X \in \{R, S, T, U\}$ of~\eqref{e:irrep-decomp} (Proposition~\ref{p:Crit-Block}). For each functional $F$ in~\eqref{e:three-functionals-intro} (with $\bmu$ replaced by the smooth auxiliary functional $f_d$ introduced below), the Hessian is therefore completely characterized by the four scalar eigenvalues $\nu^F_R, \nu^F_S, \nu^F_T, \nu^F_U$.

For $I$ and $\blam$ (parts (a)--(b)), both functionals are smooth on $\P$ and amenable to the shape-derivative framework (\cite{HenrotPierre2018, Laurain2020}; see also~\cite{BogoselBucur2024} for the polygonal version), recalled in Section~\ref{ss:Shape-prelim}. For $I$, the perturbed bodies along the four representative deformations are explicit polytopes whose volume and surface area are elementary functions of $t$, yielding the four eigenvalues
\begin{equation}\label{e:Iso-eigs}
\nu^I_R = 4d^d, 
\qquad 
\nu^I_S = 8d^d, 
\qquad 
\nu^I_T = \tfrac{4(d+1)}{3}d^d, 
\qquad 
\nu^I_U = \tfrac{8}{3}d^d \ (d \ge 3),
\end{equation}
all strictly positive (Proposition~\ref{p:Iso-irrep}).

For $\blam$, writing $\lambda(t) := \lambda_1(\Omega_t)$ along the flow $\Omega_t = (\mathrm{id} + t\bth)(\Omega)$ and $u$ for the normalized first eigenfunction, the classical second shape derivative formula
\begin{equation}\label{e:dL1-BB}
\lambda''(0)[\bth, \bth] \;=\; \int_\Omega \bigl[2|\nabla \dot u|^2 - 2\lambda(0)|\dot u|^2\bigr] dx \,-\, 2\lambda'(0)[\bth]\, \langle \dot u, u\rangle
\end{equation}
(with $\dot u$ the material derivative of $u$; see~\cite{HenrotPierre2018,Laurain2020}) reduces to four Hessian eigenvalues, which depend on $d$ and two universal series-constants. The result is:
\begin{theorem}[Hessian of $\blam$]\label{t:Lam-Hessian-d}
For all $d \ge 2$, the Hessian eigenvalues of $\blam$ at $\cube$ in symmetry-adapted coordinates, on the four blocks of~\eqref{e:irrep-decomp}, are
\begin{equation}\label{e:nu-Dirichlet}
\begin{gathered}
\nu^{\blam}_R = \frac{16\pi^2}{2^d}, \qquad
\nu^{\blam}_S = \frac{4096\,D_S}{2^d\,\pi}, \\[4pt]
\nu^{\blam}_T = \frac{1}{2^d}\Bigl[\Bigl(\tfrac{16\pi^2}{3} - 40\Bigr)(d-2) + \frac{2048\,D_U}{\pi}\Bigr], \qquad
\nu^{\blam}_U = \frac{1}{2^d}\Bigl[\frac{2048\,D_U}{\pi} - \Bigl(\tfrac{16\pi^2}{3} - 40\Bigr)\Bigr],
\end{gathered}
\end{equation}
where $D_S \approx 0.072873$ and $D_U \approx 0.069587$ are positive universal series-constants, given explicitly in~\eqref{e:Dirichlet-sums}. All four eigenvalues are strictly positive for every $d \ge 2$ (the $U$ block being empty at $d = 2$).
\end{theorem}

\noindent Theorem~\ref{t:Lam-Hessian-d} is proved in Section~\ref{s:LDshape}. 

For $\bmu$ (part (c)), the cluster $\mu_1 = \cdots = \mu_d = \pi^2$ at $\cube$ prevents direct shape differentiation. Following~\cite{OstingQuadNeumann2026}, we work with the smooth auxiliary functional
\begin{equation}\label{e:fd-defn}
f_d(\Omega) \;:=\; \tfrac{1}{d|\Omega|^{2/d}}\sum_{i=1}^d \tfrac{1}{\mu_i(\Omega)},
\end{equation}
which satisfies $f_d \le \tfrac{1}{\bmu}$ with equality at $\cube$. A strict local minimum of $f_d$ therefore implies a strict local maximum of $\bmu$ (Remark~\ref{r:fd-bound}), and $f_d$ is smooth through the cluster degeneracy (Section~\ref{ss:Aux}). The Hessian of $f_d$ is computed by degenerate perturbation theory of the eigenvalue cluster (Section~\ref{s:LNshape}); the result is:

\begin{theorem}[Hessian of $f_d$]\label{t:fd-Hessian-d}
For all $d \ge 2$, the Hessian eigenvalues of $f_d$ at $\cube$ in symmetry-adapted coordinates, on the four blocks of~\eqref{e:irrep-decomp}, are
\begin{equation}\label{e:nu-Neumann}
\begin{gathered}
\nu^{f_d}_R \;=\; \frac{16}{d \cdot 2^d \,\pi^2}, \qquad
\nu^{f_d}_S \;=\; \frac{2560 - 16\pi^4 + 1024\,\pi\, \ES}{d \cdot 2^d \,\pi^6}, \\[4pt]
\nu^{f_d}_T \;=\; \frac{16\pi^6 + (40d-80)\pi^4 + (600-120d)\pi^2 + 3840\,\pi\,\ETU}{15\, d \cdot 2^d \, \pi^6}, \\[4pt]
\nu^{f_d}_U \;=\; \frac{16\pi^6 - 40\pi^4 + 480\pi^2 + 3840\,\pi\,\ETU}{15\, d \cdot 2^d \,\pi^6} \quad (d \ge 3),
\end{gathered}
\end{equation}
where $\ES$ and $\ETU$ are positive universal series-constants, given explicitly in~\eqref{e:Neumann-sums}. All four eigenvalues are strictly positive for every $d \ge 2$ (the $U$ block being empty for $d = 2$).
\end{theorem}

\noindent Theorem~\ref{t:fd-Hessian-d} is proved in Section~\ref{s:LNshape}.

\begin{remark}[Local versus global]\label{r:local-global}
Theorem~\ref{t:Main} is a \emph{local} statement: the cube is a strict local extremizer of each functional on $\P$. This is weaker than the planar quadrilateral result~\cite{OstingQuadNeumann2026}, which is global---the square is the \emph{unique} maximizer of $\bmu$ among all convex quadrilaterals. The planar proof promotes local to global in two steps.  Kr\"oger's diameter inequality~\cite{Kroger1999, HenrotMichetti2024} confines any candidate maximizer to a compact region of the four-dimensional parameter space; that region is then covered by a computer-assisted box search, certifying on each box an upper bound on a Rayleigh--Ritz eigenvalue computed on the fixed reference square. The same strategy is available in principle for each fixed $d$, but the region to be covered is now $N(d)$-dimensional, so the number of boxes grows exponentially in $d$. This strategy therefore appears out of reach, and the global result remains open for $d \ge 3$.  
\end{remark}

The restriction to \emph{polytopes} in Theorem~\ref{t:Main} is crucial: the following lemma exhibits a simple \emph{quadratic} perturbation showing that, among all convex bodies, the cube is not even a critical point of any of the three functionals.

\begin{lemma}\label{l:dish}
Let $V_n(\x) := \tfrac{d+2}{12} - |\x|^2$ on $\partial\cube$, and for $t \ge 0$ let
\begin{equation}\label{e:dish-body}
\Omega_t \;:=\; \bigcap_{k=1}^{d}\ \bigcap_{\sigma = \pm 1}\Bigl\{\x \in \R^d \,:\, \sigma x_k \le \tfrac12 + t\Bigl(\tfrac{d-1}{12} - \textstyle\sum_{l \ne k} x_l^2\Bigr)\Bigr\}
\end{equation}
be the body obtained by bending each facet of $\cube$ into the graph of $t\,V_n$.  Then $\Omega_t$ is convex for all small $t \ge 0$, $|\Omega_t| = 1 + O(t^2)$, and, as $t \downarrow 0$,
\begin{equation}\label{e:dish-derivs}
\frac{d}{dt}\Big|_{t=0^+} I(\Omega_t) = -\,\frac{d(d-1)}{3}\,(2d)^d,
\quad
\frac{d}{dt}\Big|_{t=0^+} \blam(\Omega_t) = -\,2d(d-1),
\quad
\frac{d}{dt}\Big|_{t=0^+} \bmu(\Omega_t) = +\,2(d-1).
\end{equation}
In particular, for every $d \ge 2$, the cube is not a stationary point of $I$, $\blam$, or $\bmu$ among convex bodies.
\end{lemma}

\noindent Lemma~\ref{l:dish} is proved in Section~\ref{s:Completion}.  

On the class $\P$ of flat-facet combinatorial cubes, the cube is a strict local extremizer of all three functionals above (Theorem~\ref{t:Main}); among all convex bodies it is a critical point of none of them.  This raises the ``Goldilocks'' question of the largest class on which the cube remains critical.

\subsection{Multilinear perturbations} \label{ss:multPert} 
\begin{figure}[t]\centering
\begin{tikzpicture}[scale=2.0,every node/.style={font=\footnotesize},
  x={(1cm,0cm)},y={(0cm,1cm)},z={(0.42cm,0.34cm)}]
% front face (z=0): o a b c; back face (z=1): d e f0 g; corner f0 pressed in to fp
\coordinate (o) at (0,0,0); \coordinate (a) at (1,0,0); \coordinate (b) at (1,1,0); \coordinate (c) at (0,1,0);
\coordinate (d) at (0,0,1); \coordinate (e) at (1,0,1); \coordinate (g) at (0,1,1);
\coordinate (f0) at (1,1,1);
\coordinate (fp) at (1,0.72,1);
% the split top facet: flat triangle (c,b,g) and tilted triangle (b,fp,g)
\fill[red!6]  (c)--(b)--(g)--cycle;
\fill[red!14] (b)--(fp)--(g)--cycle;
% hidden edges (the three edges meeting the occluded back vertex d)
\draw[gray,dashed] (d)--(o) (d)--(e) (d)--(g);
% unperturbed corner and its former edges
\draw[gray,dotted,thick] (b)--(f0) (g)--(f0) (fp)--(f0);
\draw[gray,fill=white] (f0) circle (1.1pt);
% visible edges of the perturbed polytope
\draw[thick] (o)--(a)--(b)--(c)--(o) (a)--(e) (c)--(g);
\draw[thick] (e)--(fp) (b)--(fp) (g)--(fp);
% inserted crease
\draw[red,thick] (b)--(g) node[midway, sloped, above=1pt, red] {crease};
% displacement of the corner
\draw[red,->] (1.52,1.32) -- (1.52,1.10);
\node[red, anchor=west] at (1.58,1.21) {dropped corner};
\end{tikzpicture}\quad
\begin{tikzpicture}[scale=2.0,every node/.style={font=\footnotesize},
  x={(1cm,0cm)},y={(0cm,1cm)},z={(0.42cm,0.34cm)}]
\coordinate (wmmm) at (0.085,0.085,0.085);
\coordinate (wmmp) at (-0.085,-0.085,1.085);
\coordinate (wmpm) at (-0.085,1.085,-0.085);
\coordinate (wmpp) at (0.085,0.915,0.915);
\coordinate (wpmm) at (1.085,-0.085,-0.085);
\coordinate (wpmp) at (0.915,0.085,0.915);
\coordinate (wppm) at (0.915,0.915,0.085);
\coordinate (wppp) at (1.085,1.085,1.085);
\fill[red!6] (0.085,0.085,0.085)--(1.085,-0.085,-0.085)--(0.915,0.915,0.085)--(-0.085,1.085,-0.085)--cycle;
\draw[gray!55, line width=0.3pt] (0.335,0.043,0.043) -- (0.165,1.042,-0.042);
\draw[gray!55, line width=0.3pt] (0.585,0.000,0.000) -- (0.415,1.000,0.000);
\draw[gray!55, line width=0.3pt] (0.835,-0.042,-0.042) -- (0.665,0.958,0.043);
\fill[red!14] (-0.085,1.085,-0.085)--(0.915,0.915,0.085)--(1.085,1.085,1.085)--(0.085,0.915,0.915)--cycle;
\draw[gray!55, line width=0.3pt] (0.165,1.042,-0.042) -- (0.335,0.958,0.958);
\draw[gray!55, line width=0.3pt] (0.415,1.000,0.000) -- (0.585,1.000,1.000);
\draw[gray!55, line width=0.3pt] (0.665,0.958,0.043) -- (0.835,1.042,1.042);
\fill[red!9] (1.085,-0.085,-0.085)--(0.915,0.085,0.915)--(1.085,1.085,1.085)--(0.915,0.915,0.085)--cycle;
\draw[gray!55, line width=0.3pt] (1.042,-0.042,0.165) -- (0.958,0.958,0.335);
\draw[gray!55, line width=0.3pt] (1.000,0.000,0.415) -- (1.000,1.000,0.585);
\draw[gray!55, line width=0.3pt] (0.958,0.043,0.665) -- (1.042,1.042,0.835);
% unperturbed cube
\draw[gray,dotted] (0,0,0) -- (1,0,0);
\draw[gray,dotted] (0,0,0) -- (0,1,0);
\draw[gray,dotted] (0,0,0) -- (0,0,1);
\draw[gray,dotted] (0,0,1) -- (1,0,1);
\draw[gray,dotted] (0,0,1) -- (0,1,1);
\draw[gray,dotted] (0,1,0) -- (1,1,0);
\draw[gray,dotted] (0,1,0) -- (0,1,1);
\draw[gray,dotted] (0,1,1) -- (1,1,1);
\draw[gray,dotted] (1,0,0) -- (1,1,0);
\draw[gray,dotted] (1,0,0) -- (1,0,1);
\draw[gray,dotted] (1,0,1) -- (1,1,1);
\draw[gray,dotted] (1,1,0) -- (1,1,1);
% edges of the perturbed body (straight lines under multilinear interpolation)
\draw[thick] (wmmm) -- (wpmm);
\draw[thick] (wmmm) -- (wmpm);
\draw[gray,dashed] (wmmm) -- (wmmp);
\draw[gray,dashed] (wmmp) -- (wpmp);
\draw[gray,dashed] (wmmp) -- (wmpp);
\draw[thick] (wmpm) -- (wppm);
\draw[thick] (wmpm) -- (wmpp);
\draw[thick] (wmpp) -- (wppp);
\draw[thick] (wpmm) -- (wppm);
\draw[thick] (wpmm) -- (wpmp);
\draw[thick] (wpmp) -- (wppp);
\draw[thick] (wppm) -- (wppp);
% vertex markers: outward (red) / inward (blue)
\fill[blue] (0.085,0.085,0.085) circle (1.1pt);
\fill[red] (-0.085,-0.085,1.085) circle (1.1pt);
\fill[red] (-0.085,1.085,-0.085) circle (1.1pt);
\fill[blue] (0.085,0.915,0.915) circle (1.1pt);
\fill[red] (1.085,-0.085,-0.085) circle (1.1pt);
\fill[blue] (0.915,0.085,0.915) circle (1.1pt);
\fill[blue] (0.915,0.915,0.085) circle (1.1pt);
\fill[red] (1.085,1.085,1.085) circle (1.1pt);
\node[red, anchor=west, font=\scriptsize]  at (1.55,1.02) {$\sigma_1\sigma_2\sigma_3 = +1$: outward};
\node[blue, anchor=west, font=\scriptsize] at (1.55,0.84) {$\sigma_1\sigma_2\sigma_3 = -1$: inward};
\end{tikzpicture}
\caption{Exiting $\P$ at $d = 3$. 
\emph{Left}: a corner-drop perturbation, with the perturbed body taken to be the convex hull of the displaced vertices.  
\emph{Right}: the deformation $(\mathrm{id} + t\bth^{*})(\square_3)$ of~\eqref{e:top-field}  in the class $\calQ$. 
See Section~\ref{ss:multPert}.}
\label{fig:pexit}
\end{figure}
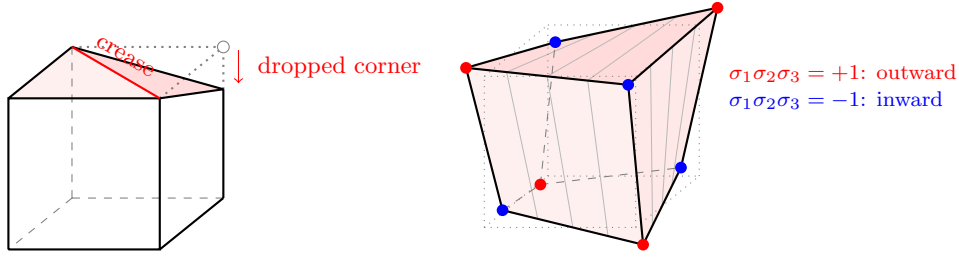

In dimensions $d\geq 3$, there are two natural ways to extend $\P$ so as to admit arbitrary perturbations of the $2^d$ vertices of $\cube$.  For simplicity, let us first consider $d = 3$ (see Figure~\ref{fig:pexit}).  Press one corner of $\square_3$ slightly into the cube, normal to one of its three facets, holding the other seven vertices fixed.  The convex hull splits the affected facet into two triangles along the diagonal joining the displaced corner's two neighbors: a crease is inserted, the polytope gains a facet, and it is no longer combinatorially a cube (Figure~\ref{fig:pexit}, left).  If instead we preserve the original combinatorics by multilinear interpolation, the facet survives as a single four-sided face, but its four corners are no longer coplanar: the bilinear patch through them is a saddle surface, and convexity is lost (Figure~\ref{fig:pexit}, right).  Either way, the perturbation exits $\P$.

We take the second route, and consider the class $\calQ$ of \emph{all} multilinear vertex-perturbations of $\cube$ (Section~\ref{s:CubeParam}), with no combinatorial restriction, together with its parameter space $\TQ \supseteq \TP$.   
While convexity of the perturbed shape is lost, smoothness of the perturbation is preserved: each flow is a polynomial family of bodies, and the shape calculus of Sections~\ref{s:Iso}--\ref{s:LNshape} applies along it.  On $\P$, Theorem~\ref{t:Main} makes the cube a strict local extremizer of all three functionals~\eqref{e:three-functionals-intro}.  On $\calQ$, a \emph{single} deformation rules out local optimality of the cube for all three functionals in high dimension. We consider the fully symmetric corner perturbation
\begin{equation}\label{e:top-field}
\bth^{*} \;:=\; 2^{\,d-1}\, \nabla (x_1 x_2 \cdots x_d) ,
\end{equation}
which displaces the vertex $\v_{\bsigma} = \tfrac12\bsigma$, $\bsigma \in \{-1, +1\}^d$, by $(\sigma_1 \cdots \sigma_d)\, \bsigma$: every vertex slides along its own diagonal, outward when $\sigma_1 \cdots \sigma_d = +1$ and inward when $\sigma_1 \cdots \sigma_d = -1$, so that every facet becomes a bilinear saddle.  For $d = 3$ this is exactly the deformation of Figure~\ref{fig:pexit} (right).

In the following statement, $F''(0)$ for $F = \bmu$ is well-defined even though $\mu_1$ is not differentiable at $\cube$: the symmetry of $\bth^{*}$ keeps the lowest Neumann cluster unsplit through second order (Lemma~\ref{l:B-cluster}), so $t \mapsto \bmu(\Omega_t)$ is smooth along this deformation.

\begin{theorem}\label{t:top-mode}
Let $d \ge 3$, let $\Omega_t := (\mathrm{id} + t \bth^{*})(\cube)$ with $\bth^{*}$ the corner perturbation~\eqref{e:top-field}, and write $F''(0) := \tfrac{d^2}{dt^2}\big|_{t=0} F(\Omega_t)$.  Each of $I$, $\blam$, $\bmu$ is critical at $t = 0$ along the flow $\Omega_t$, and
\begin{enumerate}
\item[(a)] $I''(0) = -\,8\,d(d-1)(d-3)\, 3^{\,2-d}\, I(\cube)$, which is strictly negative for every $d \ge 4$;
\item[(b)] $\blam''(0) < 0$ for every $d \ge 3$;
\item[(c)] $\bmu''(0) > 0$ for every $d \ge 6$.
\end{enumerate}
Consequently, modulo similarity, $\cube$ is not a local minimizer of $I$ on $\calQ$ for any $d \ge 4$ and not a local maximizer of $\bmu$ on $\calQ$ for any $d \ge 6$; and $\cube$ is a strict local minimizer of $\blam$ on $\calQ$ if and only if $d = 2$.
\end{theorem}

Theorem~\ref{t:top-mode} is proved in Section~\ref{s:TopMode}.  The thresholds in $d$ are sharp for this deformation ($I''(0) = 0$ at $d = 3$, and numerically $\bmu''(0) < 0$ for $d = 3, 4, 5$), and for $d = 2$ the class $\calQ$ coincides with $\P$ near the square, so Theorem~\ref{t:Main} applies.  Based on numerical computations, we conjecture that, modulo similarity, $\cube$ is a strict local minimizer of $I$ on $\calQ$ if and only if $d \in \{2, 3\}$, and a strict local maximizer of $\bmu$ on $\calQ$ if and only if $2 \le d \le 5$.

\subsection*{Outline}
Section~\ref{s:CubeParam} constructs the symmetry-adapted parameterization of $\P$, states its decomposition into irreducible $B_d$-representations, and treats the rectangular deformation, for which everything is explicit, in closed form. Sections~\ref{s:Iso}, \ref{s:LDshape}, and~\ref{s:LNshape} compute the Hessians of $I$, $\blam$, and $f_d$, respectively. Section~\ref{s:Completion} assembles the proof of Theorem~\ref{t:Main} and proves Lemma~\ref{l:dish}, and Section~\ref{s:TopMode} proves Theorem~\ref{t:top-mode}.

%%%%%%%%%%%%%%%%%%%%%%%%%%%%%%%%%%
%%%%%%%%%%%%%%%%%%%%%%%%%%%%%%%%%%
\section{Symmetry-adapted parameterization of the $d$-cube}\label{s:CubeParam}

We parameterize the class $\tP$ of combinatorial $d$-cubes near $\cube$ and record the representation-theoretic structure that diagonalizes the Hessian of any $B_d$-invariant smooth functional at the cube.

\subsection{Parameterization by vertex displacements}\label{ss:CubeParam-vertex}

The $2^d$ vertices of $\cube$ are indexed by sign vectors $\bsigma = (\sigma_1, \ldots, \sigma_{d}) \in \{-1, +1\}^d$ via $\v_{\bsigma} := \tfrac{1}{2}\bsigma \in \R^d$. 
Define, for each $\bsigma$, the multilinear ``hat function'' $\varphi_{\bsigma} \colon \cube \to [0,1]$, by $\varphi_{\bsigma}(\x) := \prod_{k=1}^{d}\bigl(\tfrac{1}{2} + \sigma_k x_k\bigr),  \x = (x_1, \ldots, x_{d}) \in \cube$. These satisfy $\varphi_{\bsigma}(\v_{\bsigma'}) = \delta_{\bsigma, \bsigma'}$ and $\sum_{\bsigma} \varphi_{\bsigma} \equiv 1$.
A \emph{vertex displacement} is a tuple $\bdelta = (\bdelta_{\bsigma})_{\bsigma \in \{-1, +1\}^d}$ with $\bdelta_{\bsigma} \in \R^d$; the space $(\R^d)^{2^d}$ of vertex displacements carries the Euclidean inner product $\langle \bdelta, \bdelta' \rangle := \sum_{\bsigma} \bdelta_{\bsigma} \cdot \bdelta'_{\bsigma}$, with norm $\|\bdelta\| := \langle \bdelta, \bdelta \rangle^{1/2}$. For a vertex displacement $\bdelta$, the multilinear interpolant $\Phi_{\bdelta}$ and the perturbed body $\Omega_{\bdelta}$ are
\begin{equation}\label{e:Phi-delta}
\Phi_{\bdelta}(\x) \;:=\; \sum_{\bsigma \in \{-1, +1\}^d} \varphi_{\bsigma}(\x)\,\bigl(\v_{\bsigma} + \bdelta_{\bsigma}\bigr), \qquad 
\Omega_{\bdelta} \;:=\; \Phi_{\bdelta}(\cube).
\end{equation}
For $\|\bdelta\|$ sufficiently small, $\Phi_{\bdelta}$ is a diffeomorphism onto its image, and we use it throughout to pull back integrals and differential operators from $\Omega_{\bdelta}$ to $\cube$. Each facet of $\Omega_{\bdelta}$ is the multilinear image of a facet of $\cube$: a $(d-1)$-dimensional patch with flat edges, in general curved (a bilinear ``saddle'' at $d = 3$), so $\Omega_{\bdelta}$ need not be convex.

For a fixed small $\varepsilon > 0$, define the classes
\begin{equation}\label{e:P-defn}
\begin{aligned}
\tQ \;&:=\; \bigl\{ \Omega_{\bdelta} = \Phi_{\bdelta}(\cube) \,:\, \|\bdelta\| < \varepsilon \bigr\},\\
\tP \;&:=\; \bigl\{ \Omega_{\bdelta} \in \tQ \,:\, \Omega_{\bdelta} \text{ is a convex polytope combinatorially equivalent to } \cube \bigr\}.
\end{aligned}
\end{equation}
Membership $\Omega_{\bdelta} \in \tP$ holds precisely when, for each facet of $\cube$, the $2^{d-1}$ displaced vertices remain coplanar; in that case every facet of $\Omega_{\bdelta}$ is flat, $\Omega_{\bdelta}$ coincides with the convex hull of the displaced vertices, and its face lattice is isomorphic to that of $\cube$. When some facet's displaced vertices are \emph{not} coplanar, $\Omega_{\bdelta}$ leaves $\tP$ but remains in $\tQ$: the multilinear body keeps the facet as a single curved face and loses convexity, while the convex hull of the same vertices would instead split the facet along a crease (Figure~\ref{fig:pexit}).

Two bodies in $\tQ$ that differ by a similarity transformation (a translation, rotation, or dilation) have the same values of $I$, $\blam$, and $\bmu$, so we work modulo similarity. The $d \cdot 2^d$-dimensional space $(\R^d)^{2^d}$ of all vertex displacements contains the $\bigl(d + \binom{d}{2} + 1\bigr)$-dimensional subspace $\sm := \bigl\{ \bigl(\mathbf{c} + (J + \lambda I)\v_{\bsigma}\bigr)_{\bsigma} \,:\, \mathbf{c} \in \R^d,\ J \in \mathfrak{so}(d),\ \lambda \in \R \bigr\}$ of \emph{similarity directions}: the vertex velocity fields of translations, rotations, and dilations, equivalently the tangent space at $\cube$ of its orbit under similarity transformations. We normalize by restricting to displacements orthogonal to $\sm$, and define the \emph{normalized classes}
\begin{equation}\label{e:M-defn}
\calQ \;:=\; \bigl\{ \Omega_{\bdelta} \in \tQ \,:\, \bdelta \perp \sm \bigr\}, \qquad
\P \;:=\; \tP \cap \calQ .
\end{equation}
Since \emph{every} small displacement produces a body in $\tQ$, the class $\calQ$ is identified with a neighborhood of $\bzero$ in the linear space
\begin{equation}\label{e:Q-defn}
\TQ \;:=\; \sm^{\perp} \;=\; \bigl\{ \bdelta \in (\R^d)^{2^d} \,:\, \langle \bdelta, \btau \rangle = 0 \text{ for all } \btau \in \sm \bigr\},
\qquad
\dim \TQ \;=\; d \cdot 2^d - d - \tbinom{d}{2} - 1 ,
\end{equation}
which is its tangent space at $\cube$. The class $\P$ is a subspace of $\calQ$ constrained by facet-coplanarity.

\begin{lemma}\label{l:Pd-dim}
Let $d \ge 2$. Near $\bdelta = \bzero$, the set $\{\bdelta : \Omega_{\bdelta} \in \tP\}$ is a smooth submanifold of $(\R^d)^{2^d}$ of dimension $2d^2$ whose tangent space at $\bzero$ contains $\sm$, and $\P$ is a smooth submanifold of $\calQ$ of dimension $N(d)$ of~\eqref{e:Nd}. Every body in $\tP$ near $\cube$ is similar to a body in $\P$.
\end{lemma}

\begin{proof}
%By~\cite[Lemma~3.1]{CubeRep2026}, 
The set $\{\bdelta : \Omega_{\bdelta} \in \tP\}$ is, near $\bdelta = \bzero$, a smooth submanifold of dimension $2d^2$ whose tangent space at $\bzero$ contains $\sm$ (the similarity orbit of $\cube$ lies in $\tP$). It therefore meets $\TQ = \sm^{\perp}$ transversally at $\bzero$, so $\P$ is a submanifold of $\calQ$ of dimension $2d^2 - \dim\sm = N(d)$. The map $(g, \Omega) \mapsto g\,\Omega$ from (similarities) $\times\, \P$ to $\tP$ has differential $\sm \oplus T_{\bzero}\P$ at $(\mathrm{id}, \cube)$, the whole tangent space of $\tP$ at $\cube$, so by the inverse function theorem it is a local diffeomorphism, and every body in $\tP$ near $\cube$ is of the form $g\,\Omega$ with $\Omega \in \P$.
\end{proof}

The tangent space of $\P$ at $\cube$, in the coordinates $\bdelta \in \TQ$, is
\begin{equation}\label{e:hatP-defn}
\TP \;:=\; T_{\bzero}\P \;\subset\; \TQ, \qquad \dim \TP \;=\; N(d).
\end{equation}
A direction $\bth \in \TQ \setminus \TP$, with $\TQ$ as in~\eqref{e:Q-defn} --- an excess of dimension $d \cdot 2^d - 2d^2$, present for every $d \ge 3$ --- breaks the coplanarity of some facet's vertices at first order, so $\Omega_{t\bth}$ leaves $\P$, while remaining in $\calQ$, for all small $t \ne 0$. A distinguished direction of this kind is the corner perturbation~\eqref{e:top-field} of the Introduction.

\subsection{Hyperoctahedral symmetry and representative vector fields}\label{ss:Bd-action}\label{ss:Reps}

The symmetry group of $\cube$ is the \emph{hyperoctahedral group} $B_d \cong (\Z/2\Z)^d \rtimes S_d$, of order $2^d \cdot d!$. An element of $B_d$ is a pair $((\varepsilon_1, \ldots, \varepsilon_{d}), \pi)$ with $\varepsilon_k \in \{-1, +1\}$ and $\pi \in S_d$, acting on $\cube$ by axis permutations composed with sign flips. Equivalently, $B_d$ acts on the vertex index set $\{-1, +1\}^d$ by permuting coordinates and flipping signs. It acts on vertex displacements by $(g \cdot \bdelta)_{\bsigma} := g\,\bdelta_{g^{-1}\bsigma}$, the action induced on vertex velocity fields; in particular a sign flip $\bepsilon$ acts by $(\bepsilon \cdot \bdelta)_{\bsigma}[k] = \varepsilon_k\, \bdelta_{\bepsilon\bsigma}[k]$, where $\bdelta_{\bsigma}[k]$ denotes the $k$th component of $\bdelta_{\bsigma}$, and an axis permutation permutes axes and vertex coordinates simultaneously. This action is orthogonal for the inner product of Section~\ref{ss:CubeParam-vertex} and preserves $\sm$, $\TQ$ and $\TP$.

\begin{remark}[Reduction to four scalars]\label{r:irrep-intuition}
The Hessian of a shape functional at $\cube$ is an $N(d) \times N(d)$ symmetric matrix  and computing it entry by entry for every $d$ is not feasible.  Because the functional is $B_d$-invariant, its Hessian \emph{commutes} with the $B_d$-action on the parameter space, and Schur's lemma forces any such symmetric form to vanish between non-isomorphic irreducible subrepresentations and to act as a scalar multiple of the identity on each one. The parameter space splits into exactly four irreducible representations (Theorem~\ref{t:irrep-decomp}), which we name by the shapes they deform the cube into: $R$ (rectangular stretches), $S$ (rhombic pure shears), and $T, U$ (two families of trapezoidal deformations). In symmetry-adapted coordinates the entire Hessian of any cube-symmetric functional is therefore captured by four numbers $\nu_R, \nu_S, \nu_T, \nu_U$, regardless of $d$, and positive definiteness reduces to four sign checks. The four representative fields introduced below are one convenient vector chosen from each block and 
are illustrated in Figure~\ref{fig:modes}.
\end{remark}

\begin{theorem}
%[{\cite[Theorem~3.4]{CubeRep2026}}]
\label{t:irrep-decomp}
Let $d \ge 2$. The space $\TP$ of~\eqref{e:hatP-defn} decomposes as in~\eqref{e:irrep-decomp}, $\TP = R^{(d)} \oplus S^{(d)} \oplus T^{(d)} \oplus U^{(d)}$, as a direct sum of four pairwise non-isomorphic irreducible $B_d$-representations, of dimensions $d-1$, $\tbinom{d}{2}$, $d$, and $d(d-2)$, characterized geometrically as follows; in the spanning vectors, components of $\bdelta_{\bsigma}$ that are not listed are zero.
\begin{itemize}
\item $R^{(d)}$, the $(d-1)$-dim space of traceless axis-aligned stretches.  Spanned by $\bdelta_{\bsigma} = \sum_k a_k \sigma_k \mathbf{e}_k$ with $\sum_k a_k = 0$.
\item $S^{(d)}$, the $\binom{d}{2}$-dim space of symmetric pairwise shears.  Spanned by $\bdelta_{\bsigma}[i] = \sigma_j, \bdelta_{\bsigma}[j] = \sigma_i$ for ordered pairs $i < j$; each generator is a traceless stretch along the diagonal $\tfrac{1}{\sqrt2}(\mathbf e_i + \mathbf e_j)$ (compression along $\tfrac{1}{\sqrt2}(\mathbf e_i - \mathbf e_j)$), carrying $\cube$ to a rhombic parallelepiped (a rhombus in the $(x_i, x_j)$-plane).  It is the symmetric counterpart of the antisymmetric rotations $\bdelta_{\bsigma}[i] = \sigma_j, \bdelta_{\bsigma}[j] = -\sigma_i$, which are trivial.
\item $T^{(d)}$, the $d$-dim space of ``symmetric trapezoid'' deformations.  Spanned, for a \emph{modulator axis} $c \in \{1, \ldots, d\}$, by $\bdelta_{\bsigma} = \sigma_c \sum_{a \ne c} \sigma_a\, \mathbf{e}_a$.
\item $U^{(d)}$, the $d(d-2)$-dim space of ``antisymmetric trapezoid'' deformations.  For each modulator axis $c$, spanned by $\bdelta_{\bsigma} = \sum_{a \ne c} t_a\, \sigma_a \sigma_c\, \mathbf{e}_a$ with $\sum_{a \ne c} t_a = 0$ (the traceless counterpart of the $T$ generators, which have all $t_a$ equal).
\end{itemize}
For $d = 2$ the summand $U^{(2)}$ is zero.
\end{theorem}

\begin{remark}
The proof of Theorem~\ref{t:irrep-decomp} is given in a companion paper, in preparation, and summarized as follows. Each summand is spanned by the displayed vectors and is invariant under coordinate permutations and sign changes, the summands are mutually orthogonal with dimensions $d-1$, $\tbinom d2$, $d$, $d(d-2)$ summing to $N(d)$, so they span $\TP$, and their irreducibility and pairwise non-isomorphism follow from the little-group criterion for the wreath product $B_d = \Z/2\Z \wr S_d$, applied to the characters of the sign subgroup $(\Z/2\Z)^d$ that occur in each summand.
\end{remark}

All Hessian eigenvalues quoted in this paper are with respect to the norm $\|\bdelta\|$ of Section~\ref{ss:CubeParam-vertex}. Equivalent invariant choices change each eigenvalue by a fixed positive constant, so the signs --- which are what the theorems require --- are independent of the choice.

Throughout the rest of the paper we work with a specific representative vector field $\bth$ on $\cube$ in each of the four blocks of~\eqref{e:irrep-decomp}:
\begin{itemize}
\item $R$-deformation: $\bth_R(\x) := \bigl((d-1)x_1, -x_2, \ldots, -x_{d}\bigr)$, a volume-preserving traceless stretch.
\item $S$-deformation (in axes $1, 2$): $\bth_S(\x) := (x_2, x_1, 0, \ldots, 0)$, a pairwise pure shear---a traceless stretch along the $\tfrac{1}{\sqrt2}(1, 1)$ diagonal of the $(x_1, x_2)$-plane---carrying $\cube$ to a rhombus $\times$ box; \emph{not} a simple shear, so the image is equilateral.
\item $T$-deformation (modulator $c = d$): $\bth_T(\x) := (4 x_1 x_{d}, 4 x_2 x_{d}, \ldots, 4 x_{d-1} x_{d}, 0)$.
\item $U$-deformation (axes $1, 2$, modulator $c = d$): $\bth_U(\x) := (4 x_1 x_{d}, -4 x_2 x_{d}, 0, \ldots, 0)$.
\end{itemize}
For a given deformation we call the axes along which $\bth$ has a nonzero component the \emph{active} axes, and the axes along which $\bth$ neither points nor varies the \emph{inactive} axes.  The inactive axes are $3, \ldots, d$ for $\bth_S$ and $3, \ldots, d-1$ for $\bth_U$; $\bth_R$ and $\bth_T$ have none.  The modulator axis $c$ of the $T$- and $U$-deformations is neither active nor inactive, and a facet $F_{k, \sigma}$ is called active, inactive, or modulator according to its normal axis $k$.

The factor $4 = (\tfrac{1}{2})^{-2}$ in $\bth_T, \bth_U$ normalizes them so that at each vertex $\v_{\bsigma} \in \{-\tfrac{1}{2}, \tfrac{1}{2}\}^d$ the nonzero components of $\bth$ have absolute value $1$. Direct calculation gives the squared $\ell^2$-norms (sums of squared components over all $2^d$ vertices):
\begin{equation}\label{e:norms}
\|\bth_R\|^2 \;=\; 2^{d-2}d(d-1), \quad
\|\bth_S\|^2 \;=\; 2^{d-1}, \quad
\|\bth_T\|^2 \;=\; 2^d(d-1), \quad
\|\bth_U\|^2 \;=\; 2^{d+1}.
\end{equation}

\subsection{Criticality and block-diagonality of the Hessian}\label{ss:Critical}

\begin{proposition}\label{p:Crit-Block}
Let $F$ be a $B_d$-invariant functional, smooth in a neighborhood of $\cube$ in the class $\calQ$ (e.g., $F = V$, $P$, $I$, $\blam$, or the auxiliary functional $f_d$ of~\eqref{e:fd-defn}; see Sections~\ref{s:Iso}--\ref{s:LNshape}).  Then, in the linear coordinates $\bdelta \in \TQ$,
\begin{enumerate}
\item[(i)] $\nabla F(\cube) = \bzero$;
\item[(ii)] the restriction of the Hessian $\Hess F(\cube)$ to the $B_d$-invariant subspace $\TP$ is, in symmetry-adapted coordinates --- the decomposition~\eqref{e:irrep-decomp} into irreducible $B_d$-representations --- block-diagonal, and on each irreducible block it is a scalar multiple of the identity;
\item[(iii)] the cube is a critical point of the restriction of $F$ to the manifold $\P$ of~\eqref{e:M-defn}, with Hessian $\Hess F(\cube)|_{\TP}$; in particular, for every smooth curve $\bdelta(t)$ in $\TQ$ with $\bdelta(0) = \bzero$ and $\bdelta'(0) = \bth \in \TP$, one has $\tfrac{d^2}{dt^2}\big|_{t=0} F(\Omega_{\bdelta(t)}) = \Hess F(\cube)[\bth, \bth]$.
\end{enumerate}
\end{proposition}

\begin{proof}
(i) By $B_d$-invariance, $F \circ g = F$ for all $g \in B_d$.  Differentiating at $\cube$ (which is fixed by every $g \in B_d$), $(Dg)^T \nabla F(\cube) = \nabla F(\cube)$ for every $g$: the gradient is a $B_d$-invariant vector in $\TQ \subset (\R^d)^{2^d}$.  We claim that the only $B_d$-invariant vertex displacements are the dilations, which lie in $\sm$ and are therefore orthogonal to $\TQ \supset \TP$; this gives $\nabla F(\cube) = \bzero$.  Let $\bdelta$ be invariant and recall the action of $B_d$ on vertex displacements from Section~\ref{ss:Bd-action}.  Invariance under the sign flip of an axis $j \ne k$ gives $\bdelta_{\bsigma}[k] = \bdelta_{\bsigma'}[k]$ whenever $\bsigma'$ differs from $\bsigma$ in the $j$th sign only, so $\bdelta_{\bsigma}[k]$ depends on $\bsigma$ through $\sigma_k$ alone; invariance under the flip of axis $k$ makes it odd in $\sigma_k$, so $\bdelta_{\bsigma}[k] = a_k \sigma_k$ for some $a_k \in \R$; and invariance under the transposition of two axes $k, l$ gives $a_k = a_l$.  Hence $\bdelta_{\bsigma} = a\,\bsigma = 2a\,\v_{\bsigma}$, a dilation, as claimed.

(ii) The Hessian $\Hess F(\cube)$ is a $B_d$-invariant symmetric bilinear form on $\TQ$, hence on the invariant subspace $\TP$.  By Theorem~\ref{t:irrep-decomp} the four blocks of~\eqref{e:irrep-decomp} are pairwise non-isomorphic irreducible representations, so by Schur's lemma an invariant form vanishes between distinct blocks and is a scalar multiple of the identity on each.

(iii) By the chain rule, $\tfrac{d^2}{dt^2}\big|_{t=0} F(\Omega_{\bdelta(t)}) = \Hess F(\cube)[\bdelta'(0), \bdelta'(0)] + \nabla F(\cube)\cdot\bdelta''(0)$, and the second term vanishes by (i).  Applied to curves in $\P$, whose velocities at $\bzero$ exhaust $\TP$, this identifies the Hessian of $F|_{\P}$ at $\cube$ with $\Hess F(\cube)|_{\TP}$.
\end{proof}

For each functional $F \in \{\blam, I, f_d\}$, Proposition~\ref{p:Crit-Block} reduces the Hessian, in symmetry-adapted coordinates, to four scalar eigenvalues $\nu^F_R, \nu^F_S, \nu^F_T, \nu^F_U$ on the blocks of~\eqref{e:irrep-decomp}, which completely characterize the Hessian on $\P$. Each is computed along the straight curve $t \mapsto \Omega_{t\bth_X}$ through the corresponding representative field of Section~\ref{ss:Reps},
\begin{equation}\label{e:nu-defn}
\nu^F_X \;=\; \frac{F''(\cube)[\bth_X, \bth_X]}{\|\bth_X\|^2}, \qquad F''(\cube)[\bth_X, \bth_X] \;:=\; \frac{d^2}{dt^2}\Big|_{t=0}F(\Omega_{t\bth_X}),
\end{equation}
which is independent of the rescaling of $\bth_X$; by Proposition~\ref{p:Crit-Block}(iii) this is the Hessian scalar whether or not the curve stays in $\P$ --- in general it does not, since the coplanarity constraints are quadratic in $\bdelta$, but for the four representative fields it does, since along each of them the deformed body is cut out by $2d$ linear inequalities (a box, a rhombic prism, a frustum, and a twisted frustum).

\subsection{The rectangular deformation in closed form}\label{ss:R-box}

For the $R$-deformation every quantity in this paper is elementary, and we record the results now: they serve as a sanity check on Sections~\ref{s:LDshape} and~\ref{s:LNshape}, where the $R$-rows of the tables are recomputed by the same methods as the other three deformations (in Section~\ref{s:Iso} the $R$-deformation computation is the box computation itself).  Since $\bth_R$ is affine and multilinear interpolation reproduces affine maps, the body $\Omega_t := \Omega_{t\bth_R} = (\mathrm{id} + t\bth_R)(\cube)$ is the axis-aligned box with side lengths $\ell_1 = 1 + (d-1)t,  \ell_2 = \cdots = \ell_d = 1 - t$, for which volume, surface area and the low spectrum are explicit:
\begin{equation}\label{e:R-box}
V(t) = \prod_{k} \ell_k, \quad
P(t) = 2\sum_{k}\prod_{j \ne k} \ell_j, \quad
\lambda_1(\Omega_t) = \pi^2 \sum_{k} \ell_k^{-2}, \quad
\{\mu_1, \ldots, \mu_d\}(\Omega_t) = \bigl\{\pi^2 \ell_k^{-2}\bigr\}_{k=1}^{d},
\end{equation}
the last for $|t|$ small, where $\pi^2 \ell_k^{-2} < 2\pi^2 \min_j \ell_j^{-2}$ so that these are the $d$ lowest nontrivial Neumann eigenvalues.

\begin{proposition}[The $R$-deformation]\label{p:R-box}
For every $d \ge 2$,
\[
\nu^I_R \;=\; 4\,d^d, \qquad
\nu^{\blam}_R \;=\; \frac{16\pi^2}{2^d}, \qquad
\nu^{f_d}_R \;=\; \frac{16}{d\cdot 2^d\,\pi^2},
\]
with $f_d$ the auxiliary functional~\eqref{e:fd-defn}.  Moreover $t \mapsto \bmu(\Omega_t)$ has one-sided derivatives $-2(d-1)\pi^2$ at $t = 0^+$ and $-2\pi^2$ at $t = 0^-$ (that is, along $-\bth_R$).
\end{proposition}

\begin{proof}
Expanding~\eqref{e:R-box} to second order,
\[
\begin{gathered}
V(t) = 1 - \tfrac{d(d-1)}{2}\,t^2 + O(t^3), \qquad
P(t) = 2d - d(d-1)(d-2)\,t^2 + O(t^3), \\
\lambda_1(\Omega_t) = d\pi^2 + 3d(d-1)\pi^2\, t^2 + O(t^3),
\end{gathered}
\]
the linear terms cancelling because $\sum_k \log \ell_k$, $\sum_k \ell_k^{-2}$ and $\sum_k \prod_{j \ne k} \ell_j$ are all stationary at $t = 0$ (which is Proposition~\ref{p:Crit-Block}(i) for $V$, $\lambda_1$ and $P$).  Hence $I = P^d V^{1-d} = (2d)^d\bigl[1 + \tfrac{d(d-1)}{2}t^2\bigr] + O(t^3)$, $\blam = V^{2/d}\lambda_1 = d\pi^2 + 2d(d-1)\pi^2 t^2 + O(t^3)$, and, from $\sum_k \mu_k^{-1} = \pi^{-2}\bigl(\ell_1^2 + (d-1)\ell_2^2\bigr) = \pi^{-2}\bigl(d + d(d-1)t^2\bigr)$ and $V^{-2/d} = 1 + (d-1)t^2 + O(t^3)$, $f_d(\Omega_t) = \tfrac{1}{d}\,V^{-2/d}\sum_{k}\tfrac{1}{\mu_k} = \tfrac{1}{\pi^2}\bigl[1 + 2(d-1)t^2\bigr] + O(t^3)$. Dividing the three second derivatives $d(d-1)(2d)^d$, $4d(d-1)\pi^2$ and $4(d-1)/\pi^2$ by $\|\bth_R\|^2 = 2^{d-2}d(d-1)$ from~\eqref{e:norms} gives the three eigenvalues by~\eqref{e:nu-defn}.  Finally $\mu_1(\Omega_t) = \pi^2 \max_k \ell_k^{-2}$ equals $\pi^2(1 + (d-1)t)^{-2}$ for $t > 0$ and $\pi^2(1-t)^{-2}$ for $t < 0$, and $V = 1 + O(t^2)$, which gives the one-sided derivatives of $\bmu$.
\end{proof}

%%%%%%%%%%%%%%%%%%%%%%%%%%%%%%%%%%
%%%%%%%%%%%%%%%%%%%%%%%%%%%%%%%%%%
\section{The isoperimetric ratio}\label{s:Iso}

We compute the Hessian of $I(\Omega) = \tfrac{P(\Omega)^d}{|\Omega|^{d-1}}$ at $\cube$, one representative deformation at a time, along the fields of Section~\ref{ss:Reps}.  Along each of these perturbed bodies is a polytope whose volume and surface area are elementary in $t$, so the computation is direct; we first record the general shape-derivative formulas at the cube, which are used for the non-polytopal deformation of Section~\ref{s:TopMode}.

\subsection{Shape derivatives at the cube}\label{ss:Shape-prelim}

We adopt the convention that a smooth velocity field $\bth$ on $\cube$ generates a family of perturbed domains $\Omega_t := (\mathrm{id} + t\bth)(\cube)$ for small $t \in \R$, and write $V(t) := |\Omega_t|$ and $P(t) := P(\Omega_t)$ for the volume and perimeter along the flow.  Each direction of the symmetry-adapted parameter space gives such a field via~\eqref{e:Phi-delta}.

\begin{lemma}\label{l:facet-derivs}
For $\cube$ with facets $F_{a, \sigma} := \{\x \in \cube \colon x_a = \tfrac{\sigma}{2}\}$ ($a \in \{1, \ldots, d\}$, $\sigma \in \{-1, +1\}$), and any smooth $\bth$,
\begin{subequations}
\label{e:derivV,P}
\begin{align}
\label{e:derivV,Pa}
V'(0)[\bth] &= \int_{\cube} \nabla \cdot \bth\, d\x, \\
\label{e:derivV,Pb}
V''(0)[\bth, \bth] &= \int_{\cube} \Bigl[(\nabla \cdot \bth)^2 - \tr\bigl((D\bth)^2\bigr)\Bigr]\, d\x, \\
\label{e:derivV,Pc}
P'(0)[\bth] &= \sum_{a, \sigma} \int_{F_{a, \sigma}} (\nabla \cdot \bth - \partial_a \bth_a)\, d\sigma_F, \\
\label{e:derivV,Pd}
P''(0)[\bth, \bth] &= \sum_{a, \sigma} \int_{F_{a, \sigma}} \Bigl[\tr(B_F) + \tfrac{1}{4} \tr(A_F)^2 - \tfrac{1}{2} \tr(A_F^2)\Bigr]\, d\sigma_F,
\end{align}
\end{subequations}
where on the facet $F_{a, \sigma}$, $A_F = (\partial_i \bth_j + \partial_j \bth_i)_{i,j \ne a}$ is the tangential symmetric strain and $B_F = ((\partial_k \bth) \cdot (\partial_k \bth))_{k \ne a}$ is the tangential gradient gram matrix. Moreover, 
$V'(0)[\bth] = P'(0)[\bth] = 0$ for every direction in $\TQ$ and
\begin{equation}\label{e:Ipp-formula}
I''(0)[\bth, \bth] \;=\; d(2d)^{d-1} P''(0)[\bth, \bth] \;-\; (d-1)(2d)^d\, V''(0)[\bth, \bth].
\end{equation}
\end{lemma}

\begin{proof}
The volume and perimeter formulas are classical (see~\cite[Ch.~5]{HenrotPierre2018} and \cite{DiazAvalosLaurain2025}). 
The criticality follows from Proposition~\ref{p:Crit-Block}(i) applied to $V$ and $P$. 
The composite formula~\eqref{e:Ipp-formula} follows from the chain rule applied to $I = P^d V^{-(d-1)}$ at $V_0 = 1$, $P_0 = 2d$, using $V'(0) = P'(0) = 0$.
\end{proof}

\subsection{Hessian of the isoperimetric ratio}\label{ss:IsoHessian}

We compute $I''(0)$ along each representative field $\bth_X$, $X \in \{R, S, T, U\}$, of Section~\ref{ss:Reps}.  Throughout, $V_0 = 1$, $P_0 = 2d$, $I_0 = (2d)^d$.

\begin{proposition}\label{p:Iso-irrep}
Let $d \ge 2$ and $X \in \{R, S, T, U\}$. The Jacobian $B_X := D\bth_X$ of the representative field (only nonzero entries listed), its traces, the second shape derivatives $V''(0) := V''(0)[\bth_X, \bth_X]$ and $P''(0) := P''(0)[\bth_X, \bth_X]$ of volume and perimeter at $\cube$, and the resulting Hessian eigenvalue of $I$ are
{\small\renewcommand{\arraystretch}{1.8}\setlength{\arraycolsep}{4pt}
\[
\begin{array}{c|l|c|c|c|c|c}
X & \multicolumn{1}{c|}{B_X} & \tr B_X & \tr(B_X^2) & V''(0) & P''(0) & \nu^I_X \\
\hline
R   & \mathrm{diag}(d{-}1, -1, \ldots, -1) & 0 & d(d-1) & -d(d-1)               & -2 d(d-1)(d-2)        & 4\, d^d \\
S   & B_{12} = B_{21} = 1 & 0 & 2 & -2                    & 12 - 4d               & 8\, d^d \\
T & \substack{B_{ii} = 4x_d,\ B_{id} = 4x_i \\ (i < d)} & 4(d{-}1)x_d & 16(d{-}1)x_d^2 & \tfrac{4(d-1)(d-2)}{3} & \tfrac{8(d-1)(d^2 - 2d + 3)}{3} & \tfrac{4(d+1)}{3}\, d^d \\
U & \substack{B_{11} = 4x_d,\ B_{22} = -4x_d, \\ B_{1d} = 4x_1,\ B_{2d} = -4x_2} & 0 & 32x_d^2 & -\tfrac{8}{3}         & \tfrac{16(3-d)}{3}    & \tfrac{8}{3}\, d^d \quad (d \ge 3)
\end{array}
\]}
The eigenvalues are as stated in~\eqref{e:Iso-eigs}; all four are strictly positive.
\end{proposition}

\begin{proof}
The Jacobians $B = B_X$ and their traces follow by differentiating the fields of Section~\ref{ss:Reps}. By~\eqref{e:derivV,Pb}, $V''(0) = \int_{\cube}\bigl[(\tr B)^2 - \tr(B^2)\bigr]$, and on the facet $F_{a,\sigma}$ write $A_F = (B_{ij} + B_{ji})_{i,j \ne a}$ and $\tr B_F = \sum_{k \ne a}|\partial_k\bth_X|^2$, so that $P''(0)$ is the sum over the $2d$ facets of $\int_{F_{a,\sigma}}\bigl[\tr B_F + \tfrac14(\tr A_F)^2 - \tfrac12\tr(A_F^2)\bigr]$ by~\eqref{e:derivV,Pd}. Both $V$ and $P$ are stationary at $t = 0$ by Proposition~\ref{p:Crit-Block}(i), so~\eqref{e:Ipp-formula} applies, and the eigenvalues follow from~\eqref{e:nu-defn} with the norms~\eqref{e:norms}. Below, $\int_{-1/2}^{1/2}x^2\,dx = \tfrac{1}{12}$ and $x_a^2 = \tfrac14$ on $F_{a,\sigma}$.

\smallskip\noindent\emph{$R$-deformation.} $V''(0) = -d(d-1)$. On $F_{1,\sigma}$, $A_F = -2I_{d-1}$ and $\tr B_F = d-1$, giving the integrand $(d-1)(d-2)$; on the $2(d-1)$ facets with $a \ge 2$, $A_F = 2\,\mathrm{diag}(d-1, -1, \ldots, -1)$ and $\tr B_F = (d-1)^2 + (d-2)$, giving $-(d-2)(d+1)$. Hence $P''(0) = -2d(d-1)(d-2)$.

\smallskip\noindent\emph{$S$-deformation.} $V''(0) = -2$. On the $2(d-2)$ facets with $a \ge 3$, $(A_F)_{12} = (A_F)_{21} = 2$ is the only nonzero entry and $\tr B_F = 2$, giving the integrand $-2$; on $F_{1,\sigma}$ and $F_{2,\sigma}$, $A_F = 0$ and $\tr B_F = 1$. Hence $P''(0) = 12 - 4d$.

\smallskip\noindent\emph{$T$-deformation.} $V''(0) = 16(d-1)(d-2)\int_{\cube}x_d^2 = \tfrac{4(d-1)(d-2)}{3}$. On the modulator facets $F_{d,\sigma}$, $A_F = 8x_d I_{d-1}$ and $\tr B_F = 16(d-1)x_d^2$, giving the integrand $16(d-1)(d-2)x_d^2 = 4(d-1)(d-2)$. On the $2(d-1)$ active facets $F_{i,\sigma}$, $i < d$, $A_F$ has entries $(A_F)_{jj} = 8x_d$ and $(A_F)_{jd} = (A_F)_{dj} = 4x_j$ for $j \ne i$, $j < d$, and $\tr B_F = 16(d-2)x_d^2 + 16\sum_{j<d}x_j^2$; the integrand is $16(d-2)(d-3)x_d^2 + 16x_i^2$, which integrates to $\tfrac{4(d-2)(d-3)}{3} + 4$. Hence $P''(0) = \tfrac{8(d-1)(d^2-2d+3)}{3}$.

\smallskip\noindent\emph{$U$-deformation ($d \ge 3$).} $V''(0) = -32\int_{\cube}x_d^2 = -\tfrac{8}{3}$. On the modulator facets $F_{d,\sigma}$, $A_F = \mathrm{diag}(8x_d, -8x_d, 0, \ldots, 0)$ and $\tr B_F = 32x_d^2$, giving the integrand $-32x_d^2 = -8$. On $F_{1,\sigma}$, the nonzero entries of $A_F$ are $(A_F)_{22} = -8x_d$ and $(A_F)_{2d} = (A_F)_{d2} = -4x_2$, and $\tr B_F = 16(x_1^2 + x_2^2 + x_d^2)$; the integrand is $16x_1^2 = 4$, and likewise on $F_{2,\sigma}$. On the $2(d-3)$ inactive facets, $A_F$ has entries $8x_d, -8x_d$ on the diagonal and $4x_1, -4x_2$ in positions $(1,d), (2,d)$ and their transposes, and $\tr B_F = 32x_d^2 + 16(x_1^2 + x_2^2)$; the integrand is $-32x_d^2$, which integrates to $-\tfrac83$. Hence $P''(0) = \tfrac{16(3-d)}{3}$.
\end{proof}

%%%%%%%%%%%%%%%%%%%%%%%%%%%%%%%%%%
%%%%%%%%%%%%%%%%%%%%%%%%%%%%%%%%%%
\section{The Laplace--Dirichlet eigenvalue}\label{s:LDshape}

Consider the volume-normalized eigenvalue $\blam = |\Omega|^{2/d}\lambda_1$ and any direction with $V'(0)[\bth] = 0$; this holds for every direction $\bth \in \TQ$, in particular for the four representative fields and for the corner perturbation $\bth^*$ of Section~\ref{s:TopMode}, by Proposition~\ref{p:Crit-Block}(i) applied to $V$. The product rule for $\blam = V^{2/d}\lambda$ then gives, at $\cube$
\begin{equation}\label{e:lam-bar-formula}
\blam''(0)[\bth, \bth] \;=\; \lambda''(0)[\bth, \bth] \;+\; \tfrac{2}{d}\,\lambda_0\, V''(0)[\bth, \bth]
\;=\; \lambda''(0)[\bth, \bth] \;+\; 2\pi^2\, V''(0)[\bth, \bth],
\end{equation}
where we used that  $V_0 = 1$ and $\lambda_0 := \lambda_1(\cube) = d\pi^2$, and the cross term vanishes since $\lambda'(0) = V'(0) = 0$. 
Here $\lambda''(0)[\bth, \bth]$ is the second shape derivative appearing in~\eqref{e:dL1-BB}; rather than evaluating~\eqref{e:dL1-BB} directly, we 
pull back the eigenvalue problem on $\Omega_t = \Phi_t(\cube)$ to the fixed cube and expand to second order in $t$.  This produces an explicit formula for $\lambda''(0)$ as the sum of a \emph{polynomial part} and a \emph{Fourier-sum part}.

The Dirichlet eigenbasis on $\cube$ is indexed by tuples $\bfn = (n_1, \ldots, n_{d}) \in \Z_{\ge 1}^d$:
\begin{equation}\label{e:phi-defn}
\Psi_{\bfn}(\x) \;=\; 2^{d/2}\prod_{k=1}^{d} \phi_{n_k}(x_k), \qquad
\phi_n(x) \;:=\; \sin\!\bigl(n\pi(x + \tfrac12)\bigr),
\end{equation}
with eigenvalues $\pi^2 |\bfn|^2$, where $|\bfn|^2 = \sum_k n_k^2$; the $\Psi_\bfn$ are $L^2(\cube)$-orthonormal. The ground state corresponds to $\mathbf{1} := (1, \ldots, 1)$: since $\phi_1(x) = \cos(\pi x)$, $\Psi_{\mathbf{1}}(\x) = 2^{d/2} \prod_{k} \phi_1(x_k) = 2^{d/2} \prod_{k} \cos(\pi x_k)$, and $\partial_{x_k} \Psi_{\mathbf{1}} = 2^{d/2}\,\phi_1'(x_k)\prod_{j \ne k}\phi_1(x_j)$.

We record here some elementary identities involving inner products of $\phi_n$, $\phi_n'$ and their moments that will be used below. Since $\phi_n(-x) = (-1)^{n+1}\phi_n(x)$, each $\phi_n$ is even in $x$ for $n$ odd and odd for $n$ even, with $\phi_n'$ of the opposite parity. Writing
\begin{equation}\label{e:weighted-pairing}
  \langle f, g\rangle_k \;:=\; \int_{-1/2}^{1/2} x^k f(x) g(x)\, dx, \qquad k = 0, 1, 2,
\end{equation}
the pairing~\eqref{e:weighted-pairing} $\langle \phi_m^{(a)}, \phi_n^{(b)}\rangle_k$ with $a, b \in \{0,1\}$ therefore vanishes unless $k + m + n + a + b$ is even.  In the remaining cases, for $m \ne n$,
\begin{equation}\label{e:mom-offdiag}
\begin{gathered}
  \langle \phi_m, \phi_n\rangle_0 = 0, \qquad
  \langle \phi_m', \phi_n\rangle_0 = \frac{2mn}{n^2-m^2}, \qquad
  \langle \phi_m', \phi_n'\rangle_0 = 0, \\[2pt]
  \langle \phi_m, \phi_n\rangle_1 = \frac{-4mn}{\pi^2(m^2-n^2)^2}, \quad
  \langle \phi_m', \phi_n\rangle_1 = \frac{-mn}{n^2-m^2}, \quad
  \langle \phi_m', \phi_n'\rangle_1 = \frac{-2mn(m^2+n^2)}{(m^2-n^2)^2},
\end{gathered}
\end{equation}
while for $m = n$,
\begin{equation}\label{e:mom-diag}
\begin{gathered}
  \langle \phi_m, \phi_m\rangle_0 = \tfrac12, \qquad
  \langle \phi_m', \phi_m'\rangle_0 = \tfrac{m^2\pi^2}{2}, \qquad
  \langle \phi_m', \phi_m\rangle_1 = -\tfrac14, \\[2pt]
  \langle \phi_m, \phi_m\rangle_2 = \tfrac{1}{24} - \tfrac{1}{4m^2\pi^2}, \qquad
  \langle \phi_m', \phi_m'\rangle_2 = \tfrac{m^2\pi^2}{24} + \tfrac14 .
\end{gathered}
\end{equation}
All of these follow from trigonometric identities.

\subsection{The pulled-back eigenvalue equation and its Taylor expansion}\label{sss:pullback-PDE}\label{ss:Material}

Let $\Phi_t = \mathrm{id} + t\bth$ for $\bth$ smooth and $t$ in a neighborhood of $0$, and set $u(t, \x) := v_t(\Phi_t(\x))$, where $v_t$ is the positive, $L^2$-normalized Dirichlet ground state on $\Omega_t = \Phi_t(\cube)$. A standard change of variables (cf.~\cite[Ch.~5]{HenrotPierre2018}) converts $-\Delta_y v_t = \lambda(t) v_t$ on $\Omega_t$ into a weighted self-adjoint eigenvalue problem on the fixed cube:
\begin{align*}
&-\nabla \cdot \bigl(A(t, \x)\, \nabla u\bigr) \;=\; \lambda(t)\, \rho(t, \x)\, u  &&\text{in } \cube, \\ 
&u = 0 && \text{on } \partial\cube,
\end{align*}
where $\rho(t, \x) := \det D\Phi_t(\x)$  and  $A(t, \x) := \rho(t, \x)\, D\Phi_t(\x)^{-1}\, D\Phi_t(\x)^{-T}$. For our four representative fields, $D\Phi_t$ has positive determinant near $t = 0$, so we may drop the absolute value.  Both $\rho$ and $A$ are real-analytic in $t$, and the family of operators $L(t) := -\nabla \cdot(A(t)\nabla\cdot)$ on the fixed Sobolev space $H^1_0(\cube)$ is of Kato type~(B)~\cite[Ch.~VII, \S 4]{Kato1995}, so the lowest eigenvalue $\lambda(t)$ and (a smooth choice of) eigenfunction $u(t)$ are real-analytic in $t$.  The same applies to the multilinear family $\Phi_{\bdelta}$ of~\eqref{e:Phi-delta}, which is polynomial in $\bdelta$: $\lambda_1(\Omega_{\bdelta})$ is real-analytic in $\bdelta \in \TQ$ near $\bzero$, and so are $V$ and $P$, hence $I$ and $\blam$ are smooth on $\calQ$, as used in Proposition~\ref{p:Crit-Block}.

We Taylor-expand around $t = 0$, where $A_0 = I$ and $\rho_0 = 1$: $A = I + t A_1 + t^2 A_2 + O(t^3), \  \rho = 1 + t \rho_1 + t^2 \rho_2 + O(t^3), \  u = u_0 + t u_1 + t^2 u_2 + O(t^3)$, and $\lambda(t) = \lambda_0 + t \lambda'(0) + \tfrac{t^2}{2} \lambda''(0) + O(t^3)$, the zeroth-order terms being the unperturbed ground state and eigenvalue, $u_0 = \Psi_{\mathbf{1}}$ and $\lambda_0 = d\pi^2$.  (We reserve subscripts on $\lambda$ for the eigenvalue index, as in $\lambda_1(\Omega)$; derivatives along the flow are written $\lambda'(0)$, $\lambda''(0)$.)  Writing $L_k u := -\nabla \cdot (A_k \nabla u)$ for $k \ge 1$, with $L_0 = -\Delta$, the order-$t^0$ equation $L_0 u_0 = \lambda_0 u_0$ is the unperturbed problem, and the order-$t^1$ and order-$t^2$ equations are
\begin{subequations}
\begin{align}
(L_0 - \lambda_0)\, u_1 \;&=\; \lambda'(0)\, u_0 \,-\, h,
  \qquad\qquad \text{where } h \;:=\; L_1 u_0 \,-\, \lambda_0 \rho_1 u_0, \label{e:h-defn}\\
(L_0 - \lambda_0)\, u_2 \;&=\; \tfrac{1}{2}\lambda''(0)\, u_0 \,+\, \lambda_0 \rho_2 u_0 \,+\, \lambda'(0) \bigl(u_1 + \rho_1 u_0\bigr) \,+\, \lambda_0 \rho_1 u_1 \,-\, L_1 u_1 \,-\, L_2 u_0.  \label{e:order2}
\end{align}
\end{subequations}
For each $\bth_X$, $\lambda'(0) = 0$ at the cube: the first Dirichlet eigenvalue is $B_d$-invariant and smooth near $\cube$, so Proposition~\ref{p:Crit-Block}(i) applies to it directly.  The order-$1$ equation~\eqref{e:h-defn} therefore reads $(L_0 - \lambda_0) u_1 = -h$, and the $\lambda'(0)$ term in~\eqref{e:order2} drops.  We impose the normalization $\langle u_1, u_0\rangle = 0$.  Taking $L^2$-inner products of \eqref{e:h-defn}--\eqref{e:order2} with $u_0$ kills both left-hand sides (since $L_0 - \lambda_0$ is self-adjoint and annihilates $u_0$) and produces the Fredholm solvability conditions
\begin{subequations}
\begin{align}
\lambda'(0) \;&=\; \langle h, u_0\rangle \;=\; \langle L_1 u_0, u_0\rangle \,-\, \lambda_0 \langle \rho_1 u_0, u_0\rangle, \label{e:lam1-formula}\\
\tfrac{1}{2}\lambda''(0) \;&=\; \langle L_2 u_0, u_0\rangle \,-\, \lambda_0 \langle \rho_2 u_0, u_0\rangle \,+\, \langle L_1 u_1, u_0\rangle \,-\, \lambda_0 \langle \rho_1 u_1, u_0\rangle. \label{e:lam2-formula}
\end{align}
\end{subequations}
Since $L_1$ is self-adjoint, the last two terms of~\eqref{e:lam2-formula} combine into a single inner product, $\langle u_1, h\rangle$. Expanding $u_1 = -\sum_{\bfn \ne \mathbf{1}} \tfrac{\langle h, \Psi_\bfn\rangle}{\pi^2(|\bfn|^2 - d)} \Psi_\bfn$ in the Dirichlet basis gives $\langle u_1, h\rangle = -\sum_{\bfn \ne \mathbf{1}} \tfrac{|\langle h, \Psi_\bfn\rangle|^2}{\pi^2(|\bfn|^2 - d)}$.  We obtain
\begin{equation}\label{e:lam2-final}
\lambda''(0) \;=\; \PiD_X(d) \,+\, \SiD_X(d),
\qquad
\begin{aligned}
\PiD_X(d) \,&:=\, 2\langle L_2 u_0, u_0\rangle \,-\, 2\lambda_0 \langle \rho_2 u_0, u_0\rangle, \\
\SiD_X(d) \,&:=\, -2\sum_{\bfn \ne \mathbf{1}} \frac{H_{\bfn}^2}{\pi^2(|\bfn|^2 - d)},
\end{aligned}
\end{equation}
where $H_{\bfn} := \langle h, \Psi_\bfn\rangle$, and the denominators $\pi^2(|\bfn|^2 - d) = \pi^2|\bfn|^2 - \lambda_0$ are the spectral gaps between $\Psi_\bfn$ and the ground state, referred to simply as gaps below.  We call $\PiD_X(d)$ the \emph{polynomial part} of $\lambda''(0)[\bth_X, \bth_X]$ and $\SiD_X(d)$ the \emph{Fourier-sum part}, analyzed in Sections~\ref{sss:Polynomial-contributions} and~\ref{sss:Dirichlet-assembly}.  The Hessian eigenvalue is then
\begin{equation}\label{e:nu-formula}
\nu^{\blam}_X \;=\; \frac{\lambda''(0)[\bth_X, \bth_X] \,+\, 2\pi^2\, V''(0)[\bth_X, \bth_X]}{\|\bth_X\|^2},
\end{equation}
where $V''(0)[\bth_X, \bth_X] = 2\int_{\cube}\rho_2(\x)\, d\x$.

The coefficients $A_1, A_2, \rho_1, \rho_2$ of this expansion are computed in the following lemma.

\begin{lemma}[Pullback coefficients]\label{l:coeff-expansion}
Let $\Phi_t = \mathrm{id} + t\bth$ and write $B := D\bth$. Then $\rho(t) = \det D\Phi_t$ and $A(t) = \rho\, D\Phi_t^{-1} D\Phi_t^{-T}$ have Taylor coefficients
\begin{equation}\label{e:coeff-universal}
\begin{aligned}
\rho_1 \;&=\; \tr B \;=\; \nabla \cdot \bth, \qquad &
\rho_2 \;&=\; \tfrac{1}{2}\bigl[(\tr B)^2 - \tr(B^2)\bigr], \\
A_1 \;&=\; (\tr B)\, I - B - B^{T}, \qquad &
A_2 \;&=\; \rho_2\, I - (\tr B)\bigl(B + B^{T}\bigr) + B^2 + (B^{T})^2 + B B^{T}.
\end{aligned}
\end{equation}
For the four representative fields $\bth_X$, $X \in \{R, S, T, U\}$, the Jacobians $B_X$ and their traces are listed in Proposition~\ref{p:Iso-irrep}, so $\rho_1$ and $\rho_2$ can be read off there; the remaining data are
{\small\renewcommand{\arraystretch}{2.0}
\[
\begin{array}{c|l|c}
X & \multicolumn{1}{c|}{A_1} & \det D\Phi_t \\
\hline
R & \mathrm{diag}(-2(d{-}1), 2, \ldots, 2) & (1 + (d-1)t)(1-t)^{d-1} \\
S & (A_1)_{12} = -2 & 1 - t^2 \\
T & \substack{(A_1)_{ii} = 4(d{-}3)x_d,\ (A_1)_{id} = -4x_i, \\ (A_1)_{dd} = 4(d{-}1)x_d \quad (i < d)} & (1 + 4tx_d)^{d-1} \\
U & \substack{(A_1)_{11} = -8x_d,\ (A_1)_{22} = 8x_d, \\ (A_1)_{1d} = -4x_1,\ (A_1)_{2d} = 4x_2} & 1 - 16t^2 x_d^2 \\
\end{array}
\]}
where only the nonzero entries are listed and $A_1$ is symmetric.
\end{lemma}

\begin{proof}
Expanding $\det(I + tB) = 1 + t\,\tr B + \tfrac{t^2}{2}[(\tr B)^2 - \tr(B^2)] + O(t^3)$ gives the $\rho$ coefficients, and $(I + tB)^{-1} = I - tB + t^2 B^2 + O(t^3)$ gives
$D\Phi_t^{-1} D\Phi_t^{-T} = I - t(B + B^T) + t^2\bigl(B^2 + (B^T)^2 + BB^T\bigr) + O(t^3)$; multiplying by $\rho(t)$ yields~\eqref{e:coeff-universal}. The table follows by substituting the $B_X$ of Proposition~\ref{p:Iso-irrep}.
\end{proof}

\subsection{The polynomial part \texorpdfstring{$\PiD_X(d)$}{Pi\_X(d)}}\label{sss:Polynomial-contributions}

We now evaluate $\PiD_X(d)$ in \eqref{e:lam2-final} for each $X \in \{R, S, T, U\}$.

\begin{lemma}\label{l:poly-table}
For the four representative fields $\bth_X$, $X \in \{R, S, T, U\}$,
\[
\renewcommand{\arraystretch}{1.5}
\setlength{\arraycolsep}{10pt}
\begin{array}{c|cccc}
X & R & S & T & U \\
\hline
\PiD_X(d) & 6\pi^2 d(d-1) & 12\pi^2 & 40(d-1)(d-2) - \tfrac{16(d-1)(d-3)\pi^2}{3} & -80 + \tfrac{64\pi^2}{3}
\end{array}
\]
\end{lemma}

\begin{proof}
By integration by parts, using $u_0 = 0$ on $\partial\cube$, the two inner products in $\PiD_X(d)$ are $\langle L_2 u_0, u_0\rangle = \int_{\cube} A_2\, \nabla u_0\cdot\nabla u_0\,d\x$ and $\langle \rho_2 u_0, u_0\rangle = \int_{\cube} \rho_2\, u_0^2\, d\x$. Substituting $B$ from the table of Lemma~\ref{l:coeff-expansion} into~\eqref{e:coeff-universal} makes $\rho_2$ and every entry $(A_2)_{ij}$ a polynomial of degree $\le 2$ in $\x$, so both integrands are products, over the axes, of $\phi_1$ and $\phi_1'$ factors multiplied by a polynomial in $\x$, and every axis integral is one of the diagonal moments~\eqref{e:mom-diag} with $m = 1$. It remains to carry this out mode by mode.

Since $u_0 = 2^{d/2}\prod_k \phi_1(x_k)$ and $\partial_k u_0 = 2^{d/2}\phi_1'(x_k)\prod_{j \ne k}\phi_1(x_j)$, the diagonal moments~\eqref{e:mom-diag} with $m = 1$ give, for distinct axes $k, l$ and $a \ne k$,
\begin{equation}\label{e:P-moments}
\begin{gathered}
\int_{\cube} (\partial_k u_0)^2 = \pi^2, \qquad
\int_{\cube} x_a^2\, (\partial_k u_0)^2 = \tfrac{\pi^2}{12} - \tfrac12, \qquad
\int_{\cube} x_k^2\, (\partial_k u_0)^2 = \tfrac{\pi^2}{12} + \tfrac12, \\
\int_{\cube} x_k x_l\, \partial_k u_0\, \partial_l u_0 = \tfrac14, \qquad
\int_{\cube} u_0^2 = 1, \qquad
\int_{\cube} x_a^2\, u_0^2 = \tfrac{1}{12} - \tfrac{1}{2\pi^2},
\end{gathered}
\end{equation}
and every integral below is a combination of these.  We write
\[
\langle L_2 u_0, u_0\rangle \;=\; \sum_{k, l}\int_{\cube}(A_2)_{kl}\,\partial_k u_0\,\partial_l u_0\,d\x .
\]

For $R$ and $S$, $B$ is constant and symmetric with $\tr B = 0$, so $A_2 = \rho_2 I + 3B^2$ is constant and diagonal, with $\tr A_2 = d\rho_2 + 3\tr(B^2)$ equal to $-\tfrac{d^2(d-1)}{2} + 3d(d-1)$ for $R$ and to $6 - d$ for $S$.  Hence $\langle L_2 u_0, u_0\rangle = \pi^2 \tr A_2$ and $\langle \rho_2 u_0, u_0\rangle = \rho_2$, which give $\PiD_R(d) = 2\pi^2 d(d-1)\bigl(3 - \tfrac{d}{2}\bigr) + 2 d\pi^2 \cdot \tfrac{d(d-1)}{2} = 6\pi^2 d(d-1)$ and $\PiD_S(d) = 2\pi^2(6 - d) + 2d\pi^2 = 12\pi^2$. The $R$ entry agrees with Proposition~\ref{p:R-box}: since $h_R = 0$ (Proposition~\ref{p:H-coeffs} below), $\lambda''(0)[\bth_R, \bth_R] = \PiD_R(d)$, and the box formula gives $\lambda_1(\Omega_t) = d\pi^2 + 3d(d-1)\pi^2 t^2 + O(t^3)$.

For $U$, substituting $B_U$ into~\eqref{e:coeff-universal} gives $\rho_2 = -16x_d^2$ and the nonzero entries
\[
\begin{gathered}
(A_2)_{11} = 32x_d^2 + 16x_1^2, \qquad (A_2)_{22} = 32x_d^2 + 16x_2^2, \qquad (A_2)_{kk} = -16x_d^2 \ (k \ge 3), \\
(A_2)_{12} = -16x_1x_2, \qquad (A_2)_{1d} = 16x_1x_d, \qquad (A_2)_{2d} = 16x_2x_d
\end{gathered}
\]
(and their transposes).  By~\eqref{e:P-moments} the two active diagonal entries contribute $4\pi^2 - 8$ each, the $d-3$ inactive diagonal entries $8 - \tfrac{4\pi^2}{3}$ each, the modulator entry $(d, d)$ contributes $-\tfrac{4\pi^2}{3} - 8$, and the three off-diagonal pairs contribute $-8$, $8$, $8$; altogether $\langle L_2 u_0, u_0\rangle = 8\pi^2 - \tfrac{4\pi^2(d-2)}{3} + 8d - 40$, while $\langle \rho_2 u_0, u_0\rangle = -16\bigl(\tfrac{1}{12} - \tfrac{1}{2\pi^2}\bigr)$.  Then $\PiD_U(d) = \tfrac{64\pi^2}{3} - 80$, the $d$-dependence cancelling between the two inner products.

For $T$, substituting $B_T$ into~\eqref{e:coeff-universal} gives $\rho_2 = 8(d-1)(d-2)x_d^2$ and, for $i \ne j < d$,
\[
\begin{gathered}
(A_2)_{ii} = 16x_i^2 + 8(d-3)(d-4)\,x_d^2, \qquad
(A_2)_{dd} = 8(d-1)(d-2)\,x_d^2, \\
(A_2)_{id} = -16(d-2)\,x_i x_d, \qquad
(A_2)_{ij} = 16\,x_i x_j .
\end{gathered}
\]
By~\eqref{e:P-moments} the $d-1$ active diagonal entries contribute $16\bigl(\tfrac{\pi^2}{12} + \tfrac12\bigr) + 8(d-3)(d-4)\bigl(\tfrac{\pi^2}{12} - \tfrac12\bigr)$ each, the modulator entry $8(d-1)(d-2)\bigl(\tfrac{\pi^2}{12} + \tfrac12\bigr)$, the $2(d-1)$ active--modulator entries $-4(d-2)$ each, and the $(d-1)(d-2)$ active--active entries $4$ each, while $\langle \rho_2 u_0, u_0\rangle = 8(d-1)(d-2)\bigl(\tfrac{1}{12} - \tfrac{1}{2\pi^2}\bigr)$.  Assembling gives $\PiD_T(d) = 40(d-1)(d-2) - \tfrac{16(d-1)(d-3)\pi^2}{3}$.
\end{proof}

\subsection{The Fourier-sum part \texorpdfstring{$\SiD_X(d)$}{Sigma\_X(d)}}\label{sss:Dirichlet-assembly}

The Fourier-sum part is expressed through two positive series with hyperbolic weights,
\begin{equation}\label{e:Dirichlet-sums}
\begin{aligned}
D_S \;&:=\; \sum_{m \ge 2\ \mathrm{even}} \frac{m^2\sqrt{m^2-2}\,\coth\bigl(\tfrac{\pi}{2}\sqrt{m^2-2}\bigr)}{(m^2-1)^4} \;\approx\; 0.072873, \\
D_U \;&:=\; \sum_{m \ge 2\ \mathrm{even}} \frac{m^2\sqrt{m^2-2}\,\tanh\bigl(\tfrac{\pi}{2}\sqrt{m^2-2}\bigr)}{(m^2-1)^4} \;\approx\; 0.069587,
\end{aligned}
\end{equation}
together with the elementary sums and Mittag-Leffler identities of the following lemma, recorded once for use here and in Section~\ref{ss:fd-out}.

\begin{lemma}[Elementary sums]\label{l:even-power-sums}
For $k \ge 1$ the sum $\sum_{p \ge 2\ \mathrm{even}} (p^2-1)^{-k}$ is a rational combination of $1, \pi^2, \ldots, \pi^{2\lfloor k/2 \rfloor}$; explicitly, for $k = 1, \ldots, 5$ it equals
\begin{equation}\label{e:even-power-sums}
\tfrac12, \qquad
\tfrac{\pi^2}{16} - \tfrac12, \qquad
\tfrac12 - \tfrac{3\pi^2}{64}, \qquad
\tfrac{\pi^4}{768} + \tfrac{5\pi^2}{128} - \tfrac12, \qquad
\tfrac12 - \tfrac{35\pi^2}{1024} - \tfrac{5\pi^4}{3072} ,
\end{equation}
respectively.  For odd indices, $\sum_{m \ge 3\ \mathrm{odd}}(m^2-1)^{-1} = \tfrac14$ and $\sum_{m \ge 3\ \mathrm{odd}}(m^2-1)^{-2} = \tfrac{\pi^2}{48} - \tfrac{3}{16}$.  Moreover, for $a > 0$,
\begin{equation}\label{e:ML-identities}
\sum_{m\ge 0}\frac{1}{(2m+1)^2+a^2} \;=\; \frac{\pi}{4a}\tanh\frac{\pi a}{2}, \qquad
\sum_{m\ge 1}\frac{1}{4m^2+a^2} \;=\; \frac{\pi}{4a}\coth\frac{\pi a}{2} - \frac{1}{2a^2}.
\end{equation}
\end{lemma}

\begin{proof}
Write $p = 2m$ and expand $(4m^2-1)^{-k} = 2^{-k}\bigl(\tfrac{1}{2m-1} - \tfrac{1}{2m+1}\bigr)^{k}$ by partial fractions into a combination of the $(2m \mp 1)^{-j}$, $1 \le j \le k$.  Summed over $m \ge 1$, the even-$j$ terms give rational combinations of $(1 - 2^{-j})\zeta(j) \in \pi^{j}\mathbb{Q}$, where $\zeta$ is the Riemann zeta function~\cite[Eq.~25.6.2]{DLMF}, and the odd-$j$ terms telescope; carrying this out for $k \le 5$ gives~\eqref{e:even-power-sums}.  For odd $m$, $m^2 - 1 = 4k(k+1)$ with $k = \tfrac{m-1}{2}$, and the two odd-index sums are elementary.  The identities~\eqref{e:ML-identities} follow from the partial-fraction expansion of $\coth$~\cite[Eq.~4.36.3]{DLMF} at $z = \tfrac{\pi a}{2}$, together with $\tanh z = 2\coth 2z - \coth z$.
\end{proof}

The Fourier sums in Proposition~\ref{p:H-coeffs} below reduce to the following two double series.

\begin{lemma}[Reduction of the Dirichlet double sums]\label{l:Dirichlet-reduction}
With $D_S, D_U$ as in~\eqref{e:Dirichlet-sums},
\begin{align}
\sum_{\substack{m_1, m_2 \ge 2 \\ m_1, m_2 \text{ even}}} \frac{m_1^2\, m_2^2}{(m_1^2 - 1)^2(m_2^2 - 1)^2(m_1^2 + m_2^2 - 2)} \;&=\; \frac{\pi^4}{1024} - \frac{\pi}{4}\,D_S, \label{e:red-S}\\
\sum_{\substack{m \ge 3,\,\text{odd}\\ m' \ge 2,\,\text{even}}} \frac{m^2\,(m')^2}{(m^2-1)^2((m')^2-1)^2(m^2 + (m')^2 - 2)} \;&=\; \frac{5\pi^4}{3072} - \frac{5\pi^2}{512} - \frac{\pi}{4}\,D_U. \label{e:red-U}
\end{align}
\end{lemma}

\begin{proof}
Fix the even index $m_2$ (respectively $m'$), and set $c := m_2^2 - 2 \ge 2$ and $u := c + 1 = m_2^2 - 1$.  In the remaining index the summand is, with $v := m_1^2$ (respectively $m^2$), $\tfrac{v}{(v-1)^2(v+c)} = \tfrac{c}{u^2}\cdot\tfrac{1}{v-1} + \tfrac{1}{u}\cdot\tfrac{1}{(v-1)^2} - \tfrac{c}{u^2}\cdot\tfrac{1}{v+c}$. For~\eqref{e:red-S} the inner index is even, and by Lemma~\ref{l:even-power-sums} and the second identity of~\eqref{e:ML-identities} with $a = \sqrt c$,
\[
\sum_{m_1 \ge 2\ \mathrm{even}} \frac{m_1^2}{(m_1^2 - 1)^2(m_1^2 + c)}
\;=\; \frac{\pi^2}{16\,u} - \frac{\pi\sqrt c}{4\,u^2}\coth\frac{\pi\sqrt c}{2} .
\]
Multiplying by $\tfrac{m_2^2}{(m_2^2-1)^2} = \tfrac{u+1}{u^2}$ and summing over even $m_2 \ge 2$, the first term gives $\tfrac{\pi^2}{16}\sum_{m_2}\bigl[u^{-2} + u^{-3}\bigr] = \tfrac{\pi^2}{16}\cdot\tfrac{\pi^2}{64}$ by~\eqref{e:even-power-sums}, and the second gives $-\tfrac{\pi}{4}D_S$.  For~\eqref{e:red-U} the inner index is odd and $\ge 3$; by the odd-index sums of Lemma~\ref{l:even-power-sums} and the first identity of~\eqref{e:ML-identities} without its $m = 0$ term, $\sum_{m \ge 3\ \mathrm{odd}}(m^2 + c)^{-1} = \tfrac{\pi}{4\sqrt c}\tanh\tfrac{\pi\sqrt c}{2} - \tfrac{1}{u}$, so, using $c = u - 1$,
\[
\sum_{m \ge 3\ \mathrm{odd}} \frac{m^2}{(m^2 - 1)^2(m^2 + c)}
\;=\; \Bigl(\frac{\pi^2}{48} + \frac{1}{16}\Bigr)\frac1u + \frac{3}{4u^2} - \frac{1}{u^3} - \frac{\pi\sqrt c}{4u^2}\tanh\frac{\pi\sqrt c}{2} .
\]
Multiplying by $\tfrac{u+1}{u^2}$ and summing over even $m' \ge 2$ with~\eqref{e:even-power-sums},
\[
\Bigl(\frac{\pi^2}{48} + \frac{1}{16}\Bigr)\frac{\pi^2}{64} + \frac34\Bigl(\frac{\pi^4}{768} - \frac{\pi^2}{128}\Bigr) - \Bigl(\frac{5\pi^2}{1024} - \frac{\pi^4}{3072}\Bigr) \;=\; \frac{5\pi^4}{3072} - \frac{5\pi^2}{512},
\]
and the hyperbolic term gives $-\tfrac{\pi}{4}D_U$.
\end{proof}

The Fourier coefficients $H_\bfn$ and the resulting sums are collected in the following proposition.

\begin{proposition}\label{p:H-coeffs}
For the $R$-deformation, $h_R \equiv 0$, so $H_\bfn^R = 0$ for every $\bfn$.  For $X \in \{S, T, U\}$, the Fourier coefficients $H_\bfn^X = \langle h_X, \Psi_\bfn\rangle$ vanish except for the following $\bfn$, where every index not displayed equals $1$, $n_1, n_2 \ge 2$ are even, $n \ge 3$ is odd, and $p \ge 2$ is even and sits on the modulator axis $d$:
\[
\renewcommand{\arraystretch}{1.7}
\begin{array}{c|l|c}
X & \multicolumn{1}{c|}{\bfn} & H^X_\bfn \\
\hline
S & (n_1, n_2, 1, \ldots, 1) & \tfrac{-64\, n_1 n_2}{(n_1^2-1)(n_2^2-1)} \\
T\ \text{(Family I)} & n \text{ in position } i < d,\ p \text{ in position } d & \tfrac{64\, n\, p}{(n^2-1)(p^2-1)} \\
T\ \text{(Family II)} & (1, \ldots, 1, p) & \tfrac{16(d-1)\, p\,(p^2-5)}{(p^2-1)^2} \\
U & n \text{ in position } i \in \{1, 2\},\ p \text{ in position } d & \tfrac{\pm 64\, n\, p}{(n^2-1)(p^2-1)}\ \ (+ \text{ for } i = 1,\ - \text{ for } i = 2)
\end{array}
\]
Consequently the Fourier-sum parts $\SiD_X(d)$ of~\eqref{e:lam2-final} are
\begin{equation}\label{e:QX}
\begin{gathered}
\SiD_R(d) \;=\; 0, \qquad
\SiD_S(d) \;=\; -8\pi^2 + \frac{2048\,D_S}{\pi}, \qquad
\SiD_U(d) \;=\; 160 - \frac{80\pi^2}{3} + \frac{4096\,D_U}{\pi}, \\
\SiD_T(d) \;=\; (d-1)\Bigl[80 - \frac{40\pi^2}{3} + \frac{2048\,D_U}{\pi}\Bigr] \;-\; 8(d-1)^2\,(10 - \pi^2).
\end{gathered}
\end{equation}
\end{proposition}

\begin{proof}
For the $R$-deformation, $\rho_1 = 0$ and $A_1$ is constant diagonal with $\tr A_1 = 0$ (Proposition~\ref{p:Iso-irrep}, Lemma~\ref{l:coeff-expansion}), so $h_R = -\sum_{k} (A_1)_{kk}\, \partial_{x_k}^2 u_0 = \pi^2 \bigl(\tr A_1\bigr)\, u_0 = 0$.  For $S$ and $U$, $\rho_1 = 0$, so $h_X = -\nabla \cdot(A_1 \nabla u_0)$; for $T$ both terms of~\eqref{e:h-defn} contribute. Substituting the $A_1$ of Lemma~\ref{l:coeff-expansion} and using $\partial_{x_k} u_0 = 2^{d/2}\phi_1'(x_k)\prod_{j \ne k}\phi_1(x_j)$ gives, up to an overall sign (immaterial, since only $H_\bfn^2$ enters~\eqref{e:lam2-final}),
\begin{align}
h_S(\x) \;&=\; -4\cdot 2^{d/2}\, \phi_1'(x_1)\,\phi_1'(x_2)\prod_{k=3}^{d}\phi_1(x_k), \label{e:hS}\\
h_T(\x) \;&=\; -8\cdot 2^{d/2}\,\sum_{i=1}^{d-1} x_i\, \phi_1'(x_i)\,\phi_1'(x_d)\!\!\prod_{k \ne i, d}\!\!\phi_1(x_k) \;+\; 8(d-1)\pi^2 \cdot 2^{d/2}\, x_d \prod_{k=1}^{d}\phi_1(x_k), \label{e:hT}\\
h_U(\x) \;&=\; -8 \cdot 2^{d/2}\,\bigl[x_1 \phi_1'(x_1)\phi_1(x_2) - x_2 \phi_1(x_1)\phi_1'(x_2)\bigr]\phi_1'(x_d)\prod_{k=3}^{d-1}\phi_1(x_k). \label{e:hU}
\end{align}
For $T$, the $(\operatorname{tr} B)\, I_d$ part of $A_1$ contributes $-(\operatorname{tr} B)\,\Delta u_0 = 4d(d-1)\pi^2 x_d u_0$, which cancels exactly against $-\lambda_0 \rho_1 u_0$; the remaining $-B - B^T$ part produces the two terms in~\eqref{e:hT}.

Each $H_\bfn = \langle h_X, \Psi_\bfn\rangle$ factors over axes: by~\eqref{e:hS}--\eqref{e:hU} the axis factors of $h_X$ are $\phi_1$, $\phi_1'$, $x\phi_1$ or $x\phi_1'$, so each axis contributes one of the pairings~\eqref{e:mom-offdiag}--\eqref{e:mom-diag} with $m = 1$. An inactive axis contributes $\langle \phi_1, \phi_{n_k}\rangle_0 = \tfrac12 \delta_{1 n_k}$ and hence forces $n_k = 1$; an axis carrying $\phi_1'$ or $x\phi_1$ forces $n_k$ even, and one carrying $x\phi_1'$ forces $n_k$ odd. This yields the stated supports, and the $2^{d/2}$ prefactors of $h_X$ and $\Psi_\bfn$ combine with the inactive factors $\tfrac12$ to give
\[
H^S_{\bfn} = -16\,\langle \phi_1', \phi_{n_1}\rangle_0 \langle \phi_1', \phi_{n_2}\rangle_0,
\qquad
H^{T, \mathrm{I}}_{\bfn} = H^{U}_{\bfn} = -32\,\langle \phi_1', \phi_{n}\rangle_1 \langle \phi_1', \phi_{p}\rangle_0 ,
\]
which are the stated values by~\eqref{e:mom-offdiag}.

At the exceptional index $\bfn = (1, \ldots, 1, p)$, for $T$ both the active-axis terms of~\eqref{e:hT}, paired with index $1$ on every active axis (contributing $\langle \phi_1', \phi_1\rangle_1 = -\tfrac14$ each) and the modulator term of~\eqref{e:hT} (contributing $16(d-1)\pi^2 \langle \phi_1, \phi_p\rangle_1$) survive, giving
\[
H^{T, \mathrm{II}}_{(1, \ldots, 1, p)}
\;=\; 8(d-1)\,\langle \phi_1', \phi_p\rangle_0 \;+\; 16(d-1)\pi^2\, \langle \phi_1, \phi_p\rangle_1
\;=\; \frac{16(d-1)p}{p^2-1} - \frac{64(d-1)p}{(p^2-1)^2},
\]
the tabulated value. For $U$, the two cases $i = 1, 2$ extend to $n = 1$, where they have equal magnitude and opposite signs (the relative $-$ sign between the $x_1$ and $x_2$ terms in $h_U$) and cancel exactly; this is why the support excludes $\bfn = (1, 1, 1, \ldots, 1, p)$.

Turning to the Fourier sums, for the $S$-deformation $|H_{(n_1, n_2, 1, \ldots, 1)}^S|^2 = \tfrac{4096\, n_1^2 n_2^2}{(n_1^2-1)^2 (n_2^2-1)^2}$. The gap is $\pi^2(|\bfn|^2 - d) = \pi^2(n_1^2 + n_2^2 - 2)$.  Summing over $n_1, n_2$ even $\ge 2$ gives $-\tfrac{8192}{\pi^2}$ times the double sum~\eqref{e:red-S}, which by Lemma~\ref{l:Dirichlet-reduction} is $-\tfrac{8192}{\pi^2}\bigl(\tfrac{\pi^4}{1024} - \tfrac{\pi}{4}D_S\bigr) = -8\pi^2 + \tfrac{2048 D_S}{\pi}$, the stated $\SiD_S(d)$.

For the $U$-deformation, $|H^U_{(n, 1, \ldots, p)}|^2 = \tfrac{4096\, n^2\, p^2}{(n^2 - 1)^2(p^2 - 1)^2}$ with gap $\pi^2(n^2 + p^2 - 2)$.  Summing yields $\tfrac{4096}{\pi^2}$ times the double sum~\eqref{e:red-U} per family, doubled by the $1$--$2$ swap symmetry and multiplied by $-2$ (from~\eqref{e:lam2-final}): by Lemma~\ref{l:Dirichlet-reduction},
\[
\SiD_U(d) \;=\; -\frac{16384}{\pi^2}\Bigl(\frac{5\pi^4}{3072} - \frac{5\pi^2}{512} - \frac{\pi}{4}D_U\Bigr) \;=\; 160 - \frac{80\pi^2}{3} + \frac{4096\,D_U}{\pi} .
\]

For the $T$-deformation, Family~\textup{I} has the same summand as the two $U$ cases, but is now summed over the $d-1$ active-axis positions $i$ of the odd index $n$, contributing $(d-1)$ times half of $\SiD_U(d)$, that is $(d-1)\bigl[80 - \tfrac{40\pi^2}{3} + \tfrac{2048 D_U}{\pi}\bigr]$.

Family~\textup{II} has $|H_{(1, \ldots, 1, p)}^{T, \mathrm{II}}|^2 = \tfrac{256(d-1)^2\, p^2 (p^2-5)^2}{(p^2-1)^4}$, and the gap at $\bfn = (1, \ldots, 1, p)$ is $\pi^2(p^2 - 1)$.  The resulting single series is elementary: with $u := p^2 - 1$, $\tfrac{p^2 (p^2-5)^2}{(p^2-1)^5} = \tfrac{(u+1)(u-4)^2}{u^5} = \tfrac{1}{u^2} - \tfrac{7}{u^3} + \tfrac{8}{u^4} + \tfrac{16}{u^5}$, so by Lemma~\ref{l:even-power-sums}
\[
\sum_{p \ge 2\ \mathrm{even}} \frac{p^2 (p^2-5)^2}{(p^2-1)^5}
\;=\; \frac{\pi^2 (10 - \pi^2)}{64} \;\approx\; 0.020109 .
\]
Family~\textup{II} therefore contributes $-2 \cdot \tfrac{256(d-1)^2}{\pi^2} \cdot \tfrac{\pi^2(10 - \pi^2)}{64} = -8(d-1)^2(10 - \pi^2)$, and adding the Family~\textup{I} contribution yields the stated total for $T$.
\end{proof}

\subsection{Proof of Theorem~\ref{t:Lam-Hessian-d}}\label{ss:Lam-Hessian}

\begin{proof}[Proof of Theorem~\ref{t:Lam-Hessian-d}]
Combining~\eqref{e:nu-formula} with~\eqref{e:lam2-final},
\[
\nu^{\blam}_X \;=\; \frac{\PiD_X(d) \,+\, \SiD_X(d) \,+\, 2\pi^2\, V''(0)[\bth_X, \bth_X]}{\|\bth_X\|^2},
\]
in which $\PiD_X(d)$ is given by Lemma~\ref{l:poly-table}, $\SiD_X(d)$ by~\eqref{e:QX}, $V''(0)[\bth_X, \bth_X]$ by Proposition~\ref{p:Iso-irrep}, and $\|\bth_X\|^2$ by~\eqref{e:norms}.  Substituting and simplifying gives~\eqref{e:nu-Dirichlet}; for $T$, every term of the numerator carries a factor $(d-1)$, which cancels against $\|\bth_T\|^2 = 2^d(d-1)$.  The value $\nu^{\blam}_R = 16\pi^2/2^d$ agrees with the box computation of Proposition~\ref{p:R-box}.

For positivity, $\tfrac{16\pi^2}{3} - 40 \approx 12.64$ is positive, being the inequality $\pi^2 > \tfrac{15}{2}$, and $D_S, D_U > 0$ as sums of positive terms, so $\nu^{\blam}_R$, $\nu^{\blam}_S$, $\nu^{\blam}_T > 0$ for every $d \ge 2$.  For $\nu^{\blam}_U$ we need $D_U > \tfrac{\pi}{2048}\bigl(\tfrac{16\pi^2}{3} - 40\bigr) \approx 0.0194$, and already the term $m = 2$ of~\eqref{e:Dirichlet-sums} gives $D_U > \tfrac{4\sqrt2}{81}\tanh\tfrac{\pi}{\sqrt2} > \tfrac{2\sqrt2}{81} > 0.034$.
\end{proof}

%%%%%%%%%%%%%%%%%%%%%%%%%%%%%%%%%%
%%%%%%%%%%%%%%%%%%%%%%%%%%%%%%%%%%
\section{The Laplace--Neumann eigenvalue}\label{s:LNshape}

The lowest part of the spectrum of $-\Delta$ on $\cube$ with Neumann boundary conditions consists of $\mu_0 = 0$ (constant eigenfunction) and the $d$-fold cluster $\mu_1 = \cdots = \mu_d = \pi^2$ attained on $\{\sqrt{2}\sin(\pi x_k)\}_{k=1}^{d}$ (recall $\cube = [-\tfrac12, \tfrac12]^d$, on which $\sin(\pi x_k)$ satisfies the Neumann condition $\partial_{x_k}u = 0$ at $x_k = \pm\tfrac12$).  Under a generic perturbation of $\cube$, this cluster splits into $d$ distinct eigenvalues.

As in Section~\ref{ss:Material}, we work on the fixed cube throughout: pulling back by a real-analytic family of diffeomorphisms $\Phi_t$ with $\Phi_0 = \mathrm{id}$ converts the Neumann eigenvalue problem on $\Omega_t = \Phi_t(\cube)$ into the weighted problem
\begin{equation}\label{e:pullback-Neumann}
\begin{aligned}
-\nabla \cdot \bigl(A(t, \x)\, \nabla u\bigr) \;&=\; \mu(t)\, \rho(t, \x)\, u && \text{in } \cube, \\
\bigl(A(t, \x)\, \nabla u\bigr) \cdot \mathbf{n} \;&=\; 0 && \text{on } \partial\cube,
\end{aligned}
\end{equation}
with the same $\rho(t, \x) = \det D\Phi_t(\x)$ and $A(t, \x) = \rho(t, \x)\, D\Phi_t(\x)^{-1} D\Phi_t(\x)^{-T}$ as there (for $|t|$ small $\rho > 0$, so the absolute value in the weight is immaterial).  Write $L(t)$ for the associated self-adjoint operator on $L^2(\cube, \rho\, d\x)$; its eigenvalues are those of the Neumann Laplacian on $\Omega_t$. 

\subsection{Auxiliary functional}\label{ss:Aux}

None of the split eigenvalues is individually smooth as a functional on $\P$; in particular $\bmu$ is not differentiable at $\cube$.  The plan of Section~\ref{s:LDshape} still applies --- pull the eigenvalue problem back to the fixed cube, expand to second order in $t$, and split the result into an elementary part and a gap-weighted series --- but to a smooth surrogate in place of the eigenvalue itself.  We use the auxiliary functional $f_d(\Omega) = \tfrac{1}{d|\Omega|^{2/d}}\sum_{i=1}^d \tfrac{1}{\mu_i(\Omega)}$ of~\eqref{e:fd-defn} --- the volume-normalized harmonic mean of the cluster --- where $\mu_1(\Omega) \le \cdots \le \mu_d(\Omega)$ are the lowest $d$ nontrivial Neumann eigenvalues.

Smoothness of $f_d$ through the cluster degeneracy is a standard consequence of analytic perturbation theory.  For a real-analytic family $\Phi_t$ with $\Phi_0 = \mathrm{id}$ --- for instance $\Phi_t = \mathrm{id} + t\bth$, or the multilinear family~\eqref{e:Phi-delta} --- both $\rho$ and $A$ are real-analytic in $t$, and $A$ is uniformly positive definite for $|t|$ small, so $L(t)$ is a real-analytic family of type~\textup{(B)} in Kato's sense~\cite[Ch.~VII, \S 4]{Kato1995}.  Its lowest nontrivial eigenvalue cluster --- of multiplicity $d$ at $t = 0$, isolated both from $0$ and from $\mu_{d+1}(0) = 2\pi^2$ --- therefore has a spectral projector $P_{\mathrm{cluster}}(t)$ analytic in $t$, by the Rellich--Kato theorem, the Riesz projector being a resolvent contour integral around the isolated cluster; $P_{\mathrm{cluster}}(0)$ is the orthogonal projection onto the $\pi^2$-eigenspace.  Consequently the \emph{cluster sum}
\begin{equation}\label{e:cluster-sum-defn}
\calS(t) \;:=\; \sum_{i=1}^d \frac{1}{\mu_i(\Omega_t)},
\end{equation}
which is the trace of the inverse of $L(t)$ on the range of $P_{\mathrm{cluster}}(t)$, is real-analytic at $t = 0$~\cite[Ch.~VII, Thm.~3.9]{Kato1995}, and $\Omega \mapsto f_d(\Omega)$ extends analytically through $\cube$, although the individual $\mu_i(\Omega_t)$ do not.

\begin{remark}[From $f_d$ to $\bmu$]\label{r:fd-bound}
Since $\mu_i(\Omega) \ge \mu_1(\Omega)$ for each $i$, summing $\tfrac{1}{\mu_i(\Omega)} \le \tfrac{1}{\mu_1(\Omega)}$ and dividing by $d|\Omega|^{2/d}$ gives $f_d(\Omega) \le \tfrac{1}{\bmu(\Omega)}$, with equality at $\cube$, where all $\mu_i$ coincide.  Hence if $f_d$ has a strict local minimum at $\cube$, then for $\Omega$ near $\cube$ and not similar to it, $\tfrac{1}{\bmu(\Omega)} \ge f_d(\Omega) > f_d(\cube) = \tfrac{1}{\bmu(\cube)}$: $\bmu$ has a strict local maximum at $\cube$.
\end{remark}

It therefore suffices to compute the Hessian of $f_d$ at $\cube$ in symmetry-adapted coordinates and verify that its four scalar eigenvalues on the blocks of~\eqref{e:irrep-decomp} are positive.

\subsection{Cluster reduction and the cluster-sum expansion}\label{ss:cluster-reduction}

Fix a representative field $\bth$, put $\bdelta = t\bth$, and apply~\eqref{e:pullback-Neumann} to the multilinear family $\Phi_{\bdelta}$ of~\eqref{e:Phi-delta}.  Let $\{\Xi_{\bfm}\}_{\bfm \in \mathbb N_0^d}$ be the Neumann eigenbasis of $\cube$, including the constant $\Xi_{\bzero} \equiv 1$, with eigenvalues $\pi^2|\bfm|^2$, and let $W := \mathrm{span}\{\psi_1, \ldots, \psi_{d}\} \;\subset\; H^1(\cube),  \psi_k := -\Xi_{\mathbf e_k} = \sqrt2\,\sin(\pi x_k)$, be the cluster eigenspace at $t = 0$; since $|\bfm|^2 = 1$ exactly when $\bfm$ is a standard basis vector, the Neumann modes outside the cluster are the $\Xi_\bfm$ with $|\bfm| \ne 1$.  In this basis the pulled-back problem~\eqref{e:pullback-Neumann} is the symmetric pencil $K(t)\mathbf u = \mu\, M(t)\mathbf u$ with
\begin{equation}\label{e:GEP}
K_{\alpha\beta}(t) \;=\; \int_{\cube} \nabla \Xi_{\alpha}^{T} A(t, \x)\, \nabla \Xi_{\beta}\, d\x,
\qquad
M_{\alpha\beta}(t) \;=\; \int_{\cube} \Xi_{\alpha} \Xi_{\beta}\, \rho(t, \x)\, d\x .
\end{equation}
Both are real-analytic in $t$, with $K(0) = \mathrm{diag}(\pi^2|\bfm|^2)$ and $M(0) = I$; in particular $K(0)|_W = \pi^2 I_d$ and $M(0)|_W = I_d$.  We write $K_1, K_2, M_1, M_2$ for the Taylor coefficients of $K(t), M(t)$ at $t = 0$, set $H := K_1 - \pi^2 M_1$, and abbreviate by $|_W$ the restriction of a matrix to the $d \times d$ cluster block.  Expanding $A(t) = I + tA_1 + t^2A_2 + O(t^3)$ and $\rho(t) = 1 + t\rho_1 + t^2\rho_2 + O(t^3)$ as in Lemma~\ref{l:coeff-expansion} --- these coefficients are geometric, so its table serves here unchanged --- and substituting into~\eqref{e:GEP}, the cluster-block entries are
\begin{equation}\label{e:cluster-coeff-formulas}
(K_j)_{kl} \;=\; \int_{\cube} (A_j \nabla\psi_k)\cdot\nabla\psi_l\,d\x, \qquad
(M_j)_{kl} \;=\; \int_{\cube} \rho_j\,\psi_k\psi_l\,d\x \qquad (j = 1, 2),
\end{equation}
where $\nabla\psi_k = \sqrt2\,\pi\cos(\pi x_k)\,\mathbf e_k$, and the outside-coupling entries against $\Xi_{\bfm} = \prod_i \chi_{m_i}(x_i)$, with $\chi_0 \equiv 1$ and $\chi_n(x) = \sqrt2\cos\bigl(n\pi(x + \tfrac12)\bigr)$ for $n \ge 1$, are
\begin{equation}\label{e:H1-coupling-formula}
H_{k, \bfm} \;=\; \int_{\cube} (A_1 \nabla\psi_k)\cdot\nabla\Xi_{\bfm}\,d\x \;-\; \pi^2\int_{\cube} \rho_1\,\psi_k\,\Xi_{\bfm}\,d\x .
\end{equation}
Every such integral factors over the axes, and the nonzero one-dimensional factors are elementary.  Writing $\int$ for $\int_{-1/2}^{1/2}\cdots dx$, the factors that occur below are
\begin{equation}\label{e:Neumann-1D}
\begin{aligned}
&\int \chi_n = \delta_{n0}, \quad
\int \cos(\pi x)\,\chi_n' = -\tfrac{\pi}{\sqrt2}\,\delta_{n1}, \quad
\int \sin(\pi x)\,\chi_n = -\tfrac{1}{\sqrt2}\,\delta_{n1}, \quad
\int \chi_n' = -2\sqrt2 \ (n \text{ odd}), \\
&\int \cos(\pi x)\,\chi_n = \begin{cases} \tfrac2\pi, & n=0,\\ \tfrac{-2\sqrt2}{\pi(n^2-1)}, & n\ \text{even} \ge 2,\end{cases} \quad
\int x\cos(\pi x)\,\chi_n = \begin{cases} -\tfrac{\sqrt2}{4\pi}, & n = 1,\\ \tfrac{\sqrt2}{\pi(n^2-1)}, & n\ \text{odd} \ge 3,\end{cases} \\
&\int x\,\chi_n = \tfrac{-2\sqrt2}{\pi^2 n^2} \ (n \text{ odd}), \quad
\int x\,\chi_n' = \sqrt2 \ (n \text{ even} \ge 2), \quad
\int x\cos(\pi x)\,\chi_n' = \tfrac{4\sqrt2\, n^2}{\pi(n^2-1)^2} \ (n \text{ even} \ge 2), \\
&\int x\sin(\pi x)\,\chi_n = \begin{cases} \tfrac{2}{\pi^2}, & n = 0,\\ \tfrac{2\sqrt2\,(n^2+1)}{\pi^2(n^2-1)^2}, & n\ \text{even} \ge 2,\end{cases}
\end{aligned}
\end{equation}
each vanishing for the parity of $n$ not listed.

Since $(K(t), M(t))$ is a real-analytic self-adjoint family and $\pi^2$ is an isolated eigenvalue of multiplicity $d$ at $t = 0$, the Rellich--Kato theorem~\cite[Ch.~VII, Thm.~3.9]{Kato1995} provides $d$ real-analytic branches $\mu_i(t) = \pi^2 + t\,\mu^{(1)}_i + t^2 \mu^{(2)}_i + O(t^3)$ with real-analytic eigenvectors $\mathbf u_i(t) = \mathbf v_i + t\,\mathbf w_i + O(t^2)$, $\mathbf v_i \in W$, $i = 1, \ldots, d$, agreeing as a set with the cluster $\{\mu_1(\Omega_t), \ldots, \mu_d(\Omega_t)\}$ but in general not ordered --- which is why the ordered eigenvalues $\mu_i(\Omega_t)$ are not individually smooth.  Since $\calS$ is a symmetric function of the cluster, the labelling is immaterial below.

The coefficients $\mu^{(1)}_i$, $\mu^{(2)}_i$ and $\mathbf v_i$ follow from the usual Rayleigh--Schr\"odinger expansion for a simple eigenvalue, applied to each branch~\cite[Ch.~II, \S 2]{Kato1995}.  Substituting the expansions into the eigenvalue equation $K(t)\mathbf u_i(t) = \mu_i(t)\,M(t)\mathbf u_i(t)$, the order-$t$ equation is $(K_0 - \pi^2)\mathbf w_i + (H - \mu^{(1)}_i)\mathbf v_i = \bzero$; its projection onto $W$ says that $\mu^{(1)}_i$ is an eigenvalue of $H|_W$ with eigenvector $\mathbf v_i$, and its projection onto $\Xi_\bfm$, $|\bfm| \ne 1$, gives $(\mathbf w_i)_\bfm = (H\mathbf v_i)_\bfm / (\pi^2(1 - |\bfm|^2))$.  Pairing the order-$t^2$ equation with $\mathbf v_i$ then gives
\begin{equation}\label{e:RS-2nd}
\mu_i^{(2)} \;=\; \mathbf v_i^{T}\bigl(K_2 - \pi^2 M_2 - \mu^{(1)}_i M_1\bigr)\big|_W\, \mathbf v_i \;+\; \sum_{|\bfm| \ne 1} \frac{\bigl(H \mathbf v_i\bigr)_{\bfm}^2}{\pi^2\bigl(1 - |\bfm|^2\bigr)},
\end{equation}
the $W$-component of $\mathbf w_i$ dropping out.  (As in Section~\ref{ss:Material}, a subscript on an expansion coefficient records the power of $t$ it carries, while superscripts $(1), (2)$ are used on $\mu_i$, whose subscript is reserved for the eigenvalue index.)  When the $\mu^{(1)}_i$ are distinct, $\mathbf v_1, \ldots, \mathbf v_d$ is the orthonormal eigenbasis of $H|_W$; in general the $\mathbf v_i$ are determined only at higher order, but $f_d$ needs only the trace of~\eqref{e:RS-2nd} over the cluster, which is the same for every orthonormal eigenbasis of $H|_W$.

\begin{lemma}[Cluster-sum expansion]\label{l:cluster-sum-formula}
With the notation above, the cluster sum $\calS$ of~\eqref{e:cluster-sum-defn} admits the expansion
\begin{equation}\label{e:cluster-sum-taylor}
\calS(t) \;=\; \frac{d}{\pi^2} \;-\; \frac{\mathrm{tr}\, H|_W}{\pi^4}\, t \;+\; \bigl[\PiN + \SiN\bigr]\, t^2 \;+\; O(t^3),
\end{equation}
where
\begin{align*}
\pi^6\,\PiN \;&=\; \mathrm{tr}\bigl((K_1|_W)^2\bigr) - \pi^2\,\mathrm{tr}\bigl(K_1|_W M_1|_W\bigr) - \pi^2\,\mathrm{tr}\, K_2|_W + \pi^4\,\mathrm{tr}\, M_2|_W, \\
\pi^6\,\SiN \;&=\; \sum_{k=1}^d \sum_{|\bfm| \ne 1} \frac{H_{k,\bfm}^2}{|\bfm|^2 - 1}.
\end{align*}
\end{lemma}

\begin{proof}
Expanding $\mu_i(t)^{-1} = \pi^{-2} - \pi^{-4}\bigl(t\,\mu^{(1)}_i + t^2\mu^{(2)}_i\bigr) + \pi^{-6}\,t^2\bigl(\mu^{(1)}_i\bigr)^2 + O(t^3)$ and summing over the cluster,
\[
\calS(t) \;=\; \frac{d}{\pi^2} \;-\; \frac{t}{\pi^4}\sum_i \mu^{(1)}_i \;+\; t^2\biggl[\frac{\sum_i \bigl(\mu^{(1)}_i\bigr)^2}{\pi^6} \,-\, \frac{\sum_i \mu^{(2)}_i}{\pi^4}\biggr] \;+\; O(t^3).
\]
The $\mu^{(1)}_i$ are the eigenvalues of $H|_W$, so $\sum_i \mu^{(1)}_i = \mathrm{tr}\, H|_W$ and $\sum_i (\mu^{(1)}_i)^2 = \mathrm{tr}\bigl((H|_W)^2\bigr)$.  Summing~\eqref{e:RS-2nd} and using $\sum_i \mathbf v_i \mathbf v_i^T = I_d$ and $\sum_i \mu^{(1)}_i\, \mathbf v_i \mathbf v_i^T = H|_W$,
\[
\sum_i \mu^{(2)}_i \;=\; \mathrm{tr}\bigl(K_2 - \pi^2 M_2\bigr)\big|_W \;-\; \mathrm{tr}\bigl(H|_W M_1|_W\bigr) \;-\; \frac{1}{\pi^2}\sum_{k=1}^d \sum_{|\bfm| \ne 1} \frac{H_{k, \bfm}^2}{|\bfm|^2 - 1}.
\]
Substituting these three sums and expanding $H|_W = K_1|_W - \pi^2 M_1|_W$ gives~\eqref{e:cluster-sum-taylor}; the terms in $\mathrm{tr}\bigl((M_1|_W)^2\bigr)$ cancel.  Both $\sum_i \mu^{(2)}_i$ and $\sum_i (\mu^{(1)}_i)^2$ are traces of operators on $W$, hence the same for every orthonormal eigenbasis of $H|_W$, as claimed.
\end{proof}

Along the representative field $\bth_X$ we write $\PiN_X(d)$ and $\SiN_X(d)$ for these two coefficients, and compute them in Lemma~\ref{l:cluster-matrices} and Section~\ref{ss:fd-out}.  They are the Neumann counterparts of the polynomial and Fourier-sum parts of the Dirichlet analysis.
One difference: since $\Xi_{\bzero}$ lies \emph{below} the cluster, one gap in $\SiN_X(d)$ is negative, so $\SiN_X(d)$ need not have a fixed sign, whereas $\SiD_X(d) \le 0$ always.

\begin{remark}[First-order behavior of $\bmu$]\label{r:one-sided}
Since $f_d$ is smooth and $B_d$-invariant, Proposition~\ref{p:Crit-Block}(i) makes it critical at $\cube$; as $V'(0) = 0$, the linear term of~\eqref{e:cluster-sum-taylor} must vanish, so $\sum_{i=1}^d \mu^{(1)}_i = \mathrm{tr}\, H|_W = 0$ for every $\bth \in \TP$. Hence $\min_i \mu^{(1)}_i \le 0 \le \max_i \mu^{(1)}_i$.  Since $\mu_1(\Omega_t)$ is the smallest of the $d$ branches, its one-sided derivatives at $t = 0$ along $+\bth$ and $-\bth$ are $\min_i \mu^{(1)}_i$ and $-\max_i \mu^{(1)}_i$, both non-positive; the same holds for $\bmu$, because $V'(0) = 0$.  This is the sense in which $\cube$ is critical for $\bmu$ in Theorem~\ref{t:Main}.  When $H|_W \ne 0$ the cluster splits at first order and $\bmu$ decreases \emph{linearly} in both directions, so $\cube$ is a corner maximum along $\bth$; the second-order analysis below is needed exactly for the directions with $H|_W = 0$, which by Lemma~\ref{l:cluster-matrices} are the $T$- and $U$-deformations.  For $\bth_R$ this is visible in closed form (Proposition~\ref{p:R-box}): the directional derivatives of $\bmu$ are $-2(d-1)\pi^2$ along $+\bth_R$ and $-2\pi^2$ along $-\bth_R$.
\end{remark}

Since $f_d = \tfrac{1}{d} V^{-2/d}\, \calS$, all first-order terms vanish at $\cube$ by Proposition~\ref{p:Crit-Block}(i), and $V(t)^{-2/d} = 1 - \tfrac{1}{d}V''(0)[\bth_X, \bth_X]\,t^2 + O(t^3)$, the expansion~\eqref{e:cluster-sum-taylor} gives the assembly formula
\begin{equation}\label{e:nu-fd-assembly}
\nu^{f_d}_X \;=\; \frac{2}{d\,\|\bth_X\|^2}\Bigl(\PiN_X(d) + \SiN_X(d) - \frac{V''(0)[\bth_X, \bth_X]}{\pi^2}\Bigr),
\end{equation}
the analogue of~\eqref{e:nu-formula}, which we use for all four blocks.

The cluster blocks $K_1|_W$, $K_2|_W$, $M_1|_W$, $M_2|_W$, and hence $\PiN_X(d)$, are computed in the following lemma for all four deformations; the coupling sums $\SiN_X(d)$ are computed in Section~\ref{ss:fd-out}.

\begin{lemma}[Cluster matrices and cluster part]\label{l:cluster-matrices}
For each $\bth_X$, $X \in \{R, S, T, U\}$, one has $M_1|_W = 0$, and
{\footnotesize\renewcommand{\arraystretch}{2.0}\setlength{\arraycolsep}{1pt}
\[
\begin{array}{c|ccc|c}
X & \mathrm{tr}\bigl((K_1|_W)^2\bigr) & \mathrm{tr}\, K_2|_W & \mathrm{tr}\, M_2|_W & \PiN_X(d) \\
\hline
R & 4d(d-1)\pi^4 & \tfrac{1}{2}d(d-1)(6-d)\,\pi^2 & -\tfrac{1}{2}d^2(d-1) & \dfrac{d(d-1)}{\pi^2} \\
S & 512 & (6-d)\,\pi^2 & -d & \dfrac{512 - 6\pi^4}{\pi^6} \\
T & 0 & -4d(d-1) + \tfrac{2}{3}(d-1)(d^2-6d+12)\,\pi^2 & \tfrac{2}{3}d(d-1)(d-2) + \tfrac{4(d-1)(d-2)}{\pi^2} & \dfrac{8(d-1)^2}{\pi^4} + \dfrac{8(d-1)(d-3)}{3\pi^2} \\
U & 0 & -8 + \tfrac{4}{3}(8-d)\,\pi^2 & -\tfrac{4d}{3} - \tfrac{8}{\pi^2} & -\dfrac{32}{3\pi^2}
\end{array}
\]}
\end{lemma}

\begin{proof}
Throughout we use the coefficients of Lemma~\ref{l:coeff-expansion}: the tabulated entries follow by substituting them into~\eqref{e:cluster-coeff-formulas} and evaluating with the one-dimensional integrals above.  For the vanishing of $K_1|_W$ and $M_1|_W$, recall in particular that $\rho_1 = 0$ for $U$ and $\rho_1 = 4(d-1)x_d$ for $T$, that for $U$ the only nonzero entries of $A_1$ are $(A_1)_{11} = -8x_d$, $(A_1)_{22} = 8x_d$ and $(A_1)_{1,d} = -4x_1$, $(A_1)_{2,d} = 4x_2$ (with $(A_1)_{kk} = 0$ for $k \ge 3$), and that for $T$, $(A_1)_{kk} = 4(d-3)x_d$ for $k < d$, $(A_1)_{dd} = 4(d-1)x_d$, $(A_1)_{k,d} = -4x_k$ for $k < d$.

\emph{(i) for $U$.} Since $\nabla\psi_k$ is supported on axis $k$,
\[
(K_1)_{k, l} \;=\; \int_{\cube} (A_1 \nabla \psi_k) \cdot \nabla \psi_l \,d\x \;=\; 2\pi^2 \int_{\cube} (A_1)_{k, l}\,\cos(\pi x_k)\cos(\pi x_l)\,d\x .
\]
For $k = l$: $(A_1)_{k,k}$ is either zero (for $k \ge 3$) or equals $\pm 8 x_d$ (for $k = 1, 2$). In the latter case, integrating $x_d$ over $[-1/2, 1/2]$ gives zero by parity. For $k \ne l$: $(A_1)_{k, l}$ is nonzero only for $\{k, l\} = \{1, d\}$ (value $-4 x_1$) or $\{2, d\}$ (value $+4 x_2$); in each case, $\int_{-1/2}^{1/2} x_k \cos(\pi x_k)\,dx_k = 0$ by parity ($x_k$ is odd, $\cos(\pi x_k)$ is even). Hence $K_1|_W = 0$.

\emph{(i) for $T$.} Here the diagonal entries $(A_1)_{kk}$ are all odd in $x_d$, the off-diagonals $(A_1)_{k,d}$ are odd in $x_k$, and $\rho_1$ is odd in $x_d$.  Hence $\int (A_1)_{kk} \cos^2(\pi x_k)\,d\mathbf x = 0$. Off-diagonals $(A_1)_{k, d} \cos(\pi x_k)\cos(\pi x_d)$ are odd in $x_k$, so the $x_k$ integral vanishes. Hence $K_1|_W = 0$. For $(M_1)_{k, l} = \int \rho_1 \psi_k \psi_l \, d\mathbf x = 2\cdot 4(d-1) \int x_d \sin(\pi x_k)\sin(\pi x_l)\,d\mathbf x$: if $k \ne d$ and $l \ne d$, the $x_d$ integral vanishes by parity; if $k = d$ (or $l = d$), we use $\int_{-1/2}^{1/2}\sin(\pi x_l)\,dx_l = 0$ (odd over symmetric interval). Hence $M_1|_W = 0$.

\emph{(ii) The cluster part.}  Since $M_1|_W = 0$, the second term in $\PiN$ drops and
\[
\PiN_X(d) \;=\; \frac{\mathrm{tr}\bigl((K_1|_W)^2\bigr) - \pi^2\,\mathrm{tr}\, K_2|_W + \pi^4\,\mathrm{tr}\, M_2|_W}{\pi^6} ;
\]
substituting the first three columns gives the last.  For $R$ the $d$-dependence cancels between the last two traces; for $S$ and $U$ it cancels entirely, leaving $d$-independent values; and for $T$, where the two traces are separately cubic in $d$, $-\pi^2\,\mathrm{tr}\, K_2|_W + \pi^4\,\mathrm{tr}\, M_2|_W = 8(d-1)^2\pi^2 + \tfrac{8}{3}(d-1)(d-3)\pi^4$. \end{proof}

\subsection{The coupling sum \texorpdfstring{$\SiN_X(d)$}{Sigma\_X(d)}}\label{ss:fd-out}

The coupling sum is expressed through two positive convergent series:
\begin{equation}\label{e:Neumann-sums}
\ES \;:=\; \sum_{n \ge 2,\ \text{even}} \frac{\tanh\!\bigl(\tfrac{\pi}{2}\sqrt{n^2-1}\bigr)}{(n^2-1)^{5/2}} \;\approx\; 0.064932,
\qquad
\ETU \;:=\; \sum_{n \ge 2} \frac{w_n}{(n^2-1)^{5/2}} \;\approx\; 0.072014,
\end{equation}
where $w_n := \tanh\bigl(\tfrac{\pi}{2}\sqrt{n^2-1}\bigr)$ for $n$ odd and $w_n := \coth\bigl(\tfrac{\pi}{2}\sqrt{n^2-1}\bigr)$ for $n$ even; both constants are positive, as sums of positive terms.

\begin{proposition}\label{p:fd-out-table}
$\SiN_R(d) = 0$, and
\[
\begin{aligned}
\SiN_S(d) \;&=\; \frac{128}{\pi^6} + \frac{256\,\ES}{\pi^5},
\qquad\qquad
\SiN_U(d) \;=\; \frac{16}{15} + \frac{16}{3\pi^2} + \frac{32}{\pi^4} + \frac{256\,\ETU}{\pi^5}, \\
\SiN_T(d) \;&=\; \frac{(d-1)\bigl(60 + 40\pi^2 + 8\pi^4\bigr)}{15\,\pi^4} - \frac{12(d-1)(d-2)}{\pi^4} + \frac{128(d-1)\,\ETU}{\pi^5}.
\end{aligned}
\]
\end{proposition}

\begin{proof}
Each $H_{k,\bfm}$ of~\eqref{e:H1-coupling-formula} is a sum of products of one-dimensional integrals, one per axis.  On an inactive axis $j$ of the deformation --- one on which neither $A_1$ nor $\rho_1$ depends and along which $A_1$ has no entry --- the factor is $\int\chi_{m_j} = \delta_{m_j 0}$, so only $\bfm$ supported on the active and modulator axes contribute; and for $k$ an inactive axis the row $(A_1)_{k,\cdot}$ vanishes, so $H_{k,\bfm} = 0$.  On the remaining axes the factors are among~\eqref{e:Neumann-1D}, whose parity conditions fix the parities of the surviving indices; the supports listed below are read off in this way, exactly as in the proof of Proposition~\ref{p:H-coeffs}.  Throughout, the gap is $|\bfm|^2 - 1$, and $\SiN_X = \pi^{-6}\sum_k\sum_{\bfm} H_{k,\bfm}^2/(|\bfm|^2 - 1)$.

For the $R$-deformation the field is affine and $\rho_1 = 0$, so $H_{k, \bfm} = (A_1)_{kk}\int_{\cube}\partial_{x_k}\psi_k\,\partial_{x_k}\Xi_{\bfm}\,d\x$.  Every axis $j \ne k$ is inactive for this entry, forcing $m_j = 0$, and the $k$th factor $\int\cos(\pi x)\,\chi_n' = -\tfrac{\pi}{\sqrt2}\delta_{n1}$ forces $\bfm = \mathbf e_k$.  So $H_{k,\bfm} = 0$ for every $\bfm$ with $|\bfm| \ne 1$ and $\SiN_R(d) = 0$.

For the $S$-deformation, $\rho_1 = 0$ and the only entries of $A_1$ are $(A_1)_{12} = (A_1)_{21} = -2$, so $H_{k, \bfm}$ vanishes unless $k \in \{1, 2\}$, and by the $1 \leftrightarrow 2$ symmetry the two values of $k$ contribute equally.  For $k = 1$, $H_{1,\bfm} = -2\sqrt2\pi\int\cos(\pi x)\chi_{m_1}\cdot\int\chi_{m_2}'$, giving the two families
\[
\renewcommand{\arraystretch}{1.3}
\begin{array}{c|c|c|c}
\bfm & \text{range} & H_{1,\bfm} & |\bfm|^2 - 1 \\ \hline
n\,\mathbf e_2 & n \text{ odd} \ge 3 & 16 & n^2 - 1 \\
a\,\mathbf e_1 + n\,\mathbf e_2 & a \text{ even} \ge 2,\ n \text{ odd} & -\tfrac{16\sqrt2}{a^2 - 1} & a^2 + n^2 - 1
\end{array}
\]
The first family sums to $256\sum_{n \ge 3\ \mathrm{odd}}(n^2 - 1)^{-1} = 64$, and the second, by the first identity of~\eqref{e:ML-identities} with $a^2 - 1$ in place of $a^2$, to
\[
\sum_{a \ge 2\ \mathrm{even}} \frac{512}{(a^2-1)^2}\sum_{n\ \mathrm{odd}}\frac{1}{n^2 + (a^2 - 1)} \;=\; 128\pi \sum_{a \ge 2\ \mathrm{even}} \frac{\tanh\bigl(\tfrac\pi2\sqrt{a^2-1}\bigr)}{(a^2-1)^{5/2}} \;=\; 128\pi\,\ES .
\]
Doubling and dividing by $\pi^6$ gives $\SiN_S$.

For the $U$-deformation, $\rho_1 = 0$, the active axes are $1, 2$ and the modulator axis is $d$, and $A_1$ has entries only in positions $(1,1)$, $(2,2)$, $(1,d)$, $(2,d)$ (Lemma~\ref{l:coeff-expansion}), so $H_{k,\bfm}$ vanishes unless $k \in \{1, 2, d\}$ and $\bfm$ is supported on an active axis $i \in \{1, 2\}$ and the modulator axis.  By~\eqref{e:Neumann-1D} the surviving families are, for $i \in \{1, 2\}$,
\[
\renewcommand{\arraystretch}{1.3}
\begin{array}{c|c|c|c|c}
k & \bfm & \text{range} & H_{k,\bfm} & |\bfm|^2 - 1 \\ \hline
i & \mathbf e_i + n\,\mathbf e_d & n \text{ odd} & \mp 4\sqrt2\bigl(1 + \tfrac{4}{n^2}\bigr) & n^2 \\
i & a\,\mathbf e_i + n\,\mathbf e_d & a \text{ odd} \ge 3,\ n \text{ odd} & \pm\tfrac{16\sqrt2}{a^2-1} & a^2 + n^2 - 1 \\
d & b\,\mathbf e_i & b \text{ even} \ge 2 & \mp 16 & b^2 - 1 \\
d & b\,\mathbf e_i + n\,\mathbf e_d & b, n \text{ even} \ge 2 & \pm\tfrac{16\sqrt2}{n^2 - 1} & b^2 + n^2 - 1
\end{array}
\]
with the upper sign for $i = 1$ (the two cluster indices $k = 1, 2$ contribute equally, as do the two halves of the $k = d$ families).  Their sums, per active axis $i$, are
\begin{align*}
\sum_{n\ \mathrm{odd}} \frac{32\,(1 + 4/n^2)^2}{n^2} &= 32\Bigl[\tfrac{\pi^2}{8} + 8\cdot\tfrac{\pi^4}{96} + 16\cdot\tfrac{\pi^6}{960}\Bigr] \;=\; 4\pi^2 + \tfrac{8}{3}\pi^4 + \tfrac{8}{15}\pi^6, \\
\sum_{a \ge 3\ \mathrm{odd}} \frac{512}{(a^2-1)^2}\sum_{n\ \mathrm{odd}}\frac{1}{n^2 + (a^2 - 1)} &= 128\pi \sum_{a \ge 3\ \mathrm{odd}} \frac{\tanh\bigl(\tfrac\pi2\sqrt{a^2-1}\bigr)}{(a^2-1)^{5/2}}, \\
\sum_{b \ge 2\ \mathrm{even}} \frac{256}{b^2 - 1} &= 128, \\
\sum_{n \ge 2\ \mathrm{even}} \frac{512}{(n^2-1)^2}\sum_{b \ge 2\ \mathrm{even}}\frac{1}{b^2 + (n^2 - 1)} &= 128\pi \sum_{n \ge 2\ \mathrm{even}} \frac{\coth\bigl(\tfrac\pi2\sqrt{n^2-1}\bigr)}{(n^2-1)^{5/2}} - 256 \sum_{n \ge 2\ \mathrm{even}}\frac{1}{(n^2-1)^3},
\end{align*}
by Lemma~\ref{l:even-power-sums} and~\eqref{e:ML-identities}.  The two hyperbolic sums are the odd- and even-index halves of $\ETU$, which is how the mixed weight $w_n$ in~\eqref{e:Neumann-sums} arises, and $128 - 256\bigl(\tfrac12 - \tfrac{3\pi^2}{64}\bigr) = 12\pi^2$.  Adding the four lines,
\begin{equation}\label{e:U-per-axis}
\text{(per active axis)} \qquad \tfrac{8}{15}\pi^6 + \tfrac{8}{3}\pi^4 + 16\pi^2 + 128\pi\,\ETU ,
\end{equation}
and doubling and dividing by $\pi^6$ gives $\SiN_U$.

For the $T$-deformation, $\rho_1 = 4(d-1)x_d \ne 0$, every axis $i < d$ is active, and $(A_1)_{ii} = 4(d-3)x_d$, $(A_1)_{dd} = 4(d-1)x_d$, $(A_1)_{id} = -4x_i$.  For $k = i < d$ and $\bfm = a\,\mathbf e_i + n\,\mathbf e_d$ ($a, n$ odd) three terms contribute: $(A_1)_{ii}$ gives $4(d-3)\sqrt2\pi\int\cos(\pi x)\chi_a'\int x\chi_n = \tfrac{8\sqrt2(d-3)}{n^2}\delta_{a1}$, the $\rho_1$ term gives $-4(d-1)\sqrt2\pi^2\int\sin(\pi x)\chi_a\int x\chi_n = -\tfrac{8\sqrt2(d-1)}{n^2}\delta_{a1}$, and $(A_1)_{di}$ gives the same as for $U$.  The $d$-dependence cancels, and the $k = i$ families coincide with those of the $U$ table; so do the $k = d$ families $\bfm = b\,\mathbf e_i$ and $b\,\mathbf e_i + n\,\mathbf e_d$, since for $b \ge 2$ the entry $(A_1)_{dd}$ and the $\rho_1$ term carry the factor $\int\chi_b = 0$.  Hence the $d-1$ active axes contribute $(d-1)$ times~\eqref{e:U-per-axis}.  In addition there are two \emph{modulator families}, supported on axis $d$ alone and present only because $\rho_1 \ne 0$: for $k = d$ and $\bfm = n\,\mathbf e_d$ with $n$ even,
\[
H_{d, n\mathbf e_d} \;=\; \begin{cases} -8\sqrt2\,(d-1), & n = 0, \\[2pt] \tfrac{16(d-1)}{n^2-1}, & n \text{ even} \ge 2, \end{cases}
\]
by~\eqref{e:Neumann-1D}.  The constant mode $\Xi_{\bzero}$ lies \emph{below} the cluster, so its gap is $|\bzero|^2 - 1 = -1$ and it contributes $-128(d-1)^2$; the family $n \ge 2$ contributes $256(d-1)^2\sum_{n\ \mathrm{even}}(n^2-1)^{-3} = 256(d-1)^2\bigl(\tfrac12 - \tfrac{3\pi^2}{64}\bigr)$, again by Lemma~\ref{l:even-power-sums}.  The two modulator families together give $-12\pi^2(d-1)^2$; adding $(d-1)$ times~\eqref{e:U-per-axis} and dividing by $\pi^6$ gives the stated $\SiN_T(d)$.
\end{proof}

\subsection{Proof of Theorem~\ref{t:fd-Hessian-d}}\label{ss:fd-proof}

\begin{proof}[Proof of Theorem~\ref{t:fd-Hessian-d}]
Substituting into the assembly formula~\eqref{e:nu-fd-assembly} the cluster parts of Lemma~\ref{l:cluster-matrices}, the coupling sums of Proposition~\ref{p:fd-out-table}, the volume derivatives $V''(0)[\bth_X, \bth_X]$ of Proposition~\ref{p:Iso-irrep}, and the norms $\|\bth_X\|^2$ of~\eqref{e:norms}, and simplifying, gives the four formulas~\eqref{e:nu-Neumann}.  The value $\nu^{f_d}_R = \tfrac{16}{d\cdot 2^d\pi^2}$ agrees with the box computation of Proposition~\ref{p:R-box}, which used no perturbation theory.

For positivity, the series terms in~\eqref{e:nu-Neumann} have positive coefficients and $\ES, \ETU > 0$, so it suffices that the polynomial part of each numerator be positive: $2560 - 16\pi^4 > 0$ since $\pi^4 < 160$; $16\pi^6 + (40d-80)\pi^4 + (600-120d)\pi^2 = 16\pi^6 + 40\pi^2\bigl[(\pi^2 - 3)d + 15 - 2\pi^2\bigr] > 0$ for $d \ge 2$, the bracket being increasing in $d$ and equal to $9$ at $d = 2$; and $16\pi^6 - 40\pi^4 + 480\pi^2 = 8\pi^2(2\pi^4 - 5\pi^2 + 60) > 0$ since $5\pi^2 < 60$.
\end{proof}

\section{Proofs of Theorem~\ref{t:Main} and Lemma~\ref{l:dish}}\label{s:Completion}

\subsection*{Proof of Theorem~\ref{t:Main}}

\begin{proof}
Let $F$ be a $B_d$-invariant functional, smooth near $\cube$ on $\calQ$.  By Proposition~\ref{p:Crit-Block} (see also Remark~\ref{r:irrep-intuition}), $F$ is critical at $\cube$, its Hessian restricted to $\TP$ is block-diagonal in symmetry-adapted coordinates, acting as the scalar $\nu^F_X$ on the block $X \in \{R, S, T, U\}$ (only three blocks at $d = 2$), and this restricted Hessian is the Hessian of $F$ on the manifold $\P$ of~\eqref{e:M-defn}, whose tangent space at $\bzero$ is $\TP$.  If the four scalars are strictly positive, Taylor's theorem on $\P$ gives $F(\Omega_{\bdelta}) \ge F(\cube) + c\,\|\bdelta\|^2$ for some $c > 0$ and all sufficiently small $\bdelta \ne \bzero$ with $\Omega_{\bdelta} \in \P$, and since every body in $\tP$ near $\cube$ is similar to a body in $\P$ (Lemma~\ref{l:Pd-dim}), this is a strict local minimum of $F$ on $\P$, and on $\tP$ modulo similarity.  The three parts of Theorem~\ref{t:Main} therefore reduce to the positivity of the four Hessian scalars for $F = I$, $F = \blam$, and $F = f_d$, respectively.

\emph{Part (a).} The volume and surface area of $\Omega_{\bdelta}$ depend smoothly on the displaced vertex positions, and $|\Omega_{\bdelta}| > 0$ near $\cube$, so $I = \tfrac{P^d}{|\Omega|^{d-1}}$ is a smooth $B_d$-invariant functional on $\P$. By Proposition~\ref{p:Iso-irrep}, its four Hessian eigenvalues at $\cube$ are given by~\eqref{e:Iso-eigs}, all strictly positive. By the discussion above, $\cube$ is a strict local minimizer of $I$ on $\P$, with $I(\cube) = (2d)^d/1^{d-1} = (2d)^d$.

\emph{Part (b).} The lowest Dirichlet eigenvalue of $\cube$ is simple, and by the pullback argument of Section~\ref{ss:Material} the map $\bdelta \mapsto \lambda_1(\Omega_{\bdelta})$ is real-analytic near $\bzero$; hence $\blam = |\Omega|^{2/d}\lambda_1$ is a smooth $B_d$-invariant functional on $\P$. By Theorem~\ref{t:Lam-Hessian-d}, its four Hessian eigenvalues at $\cube$ are given by the closed forms~\eqref{e:nu-Dirichlet}, and these are strictly positive for all $d \ge 2$. Hence $\cube$ is a strict local minimizer of $\blam$ on $\P$, with $\blam(\cube) = d\pi^2$.

\emph{Part (c).} The functional $\bmu$ itself is not differentiable at $\cube$ (the $d$-fold Neumann cluster splits under generic perturbations), so we argue through the auxiliary functional $f_d$ of~\eqref{e:fd-defn}, which is $B_d$-invariant and smooth near $\cube$ (Section~\ref{ss:Aux}). Its four Hessian eigenvalues at $\cube$ admit the closed forms of Theorem~\ref{t:fd-Hessian-d}, and all are strictly positive for every $d \ge 2$ (Section~\ref{ss:fd-proof}). Hence $f_d$ has a strict local minimum at $\cube$. Since $f_d \le \tfrac{1}{\bmu}$ on $\P$ with equality at $\cube$, Remark~\ref{r:fd-bound} converts the strict local minimum of $f_d$ into a strict local maximum of $\bmu$ at $\cube$, with $\bmu(\cube) = \pi^2$.
\end{proof}

\subsection*{Proof of Lemma~\ref{l:dish}}

\begin{proof}
\emph{Convexity and the deformation field.}  For $t \ge 0$ each set in the intersection~\eqref{e:dish-body} is a sublevel set of the convex function $\x \mapsto \sigma x_k + t\sum_{l \ne k} x_l^2$, so $\Omega_t$ is convex.  Let $\bth := \sum_k 2x_k V^{(k)}(\x)\,\mathbf e_k$ with $V^{(k)}(\x) := \tfrac{d-1}{12} - \sum_{l \ne k} x_l^2$.  On the facet $F_{a,\sigma} = \{x_a = \tfrac{\sigma}{2}\}$ the normal component of $\bth$ is
\begin{equation}\label{e:Vn-facet}
V^{(a)} \;=\; \frac{d-1}{12} \;-\; \sum_{k \ne a} x_k^2 \;=\; V_n\big|_{F_{a,\sigma}}.
\end{equation}

\emph{The isoperimetric ratio.}  Since $V^{(k)}$ does not depend on $x_k$, $\nabla \cdot \bth = \sum_k 2V^{(k)} = 2\bigl[\tfrac{d(d-1)}{12} - (d-1)|\x|^2\bigr]$, and $\int_{\cube}|\x|^2\,d\x = \tfrac{d}{12}$, so $V'(0) = 0$ by~\eqref{e:derivV,Pa} and $|\Omega_t| = 1 + O(t^2)$.  By~\eqref{e:derivV,Pc}, the integrand on $F_{a,\sigma}$ is $\nabla \cdot \bth - \partial_a \bth_a = \sum_{k \ne a} 2V^{(k)} = 2\bigl[\tfrac{(d-1)^2}{12} - (d-2)|\x|^2 - \tfrac14\bigr]$, using $x_a^2 = \tfrac14$; since $\int_{F_{a,\sigma}}|\x|^2\,dS = \tfrac14 + \tfrac{d-1}{12}$, each facet contributes $-\tfrac{d-1}{3}$.  Hence
\[
P'(0^+) \;=\; -\,\frac{2d(d-1)}{3} ,
\qquad
\frac{I'(0^+)}{I(\cube)} \;=\; d\,\frac{P'(0^+)}{2d} \;=\; -\,\frac{d(d-1)}{3} .
\]

\emph{The Dirichlet eigenvalue.}  $\lambda_1$ is simple, so Hadamard's formula along $\bth$, valid for Lipschitz domains with the normal velocity defined almost everywhere on the boundary (see~\cite[Ch.~5]{HenrotPierre2018}), gives, with $u = \prod_k \sqrt2\cos(\pi x_k)$ the normalized first eigenfunction,
\[
\lambda_1'(0^+) \;=\; -\int_{\partial\cube} |\partial_n u|^2\, V_n\, dS .
\] On $F_{a,\sigma}$, $|\partial_n u|^2 = 2\pi^2 \prod_{k \ne a} 2\cos^2(\pi x_k)$, and from $\int 2\cos^2(\pi x)\,dx = 1$ and $\int 2x^2\cos^2(\pi x)\,dx = \tfrac1{12} - \tfrac{1}{2\pi^2}$, the facet integral is $\int_{F_{a,\sigma}} |\partial_n u|^2\, V_n\, dS = 2\pi^2\bigl[\tfrac{d-1}{12} - (d-1)\bigl(\tfrac1{12} - \tfrac{1}{2\pi^2}\bigr)\bigr] = d - 1$.  Summing over the $2d$ facets, $\blam'(0^+) = \lambda_1'(0^+) = -2d(d-1)$.

\emph{The Neumann eigenvalue.}  $V_n$ is $B_d$-invariant, so the eigenvalue $\pi^2$ of $\cube$ stays $d$-fold degenerate along the flow and $\mu_1$ is differentiable at $t = 0^+$, its derivative being the average of the $d$ branches.  With $\psi_k = \sqrt2\sin(\pi x_k)$ an $L^2$-orthonormal basis of the cluster, Hadamard's formula gives $\mu_1'(0^+) = \tfrac1d \sum_{k=1}^{d} \int_{\partial\cube}\bigl(|\nabla_\tau \psi_k|^2 - \pi^2 \psi_k^2\bigr) V_n\, dS$.  On $F_{a,\sigma}$ we have $\nabla\psi_a = \sqrt2\pi\cos(\pi x_a)\mathbf{e}_a = 0$ and $\psi_a^2 = 2$, contributing $-2\pi^2$; for $k \ne a$ the gradient is tangential and the integrand is $2\pi^2(\cos^2\pi x_k - \sin^2\pi x_k) = 2\pi^2\cos(2\pi x_k)$.  Since $\cos(2\pi x_a) = -1$ on $F_{a,\sigma}$, the two cases combine into the single global density $2\pi^2 \sum_{k}\cos(2\pi x_k)$ on $\partial\cube$.  The term $k = a$ is the constant $-2\pi^2$ and integrates against~\eqref{e:Vn-facet} to zero; for $k \ne a$ only the $x_k^2$ term of~\eqref{e:Vn-facet} survives, and $\int_{-1/2}^{1/2} x^2\cos(2\pi x)\,dx = -\tfrac{1}{2\pi^2}$, so each facet contributes $d - 1$.  Summing over the $2d$ facets and dividing by $d$, $\bmu'(0^+) = \mu_1'(0^+) = 2(d-1)$.
\end{proof}

\section{Proof of Theorem~\ref{t:top-mode}}\label{s:TopMode}

We establish parts (a), (b), and (c) of Theorem~\ref{t:top-mode} in Sections~\ref{ss:top-iso}, \ref{ss:top-dir}, and~\ref{ss:top-neu}, respectively, and assemble the proof at the end of the section.  Throughout, $\Phi_t = \mathrm{id} + t\bth^{*}$ and $B := D\bth^{*}$.  The Jacobian is symmetric with zero diagonal, $B_{ij} = 2 \prod_{k \notin \{i, j\}} (2 x_k)$ for $i \ne j$ and $B_{ii} = 0$, and since $x_1 \cdots x_d$ is multilinear, hence harmonic, $\nabla \cdot \bth^{*} = \tr B = 0$.  The pullback coefficients of Lemma~\ref{l:coeff-expansion} are therefore
\begin{equation}\label{e:B-pullback}
\rho_1 = 0, \qquad \rho_2 = -\tfrac12 \tr (B^2), \qquad A_1 = -2B, \qquad A_2 = \rho_2\, I + 3 B^2 .
\end{equation}

\subsection{The isoperimetric ratio}\label{ss:top-iso}

Criticality is immediate: $V'(0) = \int_{\cube} \nabla \cdot \bth^{*} = 0$, and $P'(0) = 0$ by~\eqref{e:derivV,Pc}, since the facet integrand $\nabla \cdot \bth^{*} - \partial_a \bth^{*}_a = -B_{aa}$ vanishes identically.  Since $\tr B = 0$, formula~\eqref{e:derivV,Pb} gives $V''(0) = -\int_{\cube} \tr(B^2)$; the sum $\tr(B^2) = \sum_{i \ne j} B_{ij}^2$ has $d(d-1)$ terms, each integrating to $4 \cdot 3^{\,2-d}$ since $\int_{-1/2}^{1/2}(2x)^2\,dx = \tfrac13$, so
\begin{equation}\label{e:B-Vpp}
V''(0) \;=\; -\,4\,d(d-1)\, 3^{\,2-d} .
\end{equation}
For the perimeter we use~\eqref{e:derivV,Pd}.  On the facet $F_{a, \sigma}$ the tangential strain $A_F = 2(B_{ij})_{i, j \ne a}$ is traceless, so the integrand reduces to $\tr B_F - \tfrac12 \tr(A_F^2)$ with $\tr B_F = \sum_{k \ne a} \sum_i B_{ik}^2$ and $\tr(A_F^2) = 4\sum_{i, j \ne a} B_{ij}^2$.  Every facet is equivalent by symmetry, and $\int_{F_{a,\sigma}} B_{ij}^2 = 4 \cdot 3^{\,3-d}$ when $a \notin \{i,j\}$ and $4 \cdot 3^{\,2-d}$ when $a \in \{i, j\}$ (the fixed factor $(2x_a)^2 = 1$ drops out).  The ordered pairs $i \ne j \in \{1,\ldots,d\} \setminus \{a\}$ number $(d-1)(d-2)$, so $\int \tr(A_F^2) = 16(d-1)(d-2)\,3^{\,3-d}$; in $\tr B_F$, the pairs with $i = a$ contribute $(d-1)$ terms of $4 \cdot 3^{\,2-d}$ and the pairs with $i \ne a$ contribute $(d-1)(d-2)$ terms of $4 \cdot 3^{\,3-d}$.  Each facet integral is therefore $4(d-1)(7 - 3d)\,3^{\,2-d}$, and summing over the $2d$ facets, $P''(0) = -\,8\, d(d-1)(3d - 7)\, 3^{\,2-d}$. Since $V'(0) = P'(0) = 0$, formula~\eqref{e:Ipp-formula} applies:
\[
\frac{I''(0)}{I(\cube)} \;=\; \frac{P''(0)}{2} - (d-1)\,V''(0)
\;=\; \bigl[-4d(d-1)(3d-7) + 4d(d-1)^2\bigr] 3^{\,2-d}
\;=\; -\,8\,d(d-1)(d-3)\, 3^{\,2-d} .
\]
This proves Theorem~\ref{t:top-mode} part~(a).

\subsection{The Dirichlet eigenvalue}\label{ss:top-dir}

Criticality holds since $\lambda'(0) = 0$ by~\eqref{e:lam1-formula} and parity, every entry of $B$ being odd in each of $d - 2$ variables.  We apply the pulled-back expansion of Section~\ref{ss:Material} with the coefficients~\eqref{e:B-pullback}.  By~\eqref{e:lam2-final}, $\lambda''(0) = \PiD(d) + \SiD(d)$, where the source in $H_\bfm$ is $h = 2\nabla \cdot (B \nabla u_0)$, and $\SiD(d) \le 0$ since every gap $\pi^2(|\bfm|^2 - d)$ is positive, $u_0$ being the ground state --- no estimate on the coupling is needed.  The polynomial part $\PiD(d) = 2\langle L_2 u_0, u_0\rangle - 2\lambda_0 \langle \rho_2 u_0, u_0 \rangle$ is a sum of elementary product integrals: as in the proof of Lemma~\ref{l:poly-table}, $\langle L_2 u_0, u_0\rangle = \sum_{k,l}\int_{\cube}(A_2)_{kl}\,\partial_k u_0\,\partial_l u_0\,d\x$ and $\langle \rho_2 u_0, u_0\rangle = \int_{\cube}\rho_2\,u_0^2\,d\x$, where by~\eqref{e:B-pullback} the entries $(A_2)_{kk} = \rho_2 + 3(B^2)_{kk}$, $(A_2)_{kl} = 3(B^2)_{kl}$ and $\rho_2$ are polynomials in the $x_j$, so every integral factorizes over the axes into the one-dimensional moments $\int_{-1/2}^{1/2} (2x)^2\,dx = \tfrac13$, $\int_{-1/2}^{1/2} (2x)\,2\sin(\pi x)\cos(\pi x)\,dx = \tfrac1\pi$ and $c_2 := \int_{-1/2}^{1/2} (2x)^2\,2\cos^2(\pi x)\,dx = \tfrac{\pi^2 - 6}{3\pi^2}$.  This gives
\begin{equation}\label{e:B-lam-poly}
\PiD(d) \;=\; d(d-1)\Bigl[(d-2)\bigl(16 - \tfrac{4\pi^2}{3}\bigr) c_2^{\,d-3} \;+\; 4(d+4)\,\pi^2\, c_2^{\,d-2}\Bigr]
\;=\; 8\, d(d-1)\bigl(d + \pi^2 - 8\bigr)\, c_2^{\,d-3} ,
\end{equation}
the second equality because $c_2 = \tfrac{\pi^2 - 6}{3\pi^2}$ collapses the bracket.
Since $V'(0) = 0$, formula~\eqref{e:lam-bar-formula} gives $\blam''(0) = \lambda''(0) + 2\pi^2 V''(0)$, and with~\eqref{e:B-Vpp}, $\SiD(d) \le 0$ and $c_2^{\,d-3} = 3^{\,3-d} \eta^{\,d-3}$, $\eta := 1 - \tfrac{6}{\pi^2}$,
\begin{equation}\label{e:B-lambar-bound}
\blam''(0) \;\le\; \PiD(d) - 8\pi^2 d(d-1)\,3^{\,2-d}
\;=\; 8\,d(d-1)\, 3^{\,2-d} \Bigl[\, g(d) - \pi^2 \,\Bigr],
\qquad
g(d) := 3\bigl(d + \pi^2 - 8\bigr)\, \eta^{\,d-3} .
\end{equation}
For $d \ge 4$ the bracket is negative: $\pi^2 < 10$ gives $\eta < \tfrac25$, so $g(4) = 3(\pi^2 - 4)\eta < \tfrac65(\pi^2 - 4) < \tfrac{36}{5} < \pi^2$, and $\tfrac{g(d+1)}{g(d)} = \eta\,\tfrac{d + \pi^2 - 7}{d + \pi^2 - 8} \le \eta\,\tfrac{\pi^2 - 3}{\pi^2 - 4} < \tfrac25 \cdot \tfrac75 < 1$, so $g$ is decreasing on $d \ge 4$.  Hence $\blam''(0) < 0$ for every $d \ge 4$, with no information about $\SiD(d)$ beyond its sign.  For $d = 3$, where $\PiD(3) = 48\pi^2 - 240$, the coupling is needed, and a single mode suffices: every term of $-\SiD(3)$ is nonnegative, and for $\bfm = (2, 2, 2)$, with gap $\pi^2(12 - 3) = 9\pi^2$, the elementary product integral gives $H_{(2,2,2)} = -\tfrac{16384}{27\pi^2}$, so that $-\SiD(3) \ge \tfrac{2 H_{(2,2,2)}^2}{9\pi^2} = \tfrac{2^{29}}{3^{8}\pi^6} > 80$ by $\pi^2 < 10$.  Then $\blam''(0) = 48\pi^2 - 240 - 16\pi^2 + \SiD(3) = 32\pi^2 - 240 + \SiD(3) < 0$, since $32\pi^2 - 240 < 80$ by $\pi^2 < 10$.  This proves Theorem~\ref{t:top-mode} part~(b).

\subsection{The Neumann eigenvalue}\label{ss:top-neu}

We use the notation of Section~\ref{ss:cluster-reduction}: $\{\Xi_{\bfm}\}$ is the Neumann eigenbasis, $W$ the cluster, and $K_j, M_j, H$ the matrices of~\eqref{e:GEP}--\eqref{e:H1-coupling-formula}.  The argument has three steps.  Lemma~\ref{l:B-cluster} shows that the cluster does not split through second order and writes the common second-order coefficient as $\mu^{(2)}(d) = \PiN(d) - \SiN(d)$, with $\PiN(d)$ elementary and $\SiN(d)$ an integral over products of theta-type series; Lemma~\ref{l:B-Sigma-bound} bounds that integral in closed form; and the final step checks the resulting inequality for each $d \ge 6$.  Lemma~\ref{l:B-cluster} is phrased in terms of the kernels
\begin{equation}\label{e:theta-kernels}
\vartheta_j(t) \;:=\; \sum_{n \ \mathrm{odd}} \frac{e^{-n^2 t}}{n^{j}} \quad (j = 0, 2, 4), \qquad
\vartheta_{\mathrm e}(t) \;:=\; \frac{4}{\pi^2} + \frac{8}{\pi^2}\sum_{n \ \mathrm{even} \ \ge 2} \frac{e^{-n^2 t}}{(n^2-1)^2} .
\end{equation}

\begin{lemma}[Cluster expansion along $\bth^{*}$]\label{l:B-cluster}
Along $\bth^{*}$ one has $H|_W = 0$, and the second-order cluster matrix is scalar, $G = \mu^{(2)}(d)\, \mathrm{Id}_W$: all $d$ branches coincide through second order and $\mu_1(\Omega_t) = \pi^2 + \mu^{(2)}(d)\, t^2 + O(t^3)$.  Explicitly $\mu^{(2)}(d) = \PiN(d) - \SiN(d)$ with
\begin{align}
\PiN(d) &\;=\; 12\pi^2 (d-1)\, 3^{\,2-d} \;+\; 8 (d-1)(d-2)\, 3^{\,3-d} , \label{e:B-Pi-top} \\
\SiN(d) &\;=\; 256 (d-1) \Bigl(\frac{32}{\pi^4}\Bigr)^{d-2}  \int_0^{\infty} e^{t}\, \vartheta_{\mathrm e}(t) \Bigl[  \vartheta_0(t)\, \vartheta_4(t)^{d-2} \;+\; (d-2)\, \vartheta_2(t)^2\, \vartheta_4(t)^{d-3} \Bigr] dt . \label{e:B-Sigma-top}
\end{align}
\end{lemma}

\begin{proof}
\emph{First order.}  $M_1|_W = 0$ since $\rho_1 = 0$, and $(K_1)_{kl} = 2\pi^2 \int (A_1)_{lk} \cos(\pi x_k)\cos(\pi x_l)\,d\x$ vanishes because $(A_1)_{lk} = -2B_{lk}$ is zero for $l = k$ and, for $l \ne k$, odd in each of the $d - 2 \ge 1$ variables $x_j$, $j \notin \{k, l\}$.  Hence $H|_W = 0$, so projecting the order-$t^2$ equation of Section~\ref{ss:cluster-reduction} onto $W$ shows that the $\mu^{(2)}_i$ are the eigenvalues of the second-order cluster matrix $G := K_2|_W - \pi^2 M_2|_W - \pi^{-2}\sum_{|\bfm| \ne 1} H_{\cdot, \bfm}H_{\cdot, \bfm}^{T}/(|\bfm|^2 - 1)$.

\emph{The coupling coefficients.}  Since $\rho_1 = 0$, \eqref{e:H1-coupling-formula} gives $H_{k, \bfm} = -2\sum_{l \ne k}\int_{\cube} B_{lk}\,\partial_k\psi_k\,\partial_l\Xi_{\bfm}\,d\x$, and each term factorizes over the axes: with $\partial_k\psi_k = \sqrt2\pi\cos(\pi x_k)$ and $B_{lk} = 2\prod_{j \notin\{k,l\}}(2x_j)$, the axis-$k$ factor is $\int\cos(\pi x)\chi_{m_k} =: \alpha(m_k)$, the axis-$l$ factor is $\int \chi_{m_l}' = -2\sqrt2$ for $m_l$ odd (and $0$ otherwise), and each remaining axis contributes $\int 2x\,\chi_{m_j} =: \gamma(m_j)$, all by~\eqref{e:Neumann-1D}.  Thus $H_{k,\bfm}$ vanishes unless $m_k$ is even and every other index is odd, and then
\begin{equation}\label{e:B-H-coeff}
H_{k, \bfm} \;=\; 16\pi\, \alpha(m_k) \sum_{l \ne k}\, \prod_{j \ne k, l} \gamma(m_j)
\;=\; 16\pi\,\alpha(m_k)\prod_{j \ne k}\gamma(m_j)\,\sum_{l \ne k}\frac{1}{\gamma(m_l)},
\end{equation}
where $\alpha(0) = \tfrac{2}{\pi}$, $\alpha(n) = \tfrac{2\sqrt2}{\pi(1-n^2)}$ for $n$ even $\ge 2$, and $\gamma(n) = -\tfrac{4\sqrt2}{\pi^2 n^2}$ for $n$ odd.
Two cluster indices $k \ne k'$ require $m_k$ even and $m_k$ odd respectively, so $H_{k,\bfm}H_{k',\bfm} = 0$ for every $\bfm$: the coupling term of $G$ is diagonal.

It remains to compute $G$ itself: first its off-diagonal entries, then the two diagonal contributions.

\emph{Off-diagonal entries.}  For $k 
e l$, $(A_2)_{kl} = 3(B^2)_{kl} = 3\sum_j B_{kj}B_{jl}$ by~\eqref{e:B-pullback}, and each product $B_{kj}B_{jl}$ contains the factor $2x_k$ exactly once, so $(K_2)_{kl} = 2\pi^2\int (A_2)_{kl}\cos(\pi x_k)\cos(\pi x_l)\,d\x$ vanishes by oddness in $x_k$; likewise $(M_2)_{kl} = \int \rho_2\psi_k\psi_l\,d\x = 0$ because $\rho_2$ is even in every variable while $\psi_k\psi_l$ is odd in $x_k$.  Hence $G$ is diagonal, and by the $S_d$-symmetry of $\bth^{*}$ its diagonal entries coincide: $G = \mu^{(2)}(d)\,\mathrm{Id}_W$ with $\mu^{(2)}(d) = (K_2)_{kk} - \pi^2(M_2)_{kk} - \pi^{-2}\sum_{\bfm} H_{k,\bfm}^2/(|\bfm|^2 - 1)$ for any fixed $k$.

\emph{The cluster part.}  Here
\[
(A_2)_{kk} \;=\; \rho_2 + 3(B^2)_{kk} \;=\; -\tfrac12\sum_{i \ne j}B_{ij}^2 + 3\sum_{j \ne k}B_{kj}^2, \qquad B_{ij}^2 \;=\; 4\prod_{l \notin\{i,j\}}(2x_l)^2 .
\]
The one-dimensional moments $\int_{-1/2}^{1/2}(2x)^2\,dx = \tfrac13$ and $\int_{-1/2}^{1/2} 2\cos^2(\pi x)\,dx = \int_{-1/2}^{1/2} 2\sin^2(\pi x)\,dx = 1$, together with $c_2$ of Section~\ref{ss:top-dir} and $s_2 := \int_{-1/2}^{1/2} (2x)^2\,2\sin^2(\pi x)\,dx = \tfrac{\pi^2 + 6}{3\pi^2}$, give
\[
\int_{\cube} B_{ij}^2\,2\cos^2(\pi x_k)\,d\x \;=\; \begin{cases} 4\cdot 3^{2-d}, & k \in \{i,j\}, \\[2pt] 4c_2\,3^{3-d}, & k \notin \{i,j\}, \end{cases}
\]
and the same with $c_2$ replaced by $s_2$ for $2\sin^2(\pi x_k)$.  Counting the $2(d-1)$ ordered pairs $(i,j)$ containing $k$ and the $(d-1)(d-2)$ not containing it,
\[
\begin{gathered}
(K_2)_{kk} = \pi^2\Bigl[8(d-1)\,3^{2-d} - 2(d-1)(d-2)\,c_2\,3^{3-d}\Bigr], \\
(M_2)_{kk} = -4(d-1)\,3^{2-d} - 2(d-1)(d-2)\,s_2\,3^{3-d},
\end{gathered}
\]
and since $s_2 - c_2 = \tfrac{4}{\pi^2}$, $(K_2)_{kk} - \pi^2(M_2)_{kk} = 12\pi^2(d-1)3^{2-d} + 8(d-1)(d-2)3^{3-d} = \PiN(d)$, which is~\eqref{e:B-Pi-top}.

\emph{The coupling part.}  Write $1/(|\bfm|^2 - 1) = \int_0^{\infty} e^{-(|\bfm|^2 - 1)t}\,dt$ (the gap is at least $d - 2 > 0$) and expand the square of the last sum in~\eqref{e:B-H-coeff} into its $d - 1$ diagonal terms and $(d-1)(d-2)$ cross terms.  In a diagonal term one axis $l \ne k$ carries the weight $\gamma(m_l)^0$ and the other $d - 2$ axes carry $\gamma(m_j)^2$; in a cross term two axes carry $\gamma^1$ and the remaining $d - 3$ carry $\gamma^2$.  The sum over $\bfm$ then factorizes axis by axis into the series $\sum_{m\ \mathrm{even}}\alpha(m)^2 e^{-m^2 t} = \vartheta_{\mathrm e}(t)$ (since $\alpha(0)^2 = \tfrac{4}{\pi^2}$ and $\alpha(n)^2 = \tfrac{8}{\pi^2(n^2-1)^2}$), $\sum_{n\ \mathrm{odd}} e^{-n^2 t} = \vartheta_0(t)$, $\sum_{n\ \mathrm{odd}}\gamma(n) e^{-n^2 t} = -\tfrac{4\sqrt2}{\pi^2}\vartheta_2(t)$ and $\sum_{n\ \mathrm{odd}}\gamma(n)^2 e^{-n^2 t} = \tfrac{32}{\pi^4}\vartheta_4(t)$, and collecting the factor $(16\pi)^2/\pi^2 = 256$ gives~\eqref{e:B-Sigma-top}.
\end{proof}

The integral in~\eqref{e:B-Sigma-top} cannot be evaluated in closed form, but it can be bounded by one.  We use the Laplace transforms $\Theta_j(s) := \int_0^{\infty} e^{-st}\,\vartheta_j(t)\,dt$ of the kernels, evaluated termwise:
\begin{equation}\label{e:theta-transforms}
\begin{gathered}
\Theta_0(s) \;=\; \sum_{n \ \mathrm{odd}} \frac{1}{n^2 + s} \;=\; \frac{\pi}{4\sqrt{s}}\,\tanh\frac{\pi\sqrt{s}}{2},
\qquad
\Theta_2(s) \;=\; \sum_{n \ \mathrm{odd}} \frac{1}{n^2(n^2 + s)} \;=\; \frac{\pi^2/8 - \Theta_0(s)}{c}, \\
\Theta_{22}(s) \;:=\; \int_0^{\infty} e^{-st}\,\vartheta_2(t)^2\,dt \;=\; \sum_{m, n \ \mathrm{odd}} \frac{1}{m^2 n^2 (m^2 + n^2 + s)} ,
\end{gathered}
\end{equation}
the closed form of $\Theta_0$ being the first identity of~\eqref{e:ML-identities} with $a = \sqrt s$, and that of $\Theta_2$ following from it by partial fractions in $n^2$.  Each factor of the integrand in~\eqref{e:B-Sigma-top} is dominated by a short sum of exponentials, after which the integral reduces to finitely many values of $\Theta_0$ and $\Theta_{22}$.

\begin{lemma}[Closed-form bound for the coupling integral]\label{l:B-Sigma-bound}
Let $\beta := 96/\pi^4$, let
\[
\bar\Theta_{22}(s) \;:=\; \sum_{n \in \{1, 3, 5, 7\}} \frac{\Theta_2(n^2 + s)}{n^2} \;+\; \frac{\pi^2}{8}\,\frac{1}{81 + s}\sum_{n \ \mathrm{odd} \ \ge 9}\frac{1}{n^2}
\;\;\ge\;\; \Theta_{22}(s) ,
\]
and let
\[
(\kappa_i, r_i)_{i=1}^5 \;=\; \Bigl(\tfrac{4}{\pi^2},\, 0\Bigr),\ \Bigl(\tfrac{8}{9\pi^2},\, 4\Bigr),\ \Bigl(\tfrac{8}{225\pi^2},\, 16\Bigr),\ \Bigl(\tfrac{8}{1225\pi^2},\, 36\Bigr),\ \Bigl(\tfrac{8}{\pi^2}\Bigl[\tfrac{\pi^2 - 8}{16} - \tfrac19 - \tfrac{1}{225} - \tfrac{1}{1225}\Bigr],\, 64\Bigr),
\]
so that $\sum_i \kappa_i = \vartheta_{\mathrm e}(0) = \tfrac12$.  Then, for every $d \ge 4$,
\begin{equation}\label{e:B-Sigma-bound}
\begin{aligned}
\frac{\SiN(d)}{8(d-1)\,3^{2-d}} \le \bar\Sigma(d) :=& \; 32 \sum_{i=1}^{5} \kappa_i \Bigl[\beta^{d-2}\,\Theta_0(d - 3 + r_i) + \bigl(1 - \beta^{d-2}\bigr)\,\Theta_0(d + 5 + r_i)\Bigr] \\
&+ \frac{3072\,(d-2)}{\pi^4} \sum_{i=1}^{5} \kappa_i \Bigl[\beta^{d-3}\bar\Theta_{22}(d - 4 + r_i) + \bigl(1 - \beta^{d-3}\bigr)\bar\Theta_{22}(d + 4 + r_i)\Bigr] .
\end{aligned}
\end{equation}
\end{lemma}

\begin{proof}
All series below have positive terms, so termwise integration is justified.  Write $h(t) := \beta\, e^{t}\, \vartheta_4(t) = \beta \sum_{n \ \mathrm{odd}} n^{-4} e^{-(n^2 - 1)t}$, so that $\bigl(\tfrac{32}{\pi^4}\bigr)^k \vartheta_4^k = 3^{-k} e^{-kt} h^k$ and~\eqref{e:B-Sigma-top} becomes
\[
\frac{\SiN(d)}{8(d-1)3^{2-d}} \;=\; 32 \int_0^{\infty} e^{-(d-3)t}\, \vartheta_{\mathrm e}\, \vartheta_0\, h^{d-2}\, dt \;+\; \frac{3072\,(d-2)}{\pi^4} \int_0^{\infty} e^{-(d-4)t}\, \vartheta_{\mathrm e}\, \vartheta_2^2\, h^{d-3}\, dt .
\]
\emph{The envelope for $h^k$.}  From the series, $h/\beta = 1 + \sum_{n \ \mathrm{odd} \ge 3} n^{-4} e^{-(n^2-1)t} \le 1 + a\, e^{-8t}$ with $a := \sum_{n \ \mathrm{odd} \ge 3} n^{-4} = \tfrac{\pi^4}{96} - 1$, since $n^2 - 1 \ge 8$ there.  As $x \mapsto (1 + ax)^k$ is convex on $[0, 1]$, it lies below its chord, so with $x = e^{-8t}$,
\[
h^k \;\le\; \beta^k \bigl(1 + a e^{-8t}\bigr)^k \;\le\; \beta^k\Bigl[1 + \bigl((1+a)^k - 1\bigr) e^{-8t}\Bigr] \;=\; \beta^k + \bigl(1 - \beta^k\bigr) e^{-8t},
\]
the last equality because $\beta(1 + a) = 1$: each power $h^k$ is dominated by a \emph{convex combination} of $1$ and $e^{-8t}$, and after integration the weight $e^{-8t}$ shifts the argument of the Laplace transforms by $8$, which is the form of the brackets in~\eqref{e:B-Sigma-bound}.

\emph{The envelope for $\vartheta_{\mathrm e}$.}  Keeping the terms $n = 2, 4, 6$ of the series~\eqref{e:theta-kernels} for $\vartheta_{\mathrm e}$ and bounding the terms $n \ge 8$ by $e^{-64t}$ times their mass at $t = 0$ gives $\vartheta_{\mathrm e}(t) \le \sum_{i=1}^5 \kappa_i\, e^{-r_i t}$, with $\sum_i \kappa_i = \vartheta_{\mathrm e}(0) = \tfrac12$ by construction.

\emph{The integrals.}  Integrating termwise, $\int_0^{\infty} e^{-st}\, \vartheta_0(t)\, dt = \Theta_0(s)$ and $\int_0^{\infty} e^{-st}\, \vartheta_2(t)^2\, dt = \Theta_{22}(s)$.  Combining the two envelopes gives the first sum in~\eqref{e:B-Sigma-bound}; for the second, the partial fraction $\sum_{m \ \mathrm{odd}} \tfrac{1}{m^2(m^2 + A)} = \Theta_2(A)$ collapses the double sum to $\Theta_{22}(s) = \sum_{n \ \mathrm{odd}} n^{-2}\, \Theta_2(n^2 + s)$, and bounding the terms $n \ge 9$ via $\Theta_2(A) \le \tfrac{\pi^2}{8A}$ and $n^2 + s \ge 81 + s$ gives $\Theta_{22} \le \bar\Theta_{22}$.
\end{proof}

We can now prove Theorem~\ref{t:top-mode} part~(c).  By Lemma~\ref{l:B-cluster} the cluster does not split through second order, so $\bmu$ is differentiable along the flow with $\bmu'(0) = 0$.  Since $\bmu(\Omega_t) = |\Omega_t|^{2/d} \mu_1(\Omega_t)$ and $|\Omega_t|^{2/d} = 1 + \tfrac{1}{d}V''(0)\,t^2 + O(t^3)$, Lemma~\ref{l:B-cluster} and~\eqref{e:B-Vpp} give
\begin{equation}\label{e:B-bmu-top}
\bmu''(0) \;=\; 2\Bigl[ \PiN(d) - \SiN(d) - 4\pi^2 (d-1)\, 3^{\,2-d} \Bigr] .
\end{equation}
Since $\PiN(d) - 4\pi^2(d-1)3^{\,2-d} = 8(d-1)3^{\,2-d}\bigl[\pi^2 + 3(d-2)\bigr]$, positivity of $\bmu''(0)$ is equivalent to $\tfrac{\SiN(d)}{8(d-1)3^{\,2-d}} < \pi^2 + 3(d - 2)$, and by Lemma~\ref{l:B-Sigma-bound} it suffices to show $\bar\Sigma(d) < \pi^2 + 3(d-2)$ for $d \ge 6$.  We treat $d \ge 7$, where a crude bound suffices, and $d = 6$, where it does not, separately.

\emph{Case $d \ge 7$.}  Since $\Theta_0$ and $\bar\Theta_{22}$ are decreasing and $\sum_i \kappa_i = \tfrac12$, each bracket in~\eqref{e:B-Sigma-bound} is at most its value at the smallest argument, so
\[
\bar\Sigma(d) \;\le\; 16\,\Theta_0(d-3) \;+\; \frac{1536\,(d-2)}{\pi^4}\,\bar\Theta_{22}(d-4) .
\]
Evaluating the right side gives $23.83 < \pi^2 + 15$ for $d = 7$ and $23.41 < \pi^2 + 18$ for $d = 8$.  For $d \ge 9$ we use $\Theta_0(s) \le \tfrac{\pi}{4\sqrt s}$ and, since $\Theta_2(A) \le \tfrac{\pi^2}{8(1+A)}$, $\bar\Theta_{22}(s) \le \bigl(\tfrac{\pi^2}{8}\bigr)^2\tfrac{1}{s+2}$; hence $\bar\Sigma(d) \le \tfrac{4\pi}{\sqrt{d-3}} + 24 \le \tfrac{4\pi}{\sqrt6} + 24 < 30 < \pi^2 + 3(d-2)$.

\emph{Case $d = 6$.}  Here the threshold is thin: numerically, $\SiN(6)/(8 \cdot 5 \cdot 3^{-4}) \approx 21.746$ against $\pi^2 + 12 \approx 21.870$, i.e.\ $\bmu''(0) \approx +0.122$, so the bound of Lemma~\ref{l:B-Sigma-bound} must be used in full.  By~\eqref{e:B-Sigma-bound} and~\eqref{e:theta-transforms}, $\bar\Sigma(6)$ is an explicit finite expression in $\pi$, square roots and hyperbolic tangents; evaluating it gives
\[
\bar\Sigma(6) \;=\; 21.8227\ldots \;<\; 21.8696\ldots \;=\; \pi^2 + 12 .
\]
This proves Theorem~\ref{t:top-mode} part~(c).

\section*{Acknowledgments}
The authors used Claude (Anthropic) to assist with editing the text, verifying references, and checking derivations. All AI-assisted text and derivations were checked by the authors, who take full responsibility for the content of the paper.

\printbibliography

\end{document}